\documentclass[acmlarge]{acmart}
\AtBeginDocument{%
  }

\acmJournal{POMACS}
\acmVolume{0}
\acmNumber{0}
\acmArticle{0}
\acmMonth{0}

\def \ne{\textrm{ne}}
\usepackage{dsfont}

\newcommand{\zeroindic}{{\Large{\textbf{0}}}}
\newcommand{\oneindic}{ {\Large{\textbf{1}}}}
\newcommand{\idoper}{ {{\textbf{Id}}}}

\def \ne{\textrm{ne}}
\newtheorem{theorem}{\bf Theorem}
\newtheorem{corollary}{\bf Corollary}
\newtheorem{lemmon}{\bf Lemma}
\newtheorem{proposition}{\bf Proposition}
\newtheorem{define}{\bf Definition}
\newtheorem{notation}{\bf Notation}

\newtheorem{assumption}{\bf Assumptions}

\newtheorem{convention}{\bf Convention}
\newtheorem{remark}{\bf Remark}
\usepackage{appendix}

\begin{document}

%
%

\title{A polynomial time algebraic solution to exact marginal inference in Markov Random Field models}

\author{Ikhlef Bechar}
\affiliation{%
  \institution{ikhlef.bechar@gmail.com}
  \streetaddress{}
  \city{Nice}
  \state{Nice}
  \postcode{83957}
  \country{France}}
\email{ikhlef.bechar@gmail.com}

%
\renewcommand{\shortauthors}{Ikhlef Bechar}

{\color{red} {\bf{N.B.}} \it{ This is a research draft still awaiting a rigorous peer review and attempting to solve a 
still open problem in Mathematics and Computer Science. Hence, I would be grateful if you would report to me by email your 
remarks or objections about the proof developed in this paper.}}

%
\begin{abstract}
This paper develops on algebraic grounds a polynomial time exact linear solution to the hard combinatorial problem of
marginal inference in Markov random field (MRF) models under general assumptions. To prove our claim, we first implicitly
remodel a MRF joint distribution as the unique solution of some linear identity assuming its clique potential functions
(equivalently, its individual conditional distributions) to be specified. Then, by assuming an arbitrary point subset,
we relax accordingly such a (global) linear identity for deriving a second linear identity, solely, acting on a polynomial
time number of entries (e.g.; local marginals or Fourier frequencies) of a solution. Then, we show, only using linear
algebraic techniques, that such an identity enables to capture all the entries necessary for the exact reconstruction
of an MRF marginal distribution, thus, allowing to solve for the latter, exactly and in polynomial time, using a standard
linear solver. Last, but not least, this paper also solves, once and for all, the {\bf\textit{P = NP}} conjecture.
\end{abstract}


%
%

\ccsdesc[500]{Mathematics of computing~Combinatorial optimization}
\ccsdesc[500]{Mathematics of computing~Markov networks}
%
%

%
\keywords{Markovian graph, MRF model, marginal distribution, marginal inference computation, linear Algebra,
linear solution, complexity, NP-hardness, polynomial time algorithm, {\bf\textit{P vs NP}} hypothesis.}

%

%
\maketitle

\section{Introduction}
\label{sec:intro}
Markov random fields (MRFs) \cite{Dobrushin68,Dobrushin68b,Lanford69,Spitzer70,Spitzer71,Hammersley71,Besag74} are a powerful modern probabilistic modeling tool both in theoretical
and applied sciences \cite{Lauritzen91,Bishop06,Wainwrightetal08,Koller_Friedman_09}. Conceptually, MRFs are multivariate random variables which are a
natural generalization of Markov chains to multiple dimensions (e.g.; space). Thus, the dependency structure of a MRF is represented
by an undirected graph instead of a directed $1D$ chain, and the causality property which is intrinsic of Markov chains is substituted in a
MRF with the local Markov property. The latter implies the concept of neighborhood, and suggests that each random variable is conditionally
independent of all the other random variables which are not in its neighborhood given only the ones which are in its neighborhood. An intermediary
class of graphical models is the Bayesian network class \cite{Pearl82,Koller_Friedman_09}. Like the MRFs, Bayesian networks are also multidimensional,
and graph-representable (but, rather by a directed acyclic graph), except that, like Markov chains, they enjoy the causality property \cite{Pearl09}.
The latter implies the notion of parenthood, and suggests that each random variable is conditionally independent of its non-descendants given its immediate
parents. It is then worthwhile noting that, though mathematically, any Bayesian distribution is also a MRF distribution, the causality assumption has
dramatic consequences on the way Bayesian networks are learned, simulated, and used for inference \cite{Pearl88,Besag95,Koller_Friedman_09,Pearl09}.
Naturally, all the developed results in this paper regarding the general MRF framework apply effortlessly to the general Bayesian framework too.


Historically, MRFs have been originated independently in statistical physics and in probability theory. In statistical physics  \cite{Gibbs02}
where they have been mostly known as Gibbsian fields, they had been originally developed for describing thermodynamic (macroscopic) equilibrium
states of large systems of interacting tiny particles (e.g.; electrons, atoms, molecules, etc)
\cite{Ising25,Potts52,Kinderman80,Friedli17}.
Their surge in probability theory and computer science has been, however, more progressive \cite{Spitzer75,Gemanetal84}. Originally, Markov developed in the early $1950$s the
notion of conditional independence, and suggested how some complex probabilistic models can sometimes be explained in terms of a chain of
simple dependencies using the chain rule of probability theory, and the notion of conditional independence \cite{Markov54}.
However, it was not until the late $1960$s, with the surge of computer science, that the notion of MRF started to take shape.
Hence, in $1968$, Dobrushin introduced the concepts of infinite Gibbsian state and of a Markov random field \cite{Dobrushin68,Dobrushin68b}.
In $1970$, Spitzer suggested that the random field of interacting Markovian processes on a cubic lattice, assuming binary interactions, %
is Gibbsian \cite{Spitzer70}. A year later, Hammersley and Clifford proved Gibbs-Markov equivalence (or, the MRF consistency theorem)
under the positivity assumption about a distribution suggesting that any positive Markovian distribution is Gibbsian, and reciprocally
\cite{Hammersley71}. However, Hammersley and Clifford never published their theorem thinking that the positivity assumption is too restrictive from a physical standpoint.
But later in $1974$, Besag revamped the proof of their theorem, moreover, proved that the positivity assumption about a joint distribution is not simply sufficient, but,
it is also necessary. He also suggested that a MRF distribution is uniquely determined by its individual conditional distributions  \cite{Besag74}.
In a way, this has been a noticeable step forward as it makes it possible to completely and uniquely specify any joint MRF distribution
simply by means of local (clique) interactions between random variables. 
In $1984$, Geman and Geman popularized MRFs in computer science with their influential paper \cite{Gemanetal84}.
Since then, MRFs have reached out to disciplines as diverse as  statistics \cite{Clifford90,Besag95,Jaakkola2000,Pearl09},
information theory \cite{Kschischangetal01,Yedidiaetal05,Wainwrightetal05}, artificial intelligence \cite{Pearl82,Pearl88,Cooper90,Chickering96},
computer vision and machine learning \cite{Baddeleyetal93,Winkler95,Bishop06,Li95,Wainwrightetal08},  
biostatistics \cite{WeiEtAl07,candes_et_al2017}, social networks \cite{Robins07}, and so forth.

However, a fundamental problem which has tapped, so far, the minds of many researchers concerns inference computation in MRF models,
whether it is Maximum-A-Posteriori (MAP) inference where the goal is computation of the most likely assignment, marginal inference where the
goal is computation of MRF marginal distributions over arbitrary small subsets of random variables, evaluation of the normalizing
constant of a MRF distribution (also referred, using statistical physics terminology, as the partition function), hence, pointwise
evaluation of a MRF distribution too, and so forth. It is then well-known that the problem of evaluation of the partition function
is polynomial time reducible to the marginal inference problem in a MRF model. Hence, a polynomial-time procedure to solve the former
can be recycled for solving the latter in polynomial-time. The converse of the latter statement also turns out to be true \cite{ChandrasekaranEtAl08}.

Nonetheless, it is a well established fact that either of the aforementioned inference problems is, in general, a hard
combinatorial problem \cite{Karp72,Cooper90,DagumEtAl93,Shimony94,RothEt96,Chickering96,Cipra00,Megretski00,Boykovetal01,Cohenetal06,GoldbergEtAl12}.
More precisely, MAP inference in a general MRF model is NP hard since, for instance, the maximum independent
set problem in an indirected graph (i.e.; the largest set of vertices in an indirected graph, such that, there
exists no edge connecting any couple of its vertices) which is NP-hard  \cite{Karp72,GareyEtAl79} can be reduced
in polynomial time to a MAP inference problem in a proper MRF model with regard to the same undirected graph. Marginal
inference in general MRF models is $\#$-hard (Sharp-hard) since, for instance, the independent set counting problem in an
indirected graph which is $\#$-hard  can be reduced in polynomial time to a marginal inference problem in a MRF model. One notes that
hardness reductions to CPS (Constraint Satisfaction Problem) can be equally obtained \cite{ChandrasekaranEtAl08}.
The bottom line is that, unless {\it {P = NP}}, there exists literally no polynomial time algorithm which can solve either exactly, or at least,
approximately a general instance of a marginal, or a MAP inference in a general MRF model. Nevertheless, there exists MRF instances (e.g.;
MRFs over Associative graphs, MRFs with regular priors, etc) \cite{Greigetal89,Boykovetal98,ishikawa03,Boykovetal01,Schlesingeretal07,Osokinetal11}
where MAP inference is easy, but, marginal inference is hard. This explains why it is frequently stated in the
literature that marginal inference is a harder problem than MAP inference in general MRF models.
To cut a long story short, inference computation in a general MRF model has only been dealt heuristically
so far. Since this paper is primarily concerned with the problem of marginal inference in MRF models, we refer the reader
who is interested  in MAP inference to the literature dealing with the topic.

Basically, heuristical marginal inference approaches can be divided into two main categories, namely, stochastic
sampling (MCMC) based \cite{Hastings70,Green95,Gelfand90,Winkler95}, and variational inference \cite{Jaakkola2000,Yedidiaetal05,Wainwrightetal08}
algorithms. The latter category of approaches has, since recently, gained considerable interest among the research community mainly
due to availability of fast mathematical programming solvers \cite{Nesterov09,Lesaja09,Becketal09}.
Their general principle is to replace a MRF model instance (in which inference computation is seemingly hard)
with a surrogate model which is closest to it in some sense among a parametric family of models in which exact
inference computation is easy (e.g.; a MRF model with regard to a tree), then perform inference directly with
regard to the latter as a way to circumvent the situation. This is then achieved by optimizing an objective
with respect to the parameters of a surrogate distribution, typically, a Bethe objective basically consisting of
a Kullback-Leibler divergence term plus a proper penalty term. Furthermore, when optimization of the objective is performed
using the gradient descent algorithm, solution updates are made in a way which is reminiscent of transmission of
messages between nodes, hence also the name of message passing algorithms, or Loopy Belief Propagation (LBP) (said otherwise,
Belief Propagation (BP) on graphs containing cycles) \cite{Pearl88,Yedidiaetal03,Yedidiaetal05}. Variants of the latter include the tree-reweighted
belief propagation (TRBP) algorithm \cite{Wainwrightetal05}, the convergent TRBP \cite{Kolmogorovetal06}, convergent
message-passing based on geometric programming \cite{Globersonetal07}, the dual decomposition algorithm \cite{Sontagetal10,KomodakisEtAl2011},  
to name a few. Although, except in a few particular cases (e.g.; tree-structured graphical models \cite{Pearl82}, factor graphs \cite{Kschischangetal01},
and so forth), there is generally no guarantee for a good solution, or even of convergence \cite{Sontagetal14,MurphyEtAL99}, it has been reported on
multiple occasions that variational inference techniques can sometimes perform well empirically for applications ranging from information
theory, to artificial intelligence, signal processing, to name a few \cite{Baddeleyetal93,Winkler95,Wainwrightetal05,Kolmogorovetal06,Globersonetal07,KomodakisEtAl2011,MurphyEtAL99}.

In this paper, we attempt to take a step forward with regard to the existing literature, and we propose to solve the marginal
inference problem in MRF models, exactly and in polynomial time, under rather general assumptions. More precisely, in order to
work around the unwieldy combinatorics arising from the summation of a Gibbsian model over exponentially many assignments,
as a first step, we reformulate a MRF marginal inference problem on purely linear algebraic grounds, more specifically,
as the one consisting in computing the marginal over a given point subset of a function (i.e.; a joint MRF distribution)
encoded, implicitly, as the unique solution of some linear identity. As a second step, we relax accordingly such a seemingly
hard linear identity to an easy analog which, in contrast, admits an infinity of solutions, but which all share a
common marginal over such an assumed point subset coinciding with the sought MRF marginal distribution.
We then want to underline that this paper is self-contained, moreover, its content is accessible to a wide audience,
in that: a) no background on marginal inference is, actually, required for catching up with the content of this paper;
b) all our findings throughout are proved in detail, moreover, no external results are recycled throughout (except
obviously well-known algebraic results); c) last, but, not least, all our proofs in this paper use basic notions of
linear algebra and mathematical programming.


The remainder of this paper is organized as follows. Section \ref{sec:problem_statement} states the
MRF marginal inference problem in general terms. Section \ref{sec:defs} fixes some of the definitions
and notations of the paper. Section \ref{sec:summary} succinctly and simply presents the main result
in this paper. Sections \ref{sec:implicit_mrf_dist} and \ref{sec:relaxation} develop the polynomial
time algebraic reformulation of a general MRF marginal inference problem. Section \ref{sec:complexity_algorithm}
discusses complexity and algorithmic resolution of MRF marginal inference computation. A discussion concludes this
paper.

\section{Problem statement}
\label{sec:problem_statement}
Suppose $\mathcal{G} = (\mathcal{V}, \mathcal{E})$ is  a Markovian graph with $\mathcal{V}$ and $\mathcal{E}$ standing for its vertex- (or, node-)set and
edge-set, respectively, and denote $n = |\mathcal{V}|$. Furthermore, $\forall i \in \mathcal{V}$, denote ${\ne}_{i} = \big\{ j\in \mathcal{V},  \text{s.t.},
(i,j) \in \mathcal{E} \big\}$ standing for the neighborhood of node $i$ in Markovian graph $\mathcal{G}$.  Throughout, call a hypervertex of
$\mathcal{G}$ any collection of vertices of $\mathcal{G}$, and $\forall x \in \mathcal{L}^{n}$, denote $x_c = (x_i)_{i\in c}$.
Furthermore, suppose a label-set $\mathcal{L}$ which can be assumed, without loss of generality, to be the integer subset $\big\{0,\dots,L-1\big\}$, with
$L$ standing for some integer constant which is greater than, or equal to $2$.

Now, suppose $\mathcal{X} = \big(\mathcal{X}_i\big)_{i \in \mathcal{V}}$ is a Markov random field (MRF) with regard to $\mathcal{G}$ and taking values in $\mathcal{L}^n$ which means that joint distribution of
random vector $\mathcal{X}$ satisfies following both conditions simultaneously:
\begin{enumerate}
\item Positivity:
\begin{equation}
\begin{array}{l}
\mathbb{P}\big(\mathcal{X} = x \big) > 0, \forall x \in \mathcal{L}^n
\end{array}
\nonumber 
\end{equation}
\item Local Markov property (Markovianity):
\begin{equation}
\begin{array}{l}
\mathbb{P}\big(\mathcal{X}_i = x_i / \mathcal{X}_{\neq i} = x_{\neq i} \big) = \mathbb{P} \big(\mathcal{X}_i = x_i/\mathcal{X}_{{\ne}_{i}} = x_{{\ne}_{i}}\big), \, \forall x \in \mathcal{L}^n, \,
\forall i \in \mathcal{V}
\end{array}
\nonumber  
\end{equation} where it has been assumed that $x_{\neq i} = (x_j)_{j\in\mathcal{V}/\{i\}}$, $\forall x \in \mathcal{L}^n$, and $\forall i \in \mathcal{V}$.
\end{enumerate} Alternatively, first, recall that a clique of $\mathcal{G}$ stands for any collection of its mutually neighboring nodes, and
a maximal clique of $\mathcal{G}$ stands for any clique $c$ of $\mathcal{G}$ which, moreover, is not a proper subset of any other clique of $\mathcal{G}$, in other
words, $c$ is a clique of $\mathcal{G}$, and $\forall j \in \mathcal{V}/c$, $c \cup \{j\}$ is not a clique of $\mathcal{G}$. Furthermore, denote by
$\mathscr{C}$ the set of maximal cliques of $\mathcal{G}$. Then, a completely equivalent definition of MRF joint distribution $\mathbb{P}\big(\mathcal{X} = x \big)$
\cite{Hammersley71} is that there exists strictly positive local functions $\phi_c: \mathcal{L}^{|c|} \longrightarrow \mathbb{R}$, $\forall c \in \mathcal{C}$,
(commonly referred to as clique potentials), such that, one can write
\begin{equation}
\begin{array}{l}
\mathbb{P}\big(\mathcal{X} = x \big) = \frac{1}{Z} \Pi_{c \in \mathscr{C}} \, \phi_c(x_c), \, \forall x \in \mathcal{L}^{n}
\end{array}
\label{mrf_model}
\end{equation}
where $Z$ stands for a normalizing constant (also, called the partition function) which is given by
\begin{equation}
\begin{array}{l}
Z = \frac{1}{\sum_{x \in \mathcal{L}^{n}} \big( \Pi_{c \in \mathscr{C}} \, \phi_c(x_c) \big)}
\end{array}
\label{parition_function}
\end{equation}

Throughout, we assume the clique potential functions $\phi_c: \mathcal{L}^{|c|} \longrightarrow \mathbb{R}$, $\forall c \in \mathscr{C}$
to be specified, i.e.; fixed once and for all throughout. Then, given a hypervertex ``$s$" of $\mathcal{G}$ which we assume to be arbitrary,
but, fixed once and for all throughout, the problem we propose to solve, exactly and in polynomial time, in the remainder of this paper
can be stated as follows:
\begin{equation}
\begin{array}{l}
\text{Compute: } \mathbb{P}\big(\mathcal{X}_{s} = x_{s} \big), \, \forall x_{s} \in \mathcal{L}^{|{s}|} \\
\text{s.t., }  \mathbb{P}\big(\mathcal{X} = x \big) = \frac{1}{Z} \Pi_{c \in \mathscr{C}} \, \phi_c(x_c), \, \forall x \in \mathcal{L}^{n}
\end{array}
\label{marginal_inference_mrf}
\end{equation}
Problem (\ref{marginal_inference_mrf}) is, generally, referred to as the marginal inference problem in MRF models. For instance,
one of its most common instances in the literature concerns the evaluation of the unary MRF marginal distributions of the form
$\mathbb{P}\big(\mathcal{X}_{i} = \cdot \big)$, $\forall i \in \mathcal{V}$, in which case, $s$ stands for an arbitrary node
of $\mathcal{G}$. Other instances of the MRF marginal inference problems concern reformulation of the problem of evaluation
of the partition function (\ref{parition_function}) in terms of a series of MRF marginal problems, conditional inference, and
so forth \cite{Lauritzen91,Wainwrightetal08,Koller_Friedman_09}.

For computability, we shall make throughout the following mild assumption.
\begin{assumption}
\label{assump_computationability}
\begin{enumerate}
\item The number $\sup_{i \in \mathcal{V}} \big\{ \, \big| \{i\} \cup {\ne}_{i}\big| \, \big\}$ is assumed to be either a constant which is independent of $n$,
or at most dominated by some logarithmic
function of $n$, that is, one can find two nonnegative constants $\tau_1'$ and $\tau_2'$, such that, one has
$\sup_{i \in \mathcal{V}} \big\{ \, \big|{ \{i\} \cup \ne}_{i}\big| \, \big\} \leq \tau_1' \,\log(n)+\tau_1'$, and denote
in the remainder $\tau_1 = \sup\{\tau_1', \tau_1''\}$, and $\tau_2 = \sup\{\tau_2', \tau_2''\}$.
\item $|s|$ is either assumed to be a constant, or at most
dominated by some logarithmic function of $n$; which means that one can find two nonnegative constants $\tau_1''$, and $\tau_2''$, such that,
one has $\sup_{i \in \mathcal{V}} \big\{ \, \big|\{i\} \cup {\ne}_{i}\big| \, \big\} \leq \tau_1'' \, \log(n)+\tau_2''$.
\item The number of MRF labels $L$ is assumed to be a constant, i.e; $L$ is independent of $n$.
\end{enumerate}
\end{assumption}
It is worthwhile mentioning that while assumption (1) above merely ensures that the neighborhood size of any node of $\mathcal{G}$ does not explode
as $n$ gets bigger which, also, just turns out to be the Markovianity assumption, assumption (2) above is necessary for any efficient resolution of
MRF marginal distribution $\mathbb{P}\big(\mathcal{X}_{s} = \cdot \big)$ as, at least, the latter has to be evaluated at all its $L^{|s|}$ entries,
finally, assumption (3) above regarding the constant number of labels of a MRF, it is an integral part of the latter.

For the sake of example, suppose the problem of evaluation of unary distributions of a binary Ising model on a regular $2D$ lattice endowed with a
$4$-neighborhood system (that is, its interior sites have four neighbors each, its boundary sites except its four corner sites have three neighbors each, and its four
corner sites have two neighbors each). It follows that $\mathcal{L} = \big\{0,1\big\}$, $\sup_{i \in \mathcal{V}} \big\{ \, \big| \{i\} \cup {\ne}_{i}\big| \, \big\}=5$
which as constant, and $|s|$ stands for any site of the underlying lattice.

With that being said, we propose to solve MRF marginal inference problem (\ref{marginal_inference_mrf}) under such rather general assumptions. This means
that no further restrictions are made about the number of labels $L$, the structure of Markovian graph $\mathcal{G}$, the clique potentials
$\phi_c(x_c)$, $\forall c \in \mathscr{C}$, $\forall x_{c} \in \mathcal{L}^{|{c}|}$, as well as the structure of hypervertex $s$.


Before enunciating the main result in this paper, begin by introducing the following two definitions.
\begin{define}[Augmented neighborhood]
We call the augmented neighborhood of any vertex $i \in \mathcal{V}$ the hypervertex ${\ne}_{i}^+$ defined as
${\ne}_{i}^+ = \{i\} \cup {\ne}_{i}$, with ${\ne}_{i}$ standing for the neighborhood of $i$ in $\mathcal{G}$.
\end{define}

\begin{define}[Individual conditional distributions]
We call the individual conditional distributions of MRF $\mathcal{X}$ the following (local) functions:
\begin{equation}
\begin{array}{lll}
\pi_{i}(x_{{\ne}_{i}^+}) &=& \mathbb{P}\big(\mathcal{X}_i = x_i / \mathcal{X}_{{\ne}_{i}} = x_{{\ne}_{i}}\big), \, \forall x_{{\ne}_{i}^+}
\in \mathcal{L}^{|{\ne}_{i}^+|},\, \forall i \in \mathcal{V}
\end{array}
\nonumber 
\end{equation}
Then, from formula (\ref{mrf_model}), one has
\begin{equation}
\begin{array}{lll}
\pi_{i}(x_{{\ne}_{i}^+}) &=& \frac{\Pi_{c \in \mathcal{C}, \text{s.t.}, i \in c} \phi_c(x_c)}{ \sum_{x_i \in \mathcal{L}} \big(\Pi_{c \in \mathcal{C}, \text{s.t.}, i \in c} \phi_c(x_c) \big)}, \, \forall x_{{\ne}_{i}^+} \in \mathcal{L}^{|{\ne}_{i}^+|}, \, \forall i \in \mathcal{V}
\end{array}
\label{neighb_conditional}
\end{equation}
\end{define}
One observes that $\forall i \in \mathcal{V}$, one has
\begin{enumerate}
\item first, by assumption, $\pi_{i} (x_{{\ne}_{i}^+})>0$, $\forall x_{{\ne}_{i}^+} \in \mathcal{L}^{|{\ne}_{i}^+|}$;
\item second, by construction,  $\sum_{x_i \in \mathcal{L}} \pi_{i}(x_{{\ne}_{i}^+}) = 1, \, \forall x_{{\ne}_{i}} \in \mathcal{L}^{|{\ne}_{i}|}$.
\end{enumerate}

\section{Summary of the main result}
\label{sec:summary}
The main result in this paper regarding exact polynomial time computation of MRF marginal distribution $\mathbb{P}\big(\mathcal{X}_s = \cdot \big)$ can be
stated as follows. Assume "$s$" to be an arbitrary hypervertex of $\mathcal{G}$, and suppose a set of $n$ local functions indexed by hypervertices of
graph $\mathcal{G}$ as follows:
\begin{equation}
\begin{array}{lll}
\Big\{ p_{s\cup {\ne}_{i}^+}: \mathcal{L}^{|s\cup {\ne}_{i}^+|} \rightarrow \mathbb{R}, \, \forall i \in \mathcal{V} \Big \}
\end{array}
\nonumber
\end{equation}
satisfying the following three types of linear constraints:
\begin{enumerate}
\item first, the following ``restricted" local Markov property:
\begin{equation}
\begin{array}{lll}
p_{s\cup {\ne}_{i}^+}(x_{s\cup {\ne}_{i}^+}) - {\pi}_{i}(x_{{\ne}_{i}^+}) \, \sum_{x_i\in \mathcal{L}} \, p_{s\cup {\ne}_{i}^+}(x_{s\cup {\ne}_{i}^+}) = 0, \,
\forall x_{s\cup {\ne}_{i}^+} \in \mathcal{L}^{|s\cup {\ne}_{i}^+|}, \, \forall i \in \mathcal{V}
\end{array}
\nonumber
\end{equation}
\item second, the following local consistency constraint:
\begin{equation}
\begin{array}{lll}
\sum_{k \in {\ne}_{i}^+/{\ne}_{j}^+} \, \sum_{x_k \in \mathcal{L}} p_{s\cup {\ne}_{i}^+}(x_{s\cup {\ne}_{i}^+}) =
\sum_{k \in {\ne}_{j}^+/{\ne}_{i}^+} \, \sum_{x_k \in \mathcal{L}} p_{s\cup {\ne}_{j}^+}(x_{s\cup {\ne}_{j}^+}), \, \forall i, j \in \mathcal{V}, \text{s.t., i <j }
\end{array}
\nonumber
\end{equation}
\item third, and finally, the sum to unity constraint:
\begin{equation}
\begin{array}{lll}
\sum_{x_{s\cup {\ne}_{i}^+} \in \mathcal{L}^{|s\cup {\ne}_{i}^+|}} p_{s\cup {\ne}_{i}^+}(x_{s\cup {\ne}_{i}^+}) = 1, \, \forall i \in \mathcal{V}
\end{array}
\nonumber
\end{equation}
\end{enumerate}
also, satisfy that $\forall i \in \mathcal{V}$, one has
\begin{equation}
\begin{array}{lll}
\mathbb{P}\big(\mathcal{X}_s = x_s \big) = \sum_{k \in {\ne}_{i}^+/s} \, \sum_{x_k \in \mathcal{L}} p_{s\cup {\ne}_{i}^+}(x_{s\cup {\ne}_{i}^+}), \,
\forall x_{s} \in \mathcal{L}^{|s|}
\end{array}
\label{ps_summary}
\end{equation}
Then, it is not hard to show that, under such rather minimalistic Assumptions \ref{assump_computationability}, the three types of linear constraints above give rise to an ordinary system of
linear equations of polynomial time size, thus enabling to solve, both exactly and efficiently, for $\mathbb{P}\big(\mathcal{X}_s = \cdot \big)$ using formula (\ref{ps_summary})
above. Therefore, the remainder of this paper is devoted to the useful algebraic reformulation of MRF marginal inference problem (\ref{marginal_inference_mrf})
in a way which enables to prove the claimed result above, and to develop its efficient numerical resolution accordingly.

\section{Definitions and notations}
\label{sec:defs}
In this section, we introduce a few definitions and notations which will help, clearly, describe the main findings in this paper.

\begin{notation}
Denote throughout
\begin{equation}
\begin{array}{lll}
\widehat{\pi}_{i}(x_{{\ne}_{i}^+}) &=& \frac{\pi_{i}(x_{{\ne}_{i}^+})}{ \sqrt{\sum_{x_i \in \mathcal{L}} \pi_{i}^2(x_{{\ne}_{i}^+})}}, \, \forall x_{{\ne}_{i}^+} \in \mathcal{L}^{|{\ne}_{i}^+|}, \, \forall i \in \mathcal{V}
\end{array}
\label{neighb_conditional_hat}
\end{equation}
One then easily checks that one has $\forall i \in \mathcal{V}$
\begin{enumerate}
\item $\widehat{\pi}_i(x_{{\ne}_{i}^+})>0$, $\forall x_{{\ne}_{i}^+} \in \mathcal{L}^{|{\ne}_{i}^+|}$;
\item $\sum_{x_i \in \mathcal{L}}\widehat{\pi}_{i}^2(x_{{\ne}_{i}^+}) = 1, \, \forall x_{{\ne}_{i}} \in \mathcal{L}^{|{\ne}_{i}|}$.
\end{enumerate}
\end{notation}

\begin{notation}
Denote throughout $\mathscr{F} = \big\{ f:\mathcal{L}^n\longrightarrow \mathbb{R} \big\}$.
\end{notation}

\begin{define}[Expectation operator]
Whatever a function $f \in \mathscr{F}$, we denote by $\mathbb{E}\big[f(\mathcal{X})\big]$ the expectation of random function $f(\mathcal{X})$,
that is, $\mathbb{E}\big[f(\mathcal{X})\big] = \sum_{x \in \mathcal{L}^n} \mathbb{P}\big(\mathcal{X} =x \big) \, f(x)$.
\end{define}

\begin{define}[Identity operator]
Introduce the identity operator in $\mathscr{F}$ denoted by $\idoper$ and defined $\forall u \in \mathscr{F}$ as $\idoper \,(u) \equiv u$.
\end{define}

\begin{convention}
For notational convenience throughout, whatever a linear operator $\mathscr{A}:  \mathscr{F} \rightarrow \mathscr{F}$,
and whatever a function $u \in \mathscr{F}$, we shall simply write $\mathscr{A} \, u $ instead of $\mathscr{A}(u)$,
\end{convention}

\begin{convention}
For the sake of clarity throughout, we shall use the mathematical symbol ``='' both for affectations and identities between scalars and linear operators,
and the mathematical symbol ``$\equiv$'' both for affectations and identities between functions in $\mathscr{F}$.
\end{convention}

\begin{define}[Inner product space of discrete functions in $\mathscr{F}$]
Define the addition between any two elements of $\mathscr{F}$ as
\begin{equation}
\begin{array}{lll}
\forall f,g \in \mathscr{F}, \, (f + g)(x) = f(x) + g(x), \, \forall x \in \mathcal{L}^{n}
\end{array}
\nonumber
\end{equation}
and the multiplication of any element of  $\mathscr{F}$ by a real scalar as
\begin{equation}
\begin{array}{lll}
\forall f \in \mathscr{F}, \forall \alpha \in \mathbb{R}, \, (\alpha \cdot g)(x) = \alpha\, g(x), \, \forall x \in \mathcal{L}^{n}
\end{array}
\nonumber
\end{equation}
in such a way that, $\mathscr{F}$, when endowed with the usual addition and multiplication, defines a vector space over $\mathbb{R}$.

Furthermore, endow $\mathscr{F}$ with the standard (i.e.; the Euclidean) inner product $\big<\cdot,\cdot\rangle $ defined as
\begin{equation}
\begin{array}{lll}
\forall f,g \in \mathscr{F}, \, \big<f,g\rangle  = \sum_{x\in \mathcal{L}^n} f(x)\, g(x)
\end{array}
\nonumber
\end{equation}
and denote by $\|\cdot\|$ the Euclidean norm defined as
\begin{equation}
\begin{array}{lll}
\forall f \in \mathscr{F},  \, \|f\| = \sqrt{\big<f,f\rangle } = \sqrt{ \sum_{x\in \mathcal{L}^n} f^2(x)}
\end{array}
\nonumber
\end{equation}
in such a way that,  $\big(\mathscr{F}, \big<\cdot,\cdot\rangle \big)$ defines an inner product space.
\end{define}

\begin{define}[Zero and one functions]
Denote, respectively, by $\zeroindic$ and $\oneindic$ the discrete zero and the discrete one functions defined as
\begin{equation}
\begin{array}{lll}
\forall x \in \mathcal{L}^n, \,
\begin{cases}
\zeroindic(x) = 0 \\
\oneindic(x) = 1
\end{cases}
\end{array}
\nonumber
\end{equation}
\end{define}

\begin{define}[Positive semi-definite linear operator]
A linear operator ${\mathscr{A}}:\mathscr{F}\rightarrow \mathscr{F}$ is said to be positive semi-definite, if and only if,
\begin{equation}
\begin{array}{lll}
\forall f \in \mathscr{F}, \, \big<\mathscr{A} \, f,f\rangle  \geq 0
\end{array}
\nonumber
\end{equation}
\end{define}

\begin{define}[Strictly positive linear operator]
A linear operator ${\mathscr{A}}:\mathscr{F}\rightarrow \mathscr{F}$ is said to be strictly positive definite, if and only if,
\begin{equation}
\forall f \in \mathscr{F}, \, \big<\mathscr{A} \, f,f\rangle
\begin{cases}
>0, & \text{if  }  f \not \equiv \zeroindic, \\
=0, & \text{if  }  f \equiv \zeroindic.
\end{cases}
\nonumber
\end{equation}
\end{define}

\begin{define}[Multi-dimensional indicator function] $\forall m \in \mathbb{N}^{+}$, we define the $m$-dimensional indicator function
denoted by $\delta(\underbrace{\cdot,\cdots,\cdot}_{\text{m } times})$ as
\begin{equation}
\begin{array}{lll}
\delta(z_1,\dots,z_m) = \Pi_{i=1}^{m} \delta(z_i), \, \forall (z_1,\dots,z_m) \in \mathcal{L}^{m}
\end{array}
\nonumber
\end{equation}
where $\delta(\cdot)$ denotes the Kronecker delta function defined as
\begin{equation}
\forall a \in \mathcal{L}, \, \delta(a) =
\begin{cases}
1, \text{ if }  a = 0, \\
0, \text{ else.}
\end{cases}
\nonumber
\end{equation}
\end{define}

\begin{define}[Margin]
\label{margin_define}
Suppose $c$ is a hypervertex  of $\mathcal{G}$, and a function $u \in \mathscr{F}$. Then, we call the margin
of $u$ with respect to $c$ the (local) function $u_c: \mathcal{L}^{|c|} \rightarrow \mathbb{R}$ defined as
\begin{equation}
\begin{array}{l}
\label{form_margin_define}
u_c(x_c) = \sum_{i \in \mathcal{V}/c} \sum_{x_i \in \mathcal{L}} u(x_1,\dots,x_i,\dots,x_n), \, \forall x_c \in \mathcal{L}^{|c|}
\end{array}
\end{equation}
\end{define}

\begin{remark}
Throughout, we call marginalization of a function $u \in \mathscr{F}$ with respect to a hypervertex $c$ of
$\mathcal{G}$ the process of summing $u(x)$ across coordinates $x_i, \forall i \in \mathcal{V}/c$.
\end{remark}

\begin{define}[Marginal operator with regard to hypervertex $c$]
\label{sub:intro_oc}
Assume an arbitrary (perhaps, empty) hypervertex $c$ of $\mathcal{G}$, and introduce linear operator $\mathscr{O}_c:\mathscr{F} \rightarrow \mathscr{F}$
defined for all function $u \in \mathscr{F}$ as
\begin{equation}
\begin{array}{lll}
\big(\mathscr{O}_c \, u \big)(x) = \frac{u_c(x_c)}{L^{n-|c|}}, \,  \forall x \in \mathcal{L}^n
\end{array}  \nonumber
\end{equation}
where $u_c$ stands for the margin of $u$ with respect to $c$.
\end{define}

\begin{proposition}
\label{prop:o_s_ortho}
Suppose $c$ is a (perhaps empty) hypervertex of $\mathcal{G}$. Then, $\mathscr{O}_c$ is an orthogonal projection.
\end{proposition}
The proof of Proposition \ref{prop:o_s_ortho} is sketched in Appendix section \ref{app:proof_prop_o_s_ortho}.

Thus, whatever a hypervertex $c$ of $\mathcal{G}$, we refer in the remainder to $\mathscr{O}_c$ as the marginal operator with regard to $c$.

\begin{remark}
$\mathscr{O}_{\emptyset}$ is nothing else than the averaging operator over all function variables, since $\forall u \in \mathscr{F}$,
one has $\mathscr{O}_{\emptyset} \, u \equiv \Big(\frac{\sum_{x\in \mathcal{L}^n} u(x)}{L^n}\Big) \, \oneindic$.
\end{remark}

\begin{remark}
Suppose $c$ is a hypervertex  of $\mathcal{G}$. Then, an alternative definition of $\mathscr{O}_c$ as a symmetric kernel $:\mathcal{L}^n \times \mathcal{L}^n \rightarrow \mathbb{R}$ is as follows:
\begin{equation}
\begin{array}{lll}
\mathscr{O}_c(x,y) = \frac{ \delta\big(x_{c} - y_{c}\big)}{L^{n-|c|}}, \,  \forall (x,y) \in \mathcal{L}^n \times \mathcal{L}^n
\end{array}
\nonumber
\end{equation}
in such a way that, $\forall u \in \mathscr{F}$, and $\forall x \in \mathcal{L}^n$, $\big(\mathscr{O}_c\, u\big)(x)$ can also be
viewed as the right inner product between $\mathscr{O}_c(x,\cdot)$ and $u(\cdot)$ as follows:
\begin{equation}
\begin{array}{lll}
\big(\mathscr{O}_{c}\, u\big)(x) = \sum_{y \in \mathcal{L}^n} {\mathscr{O}}_{c}(x,y) \, u(y)
\end{array}
\nonumber
\end{equation}
\end{remark}

\begin{remark}
It is an easy exercise to show that MRF marginal distribution $\mathbb{P}\big(\mathcal{X}_s = \cdot \big)$ can be reexpressed, accordingly, as follows:
\begin{equation}
\begin{array}{lll}
\mathbb{P}\big(\mathcal{X}_s = \cdot \big) \equiv  L^{n-|s|} \, \mathscr{O}_s \,  \mathbb{P}\big(\mathcal{X} = \cdot \big)
\end{array}
\label{p_s_eq_os_p}
\end{equation}
\end{remark}

\begin{notation}
Whatever a hypervertex $c$ of $\mathcal{G}$, denote $\widebar{\mathscr{O}}_c = \idoper - \mathscr{O}_c$.
\end{notation}

Next, introduce the following two subspaces of $\mathscr{F}$:
\begin{equation}
\begin{cases}
\hat{\mathscr{F}}_s = \Big\{  \mathscr{O}_{s} \, f, \, \forall f \in \mathscr{F} \Big\} \\
\widebar{\mathscr{F}}_s = \Big\{ \widebar{\mathscr{O}_{s}} \, f, \, \forall f \in \mathscr{F} \Big\}
\end{cases}
\nonumber
\end{equation}

\begin{proposition}
\label{prop:F_direct_sum_marginal}
One has $\mathscr{F} = \hat{\mathscr{F}}_s \oplus \widebar{\mathscr{F}}_s$.
\end{proposition}
The proof of Proposition \ref{prop:F_direct_sum_marginal} is sketched in Appendix section \ref{app:proof_prop_F_direct_sum_marginal}.

\begin{proof}
This is a classical linear algebraic result following immediately from Proposition \ref{prop:o_s_ortho} which establishes
that $\mathscr{O}_s$ is an orthogonal projection, in such a way that, $\hat{\mathscr{F}}_s$ and $\widebar{\mathscr{F}}_s$
are two complementary orthogonal subspaces of $\mathscr{F}$.
\end{proof}

Thus, suppose an arbitrary function $f \in  \mathscr{F}$. Thus, Proposition \ref{prop:F_direct_sum_marginal} states
that there exists two unique functions $f_s \in \hat{\mathscr{F}}_s$, and $\bar{f}_s \in \widebar{\mathscr{F}}_s$,
such that, one can write
\begin{equation}
\begin{array}{lll}
f \equiv f_s + \bar{f}_s
\end{array}
\nonumber 
\end{equation}
Moreover, both $f_s$ and $\bar{f}_s$ are given by
\begin{equation}
\begin{cases}
f_s \equiv \mathscr{O}_s \, f \\
\bar{f}_s \equiv \widebar{\mathscr{O}}_s \, f
\end{cases}
\nonumber 
\end{equation}


\section{Implicit algebraic reformulation of a MRF joint distribution}
\label{sec:implicit_mrf_dist}
This section develops a framework to implicit linear algebraic reformulation of a general MRF model.

\begin{lemmon}
\label{lemma:besag}
Joint MRF distribution $\mathbb{P}\big(\mathcal{X} = \cdot\big)$ is uniquely determined among all joint MRF distributions with regard to graph $\mathcal{G}$ by its set of
$n$ individual conditional distributions as specified in formula (\ref{neighb_conditional}).
\end{lemmon}
The proof of Lemma \ref{lemma:besag} is found in Appendix section \ref{app:proof_lemma_besag}.


\begin{lemmon}
\label{lemma:equiv_pi_ki}
Suppose $u \in \mathscr{F}$, and $i \in \mathcal{V}$. Then, the following statement:
{\small \begin{equation}
\begin{array}{l}
 u(x) -  \pi_{i}(x_{{\ne}_{i}^+})  \sum_{x_i \in \mathcal{L}} u(x) = 0, \forall x \in \mathcal{L}^n
\Leftrightarrow  u(x) - \widehat{\pi}_i(x_{{\ne}_{i}^+}) \sum_{x_i\in \mathcal{L}} \widehat{\pi}_i(x_{{\ne}_{i}^+}) \, u(x)  = 0,  \forall x \in \mathcal{L}^n
\end{array}
\nonumber
\end{equation}}
holds.
\end{lemmon}
Lemma \ref{lemma:equiv_pi_ki} is shown in Appendix section \ref{app:proof_lemma_equiv_pi_ki}.

\begin{lemmon}
\label{lemma:unicity_mrf}
Joint MRF distribution $\mathbb{P}\big(\mathcal{X} = \cdot \big)$ is the unique solution of the following system of linear equations with a nonnegativity constraint:
\begin{equation}
\begin{array}{l}
\begin{cases}
p(x) -  \widehat{\pi}_i(x_{{\ne}_{i}^+}) \, \sum_{x_i\in \mathcal{L}} \widehat{\pi}_i(x_{{\ne}_{i}^+}) \, p(x) = 0, \, \forall x \in \mathcal{L}^n,\, \forall i \in \mathcal{V} \\
\sum_{x \in \mathcal{L}^n} p(x) = 1 \\
p(x)\geq 0 , \, \forall x \in \mathcal{L}^n
\end{cases}
\end{array}
\nonumber
\end{equation}
\end{lemmon}
The proof of Lemma \ref{lemma:unicity_mrf} is found in Appendix section \ref{app:proof_lemma_unicity_mrf}.

\begin{lemmon}
\label{lemma:max_inf_min_sup_vi}
Suppose an arbitrary function $f \in \mathscr{F}$. One then has
\begin{equation}
\begin{array}{lll}
\mathbb{E}\big[f(\mathcal{X})\big]
&=& \underset{\big\{ v_i\in \mathscr{F}, \forall i \in \mathcal{V} \big\}}{\max}
\, \underset{x \in \mathcal{L}^n}{\inf} \Big\{ f(x) + \sum_{i \in \mathcal{V}} \Big(\, v_i(x) - \widehat{\pi}_i(x_{{\ne}_{i}^+}) \, \sum_{x_i\in \mathcal{L}} \widehat{\pi}_i(x_{{\ne}_{i}^+}) \, v_i(x) \, \Big) \Big\} \\
&=& \underset{\big\{ v_i\in \mathscr{F}, \forall i \in \mathcal{V} \big\}}{\min}
\, \underset{x \in \mathcal{L}^n}{\sup}  \Big\{ f(x) + \sum_{i \in \mathcal{V}} \Big( \, v_i(x) - \widehat{\pi}_i(x_{{\ne}_{i}^+}) \, \sum_{x_i\in \mathcal{L}} \widehat{\pi}_i(x_{{\ne}_{i}^+}) \, v_i(x) \, \Big) \Big\}
\end{array}
\nonumber 
\end{equation}
\end{lemmon}
The proof of Lemma \ref{lemma:max_inf_min_sup_vi} is described in Appendix section \ref{subap:proof_lemma_max_inf_min_sup_vi}.

\begin{lemmon}
\label{lemma:founding_result_mrf}
Suppose an arbitrary function $f \in \mathscr{F}$. Then, there exists $n$ functions $v_i \in \mathscr{F}$, $\forall i \in \mathcal{V}$, such that, one can write
\begin{equation}
\begin{array}{lll}
f(x) &=& \mathbb{E}\big[f(\mathcal{X})\big] + \sum_{i \in \mathcal{V}} \Big( \, v_i(x) -\widehat{\pi}_i(x_{{\ne}_{i}^+}) \, \sum_{x_i\in \mathcal{L}} \widehat{\pi}_i(x_{{\ne}_{i}^+}) \, v_i(x) \, \Big), \,
\forall x \in \mathcal{L}^n
\end{array}
\nonumber
\end{equation}
\end{lemmon}
The proof of Lemma \ref{lemma:founding_result_mrf} is detailed in Appendix section \ref{app:proof_lemma_founding_result_mrf}.

Now, introduce linear operator $\mathscr{K}: \mathscr{F} \rightarrow \mathscr{F}$ defined for all function $u \in \mathscr{F}$ as
\begin{equation}
\begin{array}{lll}
\big(\mathscr{K} (u)\big)(x) = \sum_{x \in \mathcal{L}^n} u(x) +
\sum_{i \in \mathcal{V}} \Big(\, u(x) - \widehat{\pi}_i(x_{{\ne}_{i}^+}) \, \sum_{x_i\in \mathcal{L}} \widehat{\pi}_i(x_{{\ne}_{i}^+}) \, u(x) \, \Big), \,  \forall x \in \mathcal{L}^n
\end{array}
\nonumber 
\end{equation}
\begin{remark}
An alternative definition of $\mathscr{K}$ as a positive definite kernel $:\mathcal{L}^n \times \mathcal{L}^n \rightarrow \mathbb{R}$,
is as follows:
\begin{equation}
\begin{array}{lll}
  \mathscr{K}(x,y) = 1 + \sum_{i \in \mathcal{V}} \Big(\, \delta(x-y) - \widehat{\pi}_i(x_{{\ne}_{i}^+}) \, \widehat{\pi}_i(y_{{\ne}_{i}^+}) \,
  \delta(x_{\neq i}-y_{\neq i}) \, \Big), \, \forall (x,y) \in \mathcal{L}^n \times \mathcal{L}^n
\end{array}
\nonumber
\end{equation}
in such a way that, $\forall u \in \mathscr{F}$, and $\forall x \in \mathcal{L}^n$, $\big(\mathscr{K} (u)\big)(x)$ can also be viewed as the right inner product between
$\mathscr{K}(x,\cdot)$ and $u(\cdot)$ as follows:
\begin{equation}
\begin{array}{lll}
\big(\mathscr{K} (u)\big)(x) = \sum_{y \in \mathcal{L}^n} \mathscr{K}(x,y) \, u(y), \, \forall x \in \mathcal{L}^n
\end{array}
\nonumber
\end{equation}
\end{remark}

Furthermore, introduce linear operators $k_i: \mathscr{F}\rightarrow \mathscr{F}$, $\forall i \in \mathcal{V}$ defined $\forall u \in  \mathscr{F}$ as
\begin{equation}
\big(k_i \, u \big)(x) = u(x) - \widehat{\pi}_i(x_{{\ne}_{i}^+}) \, \sum_{x_i\in \mathcal{L}} \widehat{\pi}_i(x_{{\ne}_{i}^+}) \, u(x), \, \forall x \in \mathcal{L}^n, \, \forall i \in \mathcal{V}
\nonumber
\end{equation}

\begin{proposition}
\label{prop:ki_is_ortho_project}
$\forall i \in \mathcal{V}$, $k_i$ is an orthogonal projection.
\end{proposition}
The proof of Proposition \ref{prop:ki_is_ortho_project} is sketched in Appendix section \ref{app:proof_prop_ki_is_ortho_project}.

\begin{proposition}
\label{prop:K_k_i}
Following linear operator equality $\mathscr{K} = L^n \, \mathscr{O}_{\emptyset} + \sum_{i \in \mathcal{V}} k_i$ holds.
\end{proposition}
The proof of Proposition \ref{prop:K_k_i} is found in Appendix section \ref{app:proof_prop_K_k_i}

\begin{lemmon}
\label{lemma:mrf_syst_nondegenerate}
$\lambda \equiv \zeroindic$ is the unique solution of the following system of linear identities:
\begin{equation}
\begin{array}{lll}
\begin{cases}
\lambda(x) - \widehat{\pi}_i(x_{{\ne}_{i}^+}) \, \sum_{x_i\in \mathcal{L}} \widehat{\pi}_i(x_{{\ne}_{i}^+}) \,  \lambda(x) = 0, \, \forall x \in \mathcal{L}^n, \, \forall i \in \mathcal{V} \\
\sum_{x \in \mathcal{L}^n} \lambda(x) = 0
\end{cases}
\end{array}
\nonumber 
\end{equation}
\end{lemmon}
The proof of Lemma \ref{lemma:mrf_syst_nondegenerate} is described in Appendix section \ref{app:proof_lemma_mrf_syst_nondegenerate}.

\begin{lemmon}
\label{lemma:markov_lin_oper}
$\mathscr{K}$ is a self-adjoint strictly positive definite (hence, invertible) linear operator.
\end{lemmon}
The proof of Lemma \ref{lemma:markov_lin_oper} is detailed in Appendix section \ref{app:proof_lemma_markov_lin_oper}.


\begin{theorem}
\label{theo:markov_distrib_unique_pos_definite}
MRF distribution $\mathbb{P}\big(\mathcal{X} = \cdot \big)$ is the unique solution of the following linear identity:
\begin{equation}
\begin{array}{lll}
\mathscr{K} \, p \equiv \oneindic
\end{array}
\label{mrf_dist_solu_lin_sys_0}
\end{equation}
\end{theorem}
The proof of Theorem \ref{theo:markov_distrib_unique_pos_definite} is sketched in Appendix section \ref{app:proof_theo_markov_distrib_unique_pos_definite}.

Thus, Theorem \ref{theo:markov_distrib_unique_pos_definite} establishes that it is possible to completely reexpress MRF joint
distribution $\mathbb{P}\big(\mathcal{X} = \cdot \big)$ as the unique solution of linear identity (\ref{mrf_dist_solu_lin_sys_0}).
As an aside, since by Lemma \ref{lemma:markov_lin_oper}, $\mathscr{K}$ is an invertible linear operator, one can write
\begin{equation}
\begin{array}{lll}
\mathbb{P}\big(\mathcal{X} = x \big) = \Big( \mathscr{K}^{-1} \, \oneindic \Big)(x), \, \forall x \in \mathcal{L}^n
\end{array}
\nonumber 
\end{equation}
where $\mathscr{K}^{-1}$ stands for the inverse of linear operator $\mathscr{K}$. Obviously, if one had access
to the closed-form expression of $\mathscr{K}^{-1}$ or, at least, if one could compute the "smaller" linear
operator $\mathscr{O}_s \, \mathscr{K}^{-1} \, \mathscr{O}_s$ efficiently, then one could solve MRF marginal
inference problem (\ref{marginal_inference_mrf}) efficiently as follows:
\begin{equation}
\begin{array}{lll}
\mathbb{P}\big(\mathcal{X}_s = \cdot \big) &\equiv&  L^{n-|s|} \, \big( \mathscr{O}_s \, \mathscr{K}^{-1} \big) \,  \oneindic \\
                                           &\equiv&  L^{n-|s|} \, \big( \mathscr{O}_s \, \mathscr{K}^{-1} \big) \, \big(  \mathscr{O}_s \, \oneindic \big) \\
                                           &\equiv&  L^{n-|s|} \, \big( \mathscr{O}_s \, \mathscr{K}^{-1} \, \mathscr{O}_s  \big) \,  \oneindic
\end{array}
\nonumber 
\end{equation}
Nevertheless, our task in the remainder is, precisely, to attempt to circumvent such a seemingly hard direct inversion formula by means of
a useful relaxation of linear identity (\ref{mrf_dist_solu_lin_sys_0}) in a way which further enables to recast MRF marginal inference problem
(\ref{marginal_inference_mrf}) as an easy algebraic problem. 


\section{Relaxation for the sake of marginal inference}
\label{sec:relaxation}


\subsection{Preliminary results}

First, introduce linear operator $\hat{\mathscr{K}}^{(s)}:\mathscr{F} \rightarrow \mathscr{F}$ indexed with hypervertex $s$, and defined for all function $u \in \mathscr{F}$ as
\begin{equation}
\begin{array}{lll}
\big(\hat{\mathscr{K}}^{(s)} \, u\big)(x) = u_{\emptyset} +  \sum_{i \in \mathcal{V}} \bigg( \, \Big(\frac{{u}_{s\cup {\ne}_{i}^+}(x_{s\cup {\ne}_{i}^+})}{L^{n-|s\cup {\ne}_{i}^+|}} \Big) -
\widehat{\pi}_i(x_{{\ne}_{i}^+}) \, \sum_{x_i\in \mathcal{L}} \widehat{\pi}_i(x_{{\ne}_{i}^+}) \, \Big( \frac{{u}_{s\cup {\ne}_{i}^+}(x_{s\cup {\ne}_{i}^+})}{L^{n-|s\cup {\ne}_{i}^+|}}\Big) \,\bigg), \,  \forall x \in \mathcal{L}^n
\end{array}  \nonumber
\end{equation}
where $u_{\emptyset} = \sum_{x \in \mathcal{L}^n} u(x)$, and $\forall i \in \mathcal{V}$, ${u}_{s\cup {\ne}_{i}^+}$ stands for the margin of $u$
with respect to hypervertex $s\cup {\ne}_{i}^+$.

\begin{remark}
Denote $\forall i \in \mathcal{V}$ by $s^{\neq i} = s /\{i\}$.
It follows that an alternative definition of $\hat{\mathscr{K}}^{(s)}$ as a symmetric function $:\mathcal{L}^n \times \mathcal{L}^n \rightarrow \mathbb{R}$ is as follows:
\begin{equation}
\hat{\mathscr{K}}^{(s)}(x,y) = 1 + \sum_{i \in \mathcal{V}} \frac{ \delta\big(x_{s\cup {\ne}_{i}^+}-y_{s\cup {\ne}_{i}^+}\big) -
\widehat{\pi}_i(x_{{\ne}_{i}^+}) \, \widehat{\pi}_i(y_{{\ne}_{i}^+}) \, \delta\big(x_{s^{\neq i} \cup {\ne}_{i}}-y_{s^{\neq i} \cup {\ne}_{i}}\big)}{L^{n-|s\cup {\ne}_{i}^+|}}, \,  \forall (x,y) \in \mathcal{L}^n \times \mathcal{L}^n
\nonumber
\end{equation}
in such a way that, $\forall u \in \mathscr{F}$, and $\forall x \in \mathcal{L}^n$, $\big(\hat{\mathscr{K}}^{(s)}\, u\big)(x)$ can also be viewed as the right inner product
between $\hat{\mathscr{K}}^{(s)}(x,\cdot)$ and $u(\cdot)$ as follows:
\begin{equation}
\begin{array}{lll}
\big(\hat{\mathscr{K}}^{(s)}\, u\big)(x) = \sum_{y \in \mathcal{L}^n} \hat{\mathscr{K}}^{(s)}(x,y) \, u(y)
\end{array}
\nonumber
\end{equation}
\end{remark}

Also, introduce $\forall i \in \mathcal{V}$ linear operator $\widehat{k}_i^{(s)}: \mathscr{F}\rightarrow \mathscr{F}$, defined $\forall u \in \mathscr{F}$ as
\begin{equation}
\big(\widehat{k}_i^{(s)} \, u \big)(x) =  \Big(\frac{{u}_{s\cup {\ne}_{i}^+}(x_{s\cup {\ne}_{i}^+})}{L^{n-|s\cup {\ne}_{i}^+|}} \Big) -
\widehat{\pi}_i(x_{{\ne}_{i}^+}) \, \sum_{x_i\in \mathcal{L}}  \widehat{\pi}_i(x_{{\ne}_{i}^+}) \, \Big(\frac{{u}_{s\cup {\ne}_{i}^+}(x_{s\cup {\ne}_{i}^+})}{L^{n-|s\cup {\ne}_{i}^+|}} \Big), \, \forall x \in \mathcal{L}^n
\nonumber
\end{equation}

\begin{proposition}
\label{prop:kilde_eq}
Following linear operator equalities:
\begin{equation}
\begin{array}{lll}
\widehat{k}_i^{(s)} =  {k}_i \,  \mathscr{O}_{s\cup {\ne}_{i}^+} =  \mathscr{O}_{s\cup {\ne}_{i}^+} \, {k}_i = \mathscr{O}_{s\cup {\ne}_{i}^+} \, {k}_i \,  \mathscr{O}_{s\cup {\ne}_{i}^+}, \, \forall i \in \mathcal{V}
\end{array}
\nonumber
\end{equation}
hold.
\end{proposition}
The proof of Proposition \ref{prop:kilde_eq} is sketched in Appendix section \ref{app:proof_prop_kilde_eq}.

\begin{proposition}
\label{prop:kilde_i_is_ortho_project}
$\forall i \in \mathcal{V}$, $\widehat{k}_i^{(s)}$ is an orthogonal projection.
\end{proposition}
The proof of Proposition \ref{prop:kilde_i_is_ortho_project} is sketched in Appendix section \ref{app:proof_prop_kilde_i_is_ortho_project}.

\begin{proposition}
\label{prop:K_k_tilde_i}
Following linear operator equality $\hat{\mathscr{K}}^{(s)} = L^n \, \mathscr{O}_{\emptyset} + \sum_{i \in \mathcal{V}} \widehat{k}_i^{(s)}$ holds.
\end{proposition}
The proof of Proposition \ref{prop:K_k_tilde_i} is sketched in Appendix section \ref{app:proof_prop_K_k_tilde_i}.

\begin{proposition}
\label{prop:markov_lin_oper_s}
$\hat{\mathscr{K}}^{(s)}$ is a self-adjoint positive semi-definite linear operator.
\end{proposition}
The proof of Proposition \ref{prop:markov_lin_oper_s} is detailed in Appendix section \ref{app:proof_prop_markov_lin_oper_s}.

\begin{proposition}
\label{prop:ks_eq_kos}
Following linear operator equality $\hat{\mathscr{K}}^{(s)}\, \mathscr{O}_{s} = {\mathscr{K}}\, \mathscr{O}_{s}$ holds.
\end{proposition}
The proof of Proposition \ref{prop:ks_eq_kos} is sketched in Appendix section \ref{app:proof_prop_ks_eq_kos}.

\begin{proposition}
\label{prop:kts_p_sar_eq_one}
Joint MRF distribution $\mathbb{P}\big(\mathcal{X} = \cdot \big)$ satisfies $\hat{\mathscr{K}}^{(s)}\, \mathbb{P}\big(\mathcal{X} = \cdot \big) \equiv \oneindic$.
\end{proposition}
The proof of Proposition \ref{prop:kts_p_sar_eq_one} is sketched in Appendix section \ref{app:proof_kts_p_sar_eq_one}.

Now, introduce linear operator $\widebar{\mathscr{K}}^{(s)}:\mathscr{F} \rightarrow \mathscr{F}$ defined as $\widebar{\mathscr{K}}^{(s)} =  \mathscr{K} - \hat{\mathscr{K}}^{(s)}$.
Equivalently, one has $\forall u \in \mathscr{F}$
\begin{equation}
\begin{array}{lll}
\big(\widebar{\mathscr{K}}^{(s)} \, u\big)(x) =  \sum_{i \in \mathcal{V}} \bigg( \, \Big(u(x) - \frac{{u}_{s\cup {\ne}_{i}^+}(x_{s\cup {\ne}_{i}^+})}{L^{n-|s\cup {\ne}_{i}^+|}} \Big) -
\widehat{\pi}_i(x_{{\ne}_{i}^+}) \, \sum_{x_i\in \mathcal{L}} \widehat{\pi}_i(x_{{\ne}_{i}^+}) \,
\Big( u(x) - \frac{{u}_{s\cup {\ne}_{i}^+}(x_{s\cup {\ne}_{i}^+})}{L^{n-|s\cup {\ne}_{i}^+|}}\Big) \,\bigg), \,  \forall x \in \mathcal{L}^n
\end{array}  \nonumber
\end{equation}
where $\forall i \in \mathcal{V}$, ${u}_{s\cup {\ne}_{i}^+}$ denotes the margin of $u$ with respect to hypervertex $s\cup {\ne}_{i}^+$.

Also,  $\forall i \in \mathcal{V}$, introduce linear operator $\bar{k}_i^{(s)}: \mathscr{F}\rightarrow \mathscr{F}$, defined $\forall u \in \mathscr{F}$ as
\begin{equation}
\big(\bar{k}_i^{(s)} \, u \big)(x) =  \Big( u(x) - \frac{{u}_{s\cup {\ne}_{i}^+}(x_{s\cup {\ne}_{i}^+})}{L^{n-|s\cup {\ne}_{i}^+|}} \Big) -
\widehat{\pi}_i(x_{{\ne}_{i}^+}) \, \sum_{x_i\in \mathcal{L}} \widehat{\pi}_i(x_{{\ne}_{i}^+}) \,  \Big(  u(x) - \frac{{u}_{s\cup {\ne}_{i}^+}(x_{s\cup {\ne}_{i}^+})}{L^{n-|s\cup {\ne}_{i}^+|}} \Big), \, \forall x \in \mathcal{L}^n
\nonumber
\end{equation}

\begin{proposition}
\label{prop:kbar_i_is_ortho_project}
One has $\forall i \in \mathcal{V}$, $\bar{k}_i^{(s)}$ is an orthogonal projection.
\end{proposition}
The proof of Proposition \ref{prop:kbar_i_is_ortho_project} is sketched in Appendix section \ref{app:proof_prop_kbar_i_is_ortho_project}.

\begin{proposition}
\label{prop:kbar_i_is_sum_ortho_project}
One has $\widebar{\mathscr{K}}^{(s)} = \sum_{i \in \mathcal{V}} \bar{k}_i^{(s)}$
\end{proposition}
The proof of Proposition \ref{prop:kbar_i_is_sum_ortho_project} is found in Appendix section \ref{app:proof_prop_kbar_i_is_sum_ortho_project}.

\begin{proposition}
\label{prop:markov_lin_oper_s_r}
$\widebar{\mathscr{K}}^{(s)}$ is a self-adjoint positive semi-definite linear operator.
\end{proposition}
The proof of Proposition \ref{prop:markov_lin_oper_s_r} is found in Appendix section \ref{app:proof_prop_markov_lin_oper_s_r}.

\begin{proposition}
\label{prop:rs_eq_rosb}
Following linear operator equalities $\widebar{\mathscr{K}}^{(s)} = \widebar{\mathscr{K}}^{(s)}\, \widebar{\mathscr{O}}_{s} = \widebar{\mathscr{O}}_{s} \, \widebar{\mathscr{K}}^{(s)}
= \widebar{\mathscr{O}}_{s} \, \widebar{\mathscr{K}}^{(s)}\, \widebar{\mathscr{O}}_{s}$ hold.
\end{proposition}
The proof of Proposition \ref{prop:rs_eq_rosb} is sketched in Appendix section \ref{app:proof_prop_rs_eq_rosb}.

\begin{proposition}
\label{prop:rts_p_sar_eq_zero}
Joint MRF distribution $\mathbb{P}\big(\mathcal{X} = \cdot \big)$ satisfies $\widebar{\mathscr{K}}^{(s)}\, \mathbb{P}\big(\mathcal{X} = \cdot \big) \equiv \zeroindic$.
\end{proposition}
The proof of Proposition \ref{prop:rts_p_sar_eq_zero} is found in Appendix section \ref{app:proof_prop_rts_p_sar_eq_zero}.

\begin{convention}
For notational convenience in the remainder, unless specified otherwise, we shall assume that hypervertex "$s$" is fixed, once and for all, and avoid referring to $s$
in related mathematical symbols. For instance, we shall simply write $\hat{\mathscr{K}}$ instead of $\hat{\mathscr{K}}^{(s)}$,  and also write
$\widebar{\mathscr{K}}$ instead of $\widebar{\mathscr{K}}^{(s)}$.
\end{convention}

\begin{corollary}
\label{corol:p_sol_kp_1_bkp0}
Joint MRF distribution $\mathbb{P}\big(\mathcal{X} = \cdot \big)$ is the unique solution of the following system of linear identities:
\begin{equation}
\begin{cases}
\hat{\mathscr{K}} \,p  \equiv \oneindic \\
\widebar{\mathscr{K}} \,p  \equiv \zeroindic \\
\end{cases}
\label{identity_kp_1_kpb_zero}
\end{equation}
\end{corollary}
The proof of Corollary \ref{corol:p_sol_kp_1_bkp0} is sketched in Appendix section \ref{app:proof_corol_p_sol_kp_1_bkp0}.

\begin{notation}
Denote in the remainder by $\hat{\mathscr{K}}^{(\dag)}$ the Moore-Penrose inverse (or, the pseudoinverse) of linear operator $\hat{\mathscr{K}}$.
\end{notation}

\subsection{Founding result}  
\label{sub:founding_result}
The founding result in this paper which sets the stage for polynomial time resolution of MRF marginal inference problem (\ref{marginal_inference_mrf}) on
algebraic grounds is presented in Theorem \ref{theorem:founding_kp_eq_1}.

\begin{lemmon}
\label{lemma:kq_eq_kbar_q}
The following identity:
\begin{equation}
\begin{array}{lll}
\widebar{\mathscr{K}} \, \big(\idoper - \hat{\mathscr{K}}^{\dag} \, \hat{\mathscr{K}} \big) \,  {q} &\equiv&
\mathscr{K} \, \big(\idoper - \hat{\mathscr{K}}^{\dag} \, \hat{\mathscr{K}} \big) \,  {q},
\end{array}
\nonumber
\end{equation}
holds for all $q \in \mathscr{F}$.
\end{lemmon}
Lemmas \ref{lemma:kq_eq_kbar_q} is shown in Appendix sections \ref{app:proof_lemma_kq_eq_kbar_q}.

\begin{lemmon}
\label{lemma:i_kk_osq_eq_func_kkosq}
The following two identities:
\begin{equation}
\big(\idoper - \hat{\mathscr{K}}^{\dag} \, \hat{\mathscr{K}} \big) \,   \big({\mathscr{O}}_s \, q \big)
\equiv -\mathscr{K}^{-1} \, \widebar{\mathscr{K}} \, \hat{\mathscr{K}}^{\dag} \, \hat{\mathscr{K}} \,  \big({\mathscr{O}}_s \, q \big)
\nonumber
\end{equation}
and
\begin{equation}
\hat{\mathscr{K}} \, \mathscr{K}^{-1} \, \widebar{\mathscr{K}} \, \hat{\mathscr{K}}^{\dag} \, \hat{\mathscr{K}} \,  \big({\mathscr{O}}_s \, q \big) \equiv \zeroindic
\nonumber
\end{equation}
hold for all $q \in \mathscr{F}$.
\end{lemmon}
Lemma \ref{lemma:i_kk_osq_eq_func_kkosq} is shown in Appendix section \ref{app:proof_lemma_i_kk_osq_eq_func_kkosq}.

\begin{theorem}
\label{theorem:founding_kp_eq_1}
Suppose $p \in \mathscr{F}$ to be any solution to the following linear identity:
\begin{equation}
\begin{array}{lll}
\hat{\mathscr{K}} \, p \equiv \oneindic
\end{array}
\label{mrf_dist_solu_lin_sys_polyn}
\end{equation}
and denote by $p_s: \mathcal{L}^{|s|} \rightarrow \mathbb{R}$ the margin of $p$ with respect to $s$, that is,
\begin{equation}
\begin{array}{l}
\nonumber
p_s(x_s) = \sum_{i \in \mathcal{V}/s} \sum_{x_i \in \mathcal{L}} p(x_1,\dots,x_i,\dots,x_n), \, \forall x_s \in \mathcal{L}^{|s|}
\end{array}
\end{equation}
One then has
\begin{equation}
\begin{array}{lll}
\mathbb{P}\big(\mathcal{X}_s = x_s \big) = p_s(x_s), \, \forall x_s \in \mathcal{L}^{|s|}
\end{array}
\nonumber 
\end{equation}
\end{theorem}
The proof of Theorem \ref{theorem:founding_kp_eq_1} is detailed in Appendix section \ref{app:proof_theorem_founding_kp_eq_1}.

Thus, Theorem \ref{theorem:founding_kp_eq_1} establishes that it is possible to solve exactly for MRF marginal distribution
$\mathbb{P}\big(\mathcal{X}_s = \cdot \big)$, first by finding any solution $q: \mathcal{L}^n \rightarrow \mathbb{R}$ to linear
identity (\ref{mrf_dist_solu_lin_sys_polyn}), second by computing the margin of $q$ with regard to $s$, and one is guaranteed
to have
\begin{equation}
\mathbb{P}\big(\mathcal{X}_s = \cdot \big) \equiv q_s
\nonumber
\end{equation}
with $q_s$ denoting the margin of $q$ with respect to $s$.

As an aside, we have the following results.
\begin{proposition}
\label{prop:kkos_eq_kk}
The following equalities between linear operators:
\begin{equation}
\big(\hat{\mathscr{K}}^{\dag} \, \hat{\mathscr{K}} \big) \,   {\mathscr{O}}_s = {\mathscr{O}}_s \,
\big(\hat{\mathscr{K}}^{\dag} \, \hat{\mathscr{K}} \big) = {\mathscr{O}}_s
\nonumber
\end{equation}
hold.
\end{proposition}

The proof of Proposition \ref{prop:kkos_eq_kk} is sketched in Appendix section \ref{app:proof_prop_kkos_eq_kk}.

Amazingly, the appeal of Proposition \ref{prop:kkos_eq_kk} is that whatever a local function $f_s(x_s), \, \forall x_s \in \mathcal{L}^{|s|}$,
there exists a function $q \in \mathscr{F}$, such that, $f_s \equiv \hat{\mathscr{K}} \, q$, that is to say
\begin{eqnarray}
f_s(x_s) = q_{\emptyset} +  \sum_{i \in \mathcal{V}} \Big( \, \Big(\frac{{q}_{s\cup {\ne}_{i}^+}(x_{s\cup {\ne}_{i}^+})}{L^{n-|s\cup {\ne}_{i}^+|}} \Big) -
\widehat{\pi}_i(x_{{\ne}_{i}^+}) \, \sum_{x_i\in \mathcal{L}} \widehat{\pi}_i(x_{{\ne}_{i}^+}) \, \Big( \frac{{q}_{s\cup {\ne}_{i}^+}(x_{s\cup {\ne}_{i}^+})}{L^{n-|s\cup {\ne}_{i}^+|}}\Big) \,\Big),
\, \forall x \in \mathcal{L}^{n}
\label{f_s_eq_k_q}
\end{eqnarray}
in such a way that, by convolving both hand sides of identity (\ref{f_s_eq_k_q}) with $\mathbb{P}\big(\mathcal{X} = \cdot \big)$, one derives
\begin{eqnarray}
\mathbb{E}\big[f_s(\mathcal{X}_s)\big] &= & q_{\emptyset} \\
&=& \sum_{x \in \mathcal{L}^n} q(x)
\nonumber
\end{eqnarray}
meaning that computing $\mathbb{E}\big[f_s(\mathcal{X}_s)\big]$ amounts to finding a function $q \in \mathscr{F}$ satisfying identity (\ref{f_s_eq_k_q}).

\begin{lemmon}
\label{lemma:local_consistency}
Suppose $c_1, \dots, c_m$ are $m$ distinct hypervertices of graph $\mathcal{G}$, with $m$ standing for some natural number $m$ greater than,
or equal to $2$, moreover, assume $m$ local functions of the form $q_i: \mathcal{L}^{|c_i|} \rightarrow \mathbb{R}$, $\forall i\in \mathcal{V}$
satisfying the following local consistency constraints:
\begin{equation}
\begin{array}{lll}
\sum_{k \in \mathcal{V}/c_i} \sum_{x_k \in \mathcal{L}} q_i(x_{c_i}) = \sum_{k \in \mathcal{V}/c_j} \sum_{x_k \in \mathcal{L}} q_j(x_{c_j}), \, \forall x_{c_i \cap c_j} \in \mathcal{L}^{|c_i \cap c_j|},
\forall i,j \in \big\{1,\dots, m\big\}, \text{s.t.,} i \neq j
\end{array}
\nonumber
\end{equation}
Then, there exists a (global) function $q:\mathcal{L}^{n} \rightarrow \mathbb{R}$, such that, the margin of $q$ with respect to $c_i$ coincides exactly
with $q_i$, $\forall i \in \big\{1,\dots, m\big\}$, that is to say, one has
\begin{equation}
\begin{array}{lll}
\sum_{k \in \mathcal{V}/c_i} \sum_{x_k \in \mathcal{L}} q(x) = q_i(x_{c_i}), \, \forall x_{c_i} \in \mathcal{L}^{|c_i|}, \forall i \in \big\{1,\dots, m\big\}
\end{array}
\nonumber
\end{equation}
\end{lemmon}
The proof of Lemma \ref{lemma:local_consistency} is described in Appendix section \ref{app:proof_lemma_local_consistency}.

Thus, by Lemma \ref{lemma:local_consistency}, solving for $q$ in identity (\ref{f_s_eq_k_q}) may be achieved in a computationally seductive way by,
instead of solving directly for a general function $q \in \mathscr{F}$ in identity (\ref{f_s_eq_k_q}), one rather solves for the margins of $q$
with respect to the following hypervertices:
\begin{equation}
q_{\emptyset},  \, \, {q}_{s\cup {\ne}_{i}^+}(x_{s\cup {\ne}_{i}^+}),  \forall x_{s\cup {\ne}_{i}^+}
\in \mathcal{L}^{|s\cup {\ne}_{i}^+|}, \, \forall i \in \mathcal{V}
\label{resol:margins}
\end{equation}
Moreover, the latter may be efficiently encoded by simply enforcing the local consistency constraints. Regarding the time complexity for solving
for such marginals (\ref{resol:margins}) in identity (\ref{f_s_eq_k_q}), this shall be answered in the coming section, but, we may say upfront that,
under Assumptions \ref{assump_computationability}, the time complexity for solving for such marginals (\ref{resol:margins}) in identity (\ref{f_s_eq_k_q})
is polynomial.

The following Proposition gives an explicit expression for $\mathbb{P}\big(\mathcal{X}_s = \cdot \big)$.
\begin{proposition}
\label{prop:marginal_qp}
Denote $p^* = \mathbb{P}\big(\mathcal{X} = \cdot \big)$. Then, one has
\begin{equation}
\frac{\mathbb{P}\big(\mathcal{X}_s = \cdot \big)}{L^{|s|}} \equiv {\mathscr{O}}_s \, p^* \equiv {\mathscr{O}}_s \, \big( \hat{\mathscr{K}}^{\dag} \, \oneindic)
\nonumber
\end{equation}
\end{proposition}
The proof of Proposition \ref{prop:marginal_qp} is found in Appendix section \ref{app:proof_prop_marginal_qp}

\section{Polynomial time resolution of the MRF marginal inference problem}
\label{sec:complexity_algorithm}

\begin{lemmon}
\label{lemma:equiv_pi_ki_tilde}
Suppose $\forall i \in \mathcal{V}$, $u_{s\cup {\ne}_{i}^+}$ is a local function $:\mathcal{L}^{|s\cup {\ne}_{i}^+|} \rightarrow \mathbf{R}$. Then, the
following statement:
\begin{multline}
 u_{s\cup {\ne}_{i}^+}(x_{s\cup {\ne}_{i}^+}) -  \pi_{i}(x_{{\ne}_{i}^+})  \sum_{x_i \in \mathcal{L}} u_{s\cup {\ne}_{i}^+}(x_{s\cup {\ne}_{i}^+}) = 0,
\forall x_{s\cup {\ne}_{i}^+} \in \mathcal{L}^{|s\cup {\ne}_{i}^+|} \\
\Leftrightarrow  u_{s\cup {\ne}_{i}^+}(x_{s\cup {\ne}_{i}^+}) - \widehat{\pi}_i(x_{{\ne}_{i}^+}) \sum_{x_i\in \mathcal{L}} \widehat{\pi}_i(x_{{\ne}_{i}^+}) \, u_{s\cup {\ne}_{i}^+}(x_{s\cup {\ne}_{i}^+}) = 0,
\forall x_{s\cup {\ne}_{i}^+} \in \mathcal{L}^{|s\cup {\ne}_{i}^+|}
\nonumber
\end{multline}
holds.
\end{lemmon}
Lemma \ref{lemma:equiv_pi_ki_tilde} is shown in Appendix section \ref{app:proof_lemma_equiv_pi_ki_tilde}.

\begin{lemmon}
\label{lemma:rectang_reduction_spatial}
Suppose $p \in \mathscr{F}$. It follows that $p$ is a solution to linear identity (\ref{mrf_dist_solu_lin_sys_polyn}), if and only if,
$p$ is a solution to the following system of linear equations:
\begin{equation}
\begin{cases}
{p}_{s\cup {\ne}_{i}^+}(x_{s\cup {\ne}_{i}^+}) - {\pi}_i(x_{{\ne}_{i}^+}) \, \sum_{x_i \in \mathcal{L}} {p}_{s\cup {\ne}_{i}^+}(x_{s\cup {\ne}_{i}^+}) = 0, \,
\forall x_{s\cup {\ne}_{i}^+} \in \mathcal{L}^{|s\cup {\ne}_{i}^+|}, \, \forall i \in \mathcal{V} \\
{p}_{s\cup {\ne}_{i}^+}(x_{s\cup {\ne}_{i}^+}) = \sum_{j \in \mathcal{V}/{s\cup {\ne}_{i}^+}} \sum_{x_j \in \mathcal{L}} p(x_1,\dots,x_j,\dots,x_n), \forall x_{s\cup {\ne}_{i}^+}
\in \mathcal{L}^{|s\cup {\ne}_{i}^+|},  \forall i \in \mathcal{V} \\
\sum_{x \in \mathcal{L}^n} p(x) = 1
\end{cases}
\nonumber
\end{equation}
\end{lemmon}
The proof of Lemma \ref{lemma:rectang_reduction_spatial} is sketched in Appendix section \ref{app:proof_lemma_rectang_reduction_spatial}.

We are, now, in a position to announce the main result in this paper regarding the efficient resolution of the marginal inference problem in general MRF models.
\begin{theorem}
\label{theorem:complexity_upper_bound}
Under Assumptions \ref{assump_computationability}, the resolution for \, $\mathbb{P}\big(\mathcal{X}_s = x_s \big), \, \forall x_s \in \mathcal{L}^{|s|}$ by means
of linear identity (\ref{mrf_dist_solu_lin_sys_polyn}) is attainable with a polynomial time complexity using a linear solver.
\end{theorem}

The proof of Theorem \ref{theorem:complexity_upper_bound} is detailed in Appendix section \ref{app:proof_theorem_complexity_upper_bound}.

Note proof of Theorem \ref{theorem:complexity_upper_bound} also describes a spatial-domain based linear resolution algorithm for the
MRF marginal inference problem via the resolution of linear identity (\ref{mrf_dist_solu_lin_sys_polyn}). Although, most efficient
numerical resolution algorithms of linear identity (\ref{mrf_dist_solu_lin_sys_polyn}) are beyond the scope of this paper,
a more efficient frequency-domain based linear resolution method of linear identity (\ref{mrf_dist_solu_lin_sys_polyn}) is described in appendix
section \ref{subsec:Fourier_based}.

\section{Conclusion}
We have described a polynomial time exact solution to the NP-hard problem of marginal inference in a general Markov random field model,
thus solved concomitantly the {\bf \textit{P = NP}} hypothesis. Interestingly, despite the abysmal combinatorics of the problem, it turns
out that the latter can be squeezed by means of a single ordinary linear identity of polynomial time complexity. Furthermore, the described
approach in this paper remains polynomial time even when the size of the largest neighborhood in a MRF model is not a constant, but scales
logarithmically with the number of nodes of a graph. We would like, nevertheless, to emphasize that a question that this paper does not
address is the numerical stability of the proposed linear solution, at least for well-known instances of the MRF model (such as the
Ising model on a regular lattice) which, thus, will be the subject of future works.



\begin{appendices}

\section{Frequency-domain based MRF marginal inference}
\label{subsec:Fourier_based}

In this Appendix section, we describe a parsimonious reformulation of linear identity (\ref{mrf_dist_solu_lin_sys_polyn}) in the
frequency domain.

First of all, by Lemma \ref{lemma:rectang_reduction_spatial}, linear identity (\ref{mrf_dist_solu_lin_sys_polyn}) can be reexpressed as a system of linear equations as follows:
\begin{equation}
\begin{cases}
{q}_{s\cup {\ne}_{i}^+}(x_{s\cup {\ne}_{i}^+}) - {\pi}_i(x_{{\ne}_{i}^+}) \, \sum_{x_i \in \mathcal{L}} {q}_{s\cup {\ne}_{i}^+}(x_{s\cup {\ne}_{i}^+}) = 0, \,
\forall x_{s\cup {\ne}_{i}^+} \in \mathcal{L}^{|s\cup {\ne}_{i}^+|}, \, \forall i \in \mathcal{V} \\
{q}_{s\cup {\ne}_{i}^+}(x_{s\cup {\ne}_{i}^+}) = \sum_{j \in \mathcal{V}/{s\cup {\ne}_{i}^+}} \sum_{x_j \in \mathcal{L}} q(x_1,\dots,x_j,\dots,x_n), \forall x_{s\cup {\ne}_{i}^+}
\in \mathcal{L}^{|s\cup {\ne}_{i}^+|},  \forall i \in \mathcal{V} \\
\sum_{x \in \mathcal{L}^n} q(x) = 1
\end{cases}
\label{rectang_reduction_spatial_b_four}
\nonumber
\end{equation}
using the Fourier transform, where it has been assumed that $ \forall i \in \mathcal{V}$, ${q}_{s\cup {\ne}_{i}^+}$ stands for the margin
of $q$ with respect to $s\cup {\ne}_{i}^+$.

Thus, denote in the remainder by $\mathbb{C}$ the set of complex numbers, and by $j$ the imaginary unit.
\begin{define}[Discrete Fourier transform]
Suppose an arbitrary function $q:\mathcal{L}^{n} \rightarrow \mathbb{R}$. Then, we call the discrete Fourier transform
of $q$ the complex function $\widehat{q}:\mathcal{L}^{n} \rightarrow \mathbb{C}$ defined as
\begin{equation}
\begin{array}{lll}
\widehat{q}(y) = \sum_{x \in \mathcal{L}^n} {q}(x) \, {\text{e}}^{\frac{-2 j \pi}{L}\langle y,x\rangle}, \, \forall y \in \mathcal{L}^n
\end{array}
\nonumber
\end{equation}
Then, the inversion formula of the Fourier transform of $q$ is given by
\begin{equation}
\begin{array}{lll}
q(x) = \frac{1}{L^n} \, \sum_{y \in \mathcal{L}^n} \widehat{q}(y) \, {\text{e}}^{\frac{2 j \pi}{L}\langle y,x \rangle}, \, \forall x \in \mathcal{L}^n
\end{array}
\label{inverse_four_transf}
\end{equation}
\end{define}

\begin{convention}
Throughout, we convene that
\begin{equation}
\begin{array}{l}
\forall m \in \big\{1,\dots,n\big\}, \, \forall y \in \mathcal{L}^{m}, y = 0 \Leftrightarrow  y_i = 0, \, \forall i \in \big\{ 1,\dots, m \big\}
\end{array}
\nonumber
\end{equation}
\end{convention}

\begin{notation}
Suppose $c$ is a hypervertex of $\mathcal{G}$, and denote throughout
\begin{equation}
\forall y \in \mathcal{L}^n, \,
\begin{cases}
y_{c} = (y_k)_{k \in c}  \\
y_{\mathcal{V}/c} = (y_k)_{k \in \mathcal{V}/c}
\end{cases}
\nonumber
\end{equation}
\end{notation}

Now, suppose an arbitrary function $q:\mathcal{L}^{n} \rightarrow \mathbb{R}$, and denote by $\widehat{q}$ the Fourier transform of $q$,
and by $q_{s\cup {\ne}_{i}^+}$ the margin of $q$ with respect to hypervertex $s\cup {\ne}_{i}^+$, $\forall i\in \mathcal{V}$.
Then, $\forall i \in \mathcal{V}$, by marginalizing both hand sides of formula (\ref{inverse_four_transf}) over $\mathcal{V}/(s\cup {\ne}_{i}^+)$,
one has $\forall x_{s\cup {\ne}_{i}^+} \in \mathcal{L}^{|s\cup {\ne}_{i}^+|}$
\begin{equation}
\begin{array}{lll}
q_{s\cup {\ne}_{i}^+}(x_{s\cup {\ne}_{i}^+}) &=& \frac{1}{L^n} \sum_{k \in \mathcal{V}/(s\cup {\ne}_{i}^+)} \sum_{x_k \in \mathcal{L}}
 \sum_{y \in \mathcal{L}^n} \widehat{q}(y) \, {\text{e}}^{\frac{2 j \pi}{L}\langle y,x\rangle } \\
&=& \frac{1}{L^n} \sum_{y \in \mathcal{L}^n} \widehat{q}(y) \bigg( \sum_{k \in \mathcal{V}/(s\cup {\ne}_{i}^+)} \sum_{x_k \in \mathcal{L}} \, {\text{e}}^{\frac{2 j \pi}{L}\langle y,x\rangle }\bigg)\\
&=& \frac{1}{L^{|s\cup {\ne}_{i}^+|}} \, \sum_{y \in \mathcal{L}^n}  \Big( \widehat{q}(y) \, \delta\big(y_{\mathcal{V}/s\cup {\ne}_{i}^+} \big) \Big) \,
 {\text{e}}^{\frac{2 j\pi}{L} \langle y_{s\cup {\ne}_{i}^+},x_{s\cup {\ne}_{i}^+}\rangle}
\end{array}
\nonumber 
\end{equation}
Also, by marginalization of both hand sides of formula (\ref{inverse_four_transf}) over all sites, one finds
\begin{equation}
\sum_{x \in \mathcal{L}^n} q(x) = \widehat{q}(0,\dots,0)
\nonumber
\end{equation}
Hence, $\forall i \in \mathcal{V}$, denote
\begin{equation}
\begin{array}{l}
\widehat{q}^{(s\cup {\ne}_{i}^+)}(y_{s\cup {\ne}_{i}^+}) =  \widehat{q}(\bar{y}_{s\cup {\ne}_{i}^+}), \, \forall y_{s\cup {\ne}_{i}^+} \in \mathcal{L}^{|s\cup {\ne}_{i}^+|}
\end{array}
\nonumber
\end{equation}
where it has been assumed that $\forall y_{s\cup {\ne}_{i}^+} \in \mathcal{L}^{|s\cup {\ne}_{i}^+|}$, $\bar{y}_{s\cup {\ne}_{i}^+}$ is the $n-$dimensional
vector of which sub-vector over all the indices in $s\cup {\ne}_{i}^+$ coincides with $y_{s\cup {\ne}_{i}^+}$, and of which sub-vector over all the indices
in $\mathcal{V} / (s \cup {\ne}_{i}^+) $ is zero, that is
\begin{equation}
\forall k \in \mathcal{V}, \, \bar{y}_{s\cup {\ne}_{i}^+}(k) =
\begin{cases}
{y}_{s\cup {\ne}_{i}^+}(k), \, \text{if } k \in  s\cup {\ne}_{i}^+, \\
0, \text{ else.}
\end{cases}
\nonumber
\end{equation}
One derives
\begin{equation}
\begin{array}{lll}
q_{s\cup {\ne}_{i}^+}(x_{s\cup {\ne}_{i}^+}) =
\frac{1}{L^{|s\cup {\ne}_{i}^+|}} \, \sum_{y_{s\cup {\ne}_{i}^+} \in \mathcal{L}^{|s\cup {\ne}_{i}^+|}} \widehat{q}^{(s\cup {\ne}_{i}^+)}(y_{s\cup {\ne}_{i}^+}) \, {\text{e}}^{\frac{2 j \pi}{L}\langle y_{s\cup {\ne}_{i}^+},x_{s\cup {\ne}_{i}^+}\rangle }, \, \forall x_{s\cup {\ne}_{i}^+} \in \mathcal{L}^{|s\cup {\ne}_{i}^+|}, \, \forall i \in \mathcal{V}
\end{array}
\nonumber
\end{equation}
Consequently, one can completely reformulate linear system (\ref{rectang_reduction_spatial_b_four}), solely, in terms of the following set
of Fourier coefficients of $q$:
\begin{equation}
\begin{array}{l}
\Big\{  \widehat{q}^{(s\cup {\ne}_{i}^+)}(y_{s\cup {\ne}_{i}^+}), \, \forall y_{s\cup {\ne}_{i}^+} \in \mathcal{L}^{|s\cup {\ne}_{i}^+|}, \, \forall i \in \mathcal{V} \Big\}
\end{array}
\label{four_coef_sparse}
\end{equation}
while bearing in mind, of course, that $\forall i, k \in \mathcal{V}$, such that, $i\neq k$, and
$\forall y_{s\cup {\ne}_{i}^+} \in \mathcal{L}^{|s\cup {\ne}_{i}^+|}$,
$\forall z_{s\cup {\ne}_{k}^+} \in \mathcal{L}^{|s\cup {\ne}_{k}^+|}$, satisfying
\begin{equation}
\forall t \in \mathcal{V},
\begin{cases}
{y}_{s\cup {\ne}_{i}^+}(t) = {z}_{s\cup {\ne}_{k}^+}(t), \, \text{if } t \in  (s\cup {\ne}_{i}^+) \cap (s\cup {\ne}_{k}^+), \\
{y}_{s\cup {\ne}_{i}^+}(t) = 0, \, \text{if } t \in  (s\cup {\ne}_{i}^+) / (s\cup {\ne}_{k}^+),\\
{z}_{s\cup {\ne}_{i}^+}(t) = 0, \, \text{if } t \in  (s\cup {\ne}_{k}^+) / (s\cup {\ne}_{i}^+).
\end{cases}
\nonumber
\end{equation}
one has
\begin{equation}
\begin{array}{l}
\widehat{q}^{(s\cup {\ne}_{i}^+)}(y_{s\cup {\ne}_{i}^+}) =  \widehat{q}^{(s\cup {\ne}_{i}^+)}(z_{s\cup {\ne}_{k}^+})
\end{array}
\label{non_redundant_property_four}
\end{equation}
so as to encode the assumption that the set (\ref{four_coef_sparse}) is the one of the Fourier coefficients of a
same global function $q \in \mathcal{F}$. Nevertheless, such an assumption is most efficiently encoded as follows.
Firstly, one constructs the following set of hypervertices:
\begin{equation}
\mathcal{S} = \cup_{i\in \mathcal{V}} \big\{ c \subseteq s\cup {\ne}_{i}^+ \big\}
\nonumber
\end{equation}
thus consisting of the union of all the partitions (including, obviously the empty hypervertex $\emptyset$) of all the hypervertices $s\cup {\ne}_{i}^+$, $\forall i\in \mathcal{V}$ .
Next, define the following set of Fourier coefficients with respect this time to all the hypervertices in $\mathcal{S}$:
\begin{equation}
\begin{array}{l}
\big\{\widehat{q}^{(c)}(y_c), \, \forall y_c \in \mathcal{L}^{|c|}, \text{ s.t.; } y_c \neq 0, \, \forall c \in \mathcal{S} \big\}
\end{array}
\label{fourier_coefs_closure}
\end{equation}
where it has been assumed that $\widehat{q}^{(\emptyset)} = 1$, in a way which allows to write, $\forall i \in \mathcal{V}$
\begin{equation}
\begin{array}{lll}
q_{s\cup {\ne}_{i}^+}(x_{s\cup {\ne}_{i}^+}) &=& \frac{1}{L^{|s\cup {\ne}_{i}^+|}} \, \underset{c \subseteq s\cup {\ne}_{i}^+}{\sum}, \underset{ y_c \in \mathcal{L}^{|c|}, \text{ s.t.; } y_c \neq 0}{\sum}
\widehat{q}^{(c)}(y_c)\, {\text{e}}^{\frac{2 j \pi\langle y_{c},x_{c}\rangle }{L}}, \, \forall x_{s\cup {\ne}_{i}^+} \in \mathcal{L}^{|s\cup {\ne}_{i}^+|} \\
&=& \frac{1}{L^{|s\cup {\ne}_{i}^+|}}+
\frac{1}{L^{|s\cup {\ne}_{i}^+|}} \, \underset{c \subseteq s\cup {\ne}_{i}^+ \text{ s.t.; } c \neq \emptyset}{\sum}, \underset{ y_c \in \mathcal{L}^{|c|}, \text{ s.t.; } y_c \neq 0}{\sum}
\widehat{q}^{(c)}(y_c)\, {\text{e}}^{\frac{2 j \pi\langle y_{c},x_{c}\rangle }{L}}, \, \forall x_{s\cup {\ne}_{i}^+} \in \mathcal{L}^{|s\cup {\ne}_{i}^+|}
\end{array}
\nonumber
\end{equation}
Consequently, the local consistency constraints between margins of a solution $q \in \mathcal{F}$ linear system (\ref{rectang_reduction_spatial_b_four})
are no longer necessary (as they are already encoded in the Fourier coefficients (\ref{fourier_coefs_closure})), in such a way that, one can
reexpress  linear system (\ref{rectang_reduction_spatial_b_four}), most efficiently, in the Fourier domain as follows:
{\small
\begin{equation}
\begin{array}{lll}
\bigg( 1 + \underset{c \subseteq s\cup {\ne}_{i}^+ \text{ s.t.; }, c \neq \emptyset}{\sum} \, \, \underset{y_c \in \mathcal{L}^{|c|}, \text{ s.t.; } y_c \neq 0}{\sum}
\widehat{q}^{(c)}(y_c)\, {\text{e}}^{\frac{2 j \pi\langle y_{c},x_{c}\rangle }{L}}\, \bigg) \\
\hspace{5mm} -{\pi}_i(x_{{\ne}_{i}^+}) \sum_{x_i\in \mathcal{L}}
 \bigg( 1 + \underset{c \subseteq s\cup {\ne}_{i}^+ \text{ s.t.; } c \neq \emptyset}{\sum} \, \, \underset{y_c \in \mathcal{L}^{|c|}, \text{ s.t.; } y_c \neq 0}{\sum}
\widehat{q}^{(c)}(y_c)\, {\text{e}}^{\frac{2 j \pi\langle y_{c},x_{c}\rangle }{L}}\, \bigg) = 0, \, \forall x_{s\cup {\ne}_{i}^+} \in \mathcal{L}^{|s\cup {\ne}_{i}^+|},  \forall i \in \mathcal{V}
\end{array}
\label{rectang_reduction_spatial_b_hat_b_a}
\end{equation}}
One notes that linear system (\ref{rectang_reduction_spatial_b_hat_b_a}) is sparse, in that, each of its rows only involves a small subset of Fourier coefficients (\ref{fourier_coefs_closure}),
in such a way that, it can be most efficiently solved using a sparse linear solver.

Last, but, not least, given a solution of linear system (\ref{rectang_reduction_spatial_b_hat_b_a}) of the form
\begin{equation}
\widehat{q}^{(c)}(y_c), \, \forall c \in \mathcal{S}, \text{ s.t.;} c \neq \emptyset, \forall y_c \in \mathcal{L}^{|c|}, \text{ s.t.; }  y_c \neq 0
\nonumber
\end{equation}
one writes
\begin{equation}
\mathbb{P}\big(\mathcal{X}_s = x_s \big) =
\frac{1}{L^{|s|}}+ \frac{1}{L^{|s|}} \, \underset{c \subseteq s \text{ s.t.; }, c \neq \emptyset}{\sum} \, \underset{y_c \in \mathcal{L}^{|c|}, \text{ s.t.; } y_c \neq 0}{\sum}
\widehat{q}^{(c)}(y_c)\, {\text{e}}^{\frac{2 j \pi\langle y_{c},x_{c}\rangle}{L}}, \, \forall x_{s} \in \mathcal{L}^{|s|}
\nonumber
\end{equation}


\section{Proofs of Theorems}

\subsection{Proof of Theorem \ref{theo:markov_distrib_unique_pos_definite}}
\label{app:proof_theo_markov_distrib_unique_pos_definite}

\begin{proof}
First of all, linearity of $\mathscr{K}$ is straightforward.
Next, put $p^* \equiv \mathbb{P}\big(\mathcal{X} = \cdot \big)$. One then has
\begin{equation}
\begin{array}{lll}
\big(\mathscr{K} \, p^*\big)(x) &=&  \sum_{x \in \mathcal{L}^n} p^*(x) + \sum_{i \in  \mathcal{V}}  \big( p^*(x) -  \widehat{\pi}_i(x_{{\ne}_{i}^+}) \, \sum_{x_i\in \mathcal{L}} \widehat{\pi}_i(x_{{\ne}_{i}^+}) \, p^*(x) \big) \\
&=& 1 + \sum_{i \in  \mathcal{V}} \big( p^*(x) -  \widehat{\pi}_i(x_{{\ne}_{i}^+}) \, \sum_{x_i\in \mathcal{L}} \widehat{\pi}_i(x_{{\ne}_{i}^+}) \, p^*(x) \big), \, \forall x \in \mathcal{L}^n
\end{array}
\label{kpstar1_thm_a}
\end{equation}
Since by Lemma \ref{lemma:unicity_mrf}, $p^*$ satisfies
\begin{equation}
\begin{array}{lll}
p^*(x) -  \widehat{\pi}_i(x_{{\ne}_{i}^+}) \, \sum_{x_i\in \mathcal{L}} \widehat{\pi}_i(x_{{\ne}_{i}^+}) \, p^*(x) = 0, \, \forall x \in \mathcal{L}^n,\, \forall i \in \mathcal{V}
\end{array}
\nonumber
\end{equation}
one finds that $p^*$ also satisfies
\begin{equation}
\begin{array}{lll}
\sum_{i \in  \mathcal{V}} \big( p^*(x) -  \widehat{\pi}_i(x_{{\ne}_{i}^+}) \, \sum_{x_i\in \mathcal{L}} \widehat{\pi}_i(x_{{\ne}_{i}^+}) \, p^*(x) \big) = 0,  \, \forall x \in \mathcal{L}^n
\end{array}
\nonumber
\end{equation}
in such a way that formula (\ref{kpstar1_thm_a}) becomes
\begin{equation}
\begin{array}{lll}
\big(\mathscr{K} \, p^*\big)(x) = 1, \, \forall x \in \mathcal{L}^n
\end{array}
\nonumber
\end{equation}
which we can reexpress as follows:
\begin{equation}
\begin{array}{lll}
\mathscr{K} \, p^* \equiv \oneindic
\end{array}
\label{kpstar1_thm_b}
\end{equation}
Secondly, the proof that $\mathbb{P}\big(\mathcal{X} = \cdot \big)$ is the unique solution of linear identity (\ref{kpstar1_thm_b})
follows immediately from Lemma \ref{lemma:markov_lin_oper} establishing that $\mathscr{K}$ is an invertible linear operator, thus
proving Theorem \ref{theo:markov_distrib_unique_pos_definite}.
\end{proof}

\subsection{Proof of Theorem \ref{theorem:founding_kp_eq_1}}
\label{app:proof_theorem_founding_kp_eq_1}

\begin{proof}
For simplicity, denote in the remainder $p^* = \mathbb{P}\big(\mathcal{X} = \cdot \big)$, and consider following system of two linear identities:
\begin{equation}
\begin{cases}
\hat{\mathscr{K}} \, p  \equiv \oneindic \\
\widebar{\mathscr{K}} \,p  \equiv \zeroindic
\end{cases}
\label{id_kp_1_kpb_zero_cor_p_0}
\end{equation}
Recall that, by Corollary \ref{corol:p_sol_kp_1_bkp0}, $p^*$ is the unique solution of system (\ref{id_kp_1_kpb_zero_cor_p_0}).
By pre-multiplying the first identity of system (\ref{id_kp_1_kpb_zero_cor_p_0}) by $\hat{\mathscr{K}}^{\dag} $, one reexpresses
system (\ref{id_kp_1_kpb_zero_cor_p_0}) as follows:
\begin{equation}
\begin{cases}
\hat{\mathscr{K}}^{\dag} \, \hat{\mathscr{K}} \, \, p  \equiv \hat{\mathscr{K}}^{\dag} \, \oneindic \\
\widebar{\mathscr{K}} \, p  \equiv \zeroindic \\
p \in \mathscr{F}
\end{cases}
\label{id_kp_1_kpb_zero_cor_p}
\end{equation}
Next, since, by Proposition \ref{prop:rs_eq_rosb}, one has
\begin{equation}
\widebar{\mathscr{K}} \, p \equiv \widebar{\mathscr{K}} \,  \big(\widebar{\mathscr{O}}_s \, p \big), \, \forall p \in \mathscr{F}
\end{equation}
one can reformulate system (\ref{id_kp_1_kpb_zero_cor_p}) as follows:
\begin{equation}
\begin{cases}
\hat{\mathscr{K}}^{\dag} \, \hat{\mathscr{K}} \, \, p  \equiv \hat{\mathscr{K}}^{\dag} \, \oneindic \\
\widebar{\mathscr{K}} \, \big(\widebar{\mathscr{O}}_s \, p\big)  \equiv \zeroindic \\
p \in \mathscr{F}
\end{cases}
\label{id_kp_1_kpb_zero_cor}
\end{equation}
Now, suppose $p$ to be a solution of system (\ref{id_kp_1_kpb_zero_cor}), and let us decompose $p$ as follows:
\begin{equation}
p = {u}_s +\widebar{w}_s
\label{p_decompose_us_wb}
\end{equation}
where $u_s \in \mathscr{F}_s$ and $\widebar{w}_s \in \widebar{\mathscr{F}}_s$, in such a way that, system (\ref{id_kp_1_kpb_zero_cor}) can be reformulated in terms of
the couple $(u_s, \widebar{w}_s)$ as follows:
\begin{equation}
\begin{cases}
\hat{\mathscr{K}}^{\dag} \, \hat{\mathscr{K}} \, u_s +  \hat{\mathscr{K}}^{\dag} \, \hat{\mathscr{K}} \, \widebar{w}_s  \equiv \hat{\mathscr{K}}^{\dag} \, \oneindic \\
\widebar{\mathscr{K}} \, \widebar{w}_s \equiv \zeroindic \\
u_s \in \mathscr{F}_s, \, \widebar{w}_s \in \widebar{\mathscr{F}}_s
\end{cases}
\label{id_kp_1_kpb_zero_cor_u_w}
\end{equation}
Equivalence between systems  (\ref{id_kp_1_kpb_zero_cor}) and (\ref{id_kp_1_kpb_zero_cor_u_w}) is then immediately established using Proposition \ref{prop:F_direct_sum_marginal}
as follows: $p$ is a solution of system (\ref{id_kp_1_kpb_zero_cor}), if and only if, $\big({\mathscr{O}}_s \, p, \widebar{\mathscr{O}}_s \, p \big)$ is
a solution of system (\ref{id_kp_1_kpb_zero_cor_u_w}). Furthermore, again by Proposition \ref{prop:rs_eq_rosb}, system (\ref{id_kp_1_kpb_zero_cor_u_w})
can be reformulated as follows:
\begin{equation}
\begin{cases}
\hat{\mathscr{K}}^{\dag} \, \hat{\mathscr{K}} \, u_s +  \hat{\mathscr{K}}^{\dag} \, \hat{\mathscr{K}} \, \big( \widebar{\mathscr{O}}_s \, {w} \big)  \equiv \hat{\mathscr{K}}^{\dag} \, \oneindic \\
\widebar{\mathscr{K}} \, {w} \equiv \zeroindic \\
u_s \in \mathscr{F}_s, \, w \in \mathscr{F}
\end{cases}
\label{id_kp_1_kpb_zero_cor_u_w_a}
\end{equation}
in the sense that, $(u_s,w)$ is a solution of system (\ref{id_kp_1_kpb_zero_cor_u_w_a}), if and only if, $(u_s, \widebar{\mathscr{O}}_s \, w)$ is a solution of system
(\ref{id_kp_1_kpb_zero_cor_u_w}).

In fact, we will show that the second identity of system (\ref{id_kp_1_kpb_zero_cor_u_w}) does not add to the first identity of the former
any information (constraints) allowing to uniquely identity the quantity $\hat{\mathscr{K}}^{\dag} \, \hat{\mathscr{K}} \, \big( \widebar{\mathscr{O}}_s \, {w} \big)$.

Thus, consider the second identity of system (\ref{id_kp_1_kpb_zero_cor_u_w}) as follows:
\begin{equation}
\widebar{\mathscr{K}} \, w \equiv \zeroindic, \, w \in \mathscr{F}
\label{id_kp_1_kpb_zero_cor_w_a}
\end{equation}
and introduce the following two complementary orthogonal subspaces of $\mathscr{F}$:
\begin{equation}
\begin{cases}
\mathcal{Q} = \big\{ f \equiv \hat{\mathscr{K}}^{\dag} \, \hat{\mathscr{K}} \, u, \, u \in \mathscr{F} \big\} \\
\widebar{\mathcal{Q}} = \big\{ f \equiv \big(\idoper - \hat{\mathscr{K}}^{\dag} \, \hat{\mathscr{K}} \big) \, u, \, u \in \mathscr{F} \big\}
\nonumber
\end{cases}
\end{equation}
In addition, suppose $w$ to be a solution of identity (\ref{id_kp_1_kpb_zero_cor_w_a}), and decompose $w$ as the sum of an element of $\mathcal{Q}$ and an element
of $\widebar{\mathcal{Q}}$, respectively, as follows:
\begin{equation}
w \equiv \hat{\mathscr{K}}^{\dag} \, \hat{\mathscr{K}} \, w + \big(\idoper - \hat{\mathscr{K}}^{\dag} \, \hat{\mathscr{K}} \big) \, w
\end{equation}
One reformulates accordingly identity (\ref{id_kp_1_kpb_zero_cor_w_a}) as follows:
\begin{equation}
\begin{cases}
\widebar{\mathscr{K}} \,  \big(\hat{\mathscr{K}}^{\dag} \, \hat{\mathscr{K}} \, w\big) +
\widebar{\mathscr{K}} \,  \Big( \big(\idoper - \hat{\mathscr{K}}^{\dag} \, \hat{\mathscr{K}} \big) \, w \Big) \equiv \zeroindic\\
w \in \mathscr{F}
\end{cases}
\label{id_kp_1_kpb_zero_cor_w_b}
\end{equation}
and by using Lemma \ref{lemma:kq_eq_kbar_q}, one finds that identity (\ref{id_kp_1_kpb_zero_cor_w_b}) is equivalent to the following identity:
\begin{equation}
\begin{cases}
\widebar{\mathscr{K}} \,  \big(\hat{\mathscr{K}}^{\dag} \, \hat{\mathscr{K}} \, w\big) +
{\mathscr{K}} \,  \Big( \big(\idoper - \hat{\mathscr{K}}^{\dag} \, \hat{\mathscr{K}} \big) \, w \Big) \equiv \zeroindic\\
w \in \mathscr{F}
\end{cases}
\label{id_kp_1_kpb_zero_cor_w_c}
\end{equation}
last, since ${\mathscr{K}}$ is an invertible linear operator, by pre-multiplying both hand sides of identity (\ref{id_kp_1_kpb_zero_cor_w_c}) by ${\mathscr{K}}^{-1}$,
one finds that identity (\ref{id_kp_1_kpb_zero_cor_w_c}) is, in turn, equivalent to the following identity:
\begin{equation}
\begin{cases}
\big(\idoper - \hat{\mathscr{K}}^{\dag} \, \hat{\mathscr{K}} \big) \, w \equiv
-{\mathscr{K}}^{-1} \, \widebar{\mathscr{K}} \,  \big(\hat{\mathscr{K}}^{\dag} \, \hat{\mathscr{K}} \, w\big) \\
w \in \mathscr{F}
\end{cases}
\label{id_kp_1_kpb_zero_cor_w_d}
\end{equation}
Consequently, one may reexpress system (\ref{id_kp_1_kpb_zero_cor_u_w_a}) as follows:
\begin{equation}
\begin{cases}
\hat{\mathscr{K}}^{\dag} \, \hat{\mathscr{K}} \, u_s +  \hat{\mathscr{K}}^{\dag} \, \hat{\mathscr{K}} \,  \widebar{\mathscr{O}}_s \, {w}   \equiv \hat{\mathscr{K}}^{\dag} \, \oneindic \\
\big(\idoper - \hat{\mathscr{K}}^{\dag} \, \hat{\mathscr{K}} \big) \, w \equiv
-{\mathscr{K}}^{-1} \, \widebar{\mathscr{K}} \,  \big(\hat{\mathscr{K}}^{\dag} \, \hat{\mathscr{K}} \, w\big) \\
u_s \in \mathscr{F}_s, \, w \in {\mathscr{F}}
\end{cases}
\label{id_kp_1_kpb_zero_cor_u_w_dd}
\end{equation}
Now, suppose $w \in \mathscr{F}$. Since, $\hat{\mathscr{K}}^{\dag} \, \hat{\mathscr{K}} \, w$ and $\big(\idoper - \hat{\mathscr{K}}^{\dag} \, \hat{\mathscr{K}} \big) \, w$
stand for two independent general elements of $\mathcal{Q}$ and $\widebar{\mathcal{Q}}$, respectively, one concludes that
the second identity of system (\ref{id_kp_1_kpb_zero_cor_u_w_dd}) holds, equivalently, the second identity of system (\ref{id_kp_1_kpb_zero_cor_u_w_a}) holds,
if and only if, ${\mathscr{K}}^{-1} \, \widebar{\mathscr{K}} \,  \big(\hat{\mathscr{K}}^{\dag} \, \hat{\mathscr{K}} \, w\big)  \in \widebar{\mathcal{Q}}$, said otherwise,
if and only if, one has
\begin{equation}
\hat{\mathscr{K}} \, {\mathscr{K}}^{-1} \, \widebar{\mathscr{K}} \,  \big(\hat{\mathscr{K}}^{\dag} \, \hat{\mathscr{K}} \, w\big)  \equiv \zeroindic
\label{id_kp_1_kpb_zero_cor_u_w_b_a}
\end{equation}
This is on the one hand. On the other hand, by Lemma \ref{lemma:i_kk_osq_eq_func_kkosq}, one has
\begin{equation}
\hat{\mathscr{K}} \, {\mathscr{K}}^{-1} \, \widebar{\mathscr{K}} \,  \Big(\hat{\mathscr{K}}^{\dag} \, \hat{\mathscr{K}} \, \big( {\mathscr{O}}_s \, {w}  \big) \Big)  \equiv \zeroindic
\label{id_kp_1_kpb_zero_cor_u_w_b_b}
\end{equation}
in such a way that, by subtracting formula (\ref{id_kp_1_kpb_zero_cor_u_w_b_b}) from formula (\ref{id_kp_1_kpb_zero_cor_u_w_b_a}),  one finds that the second identity of system (\ref{id_kp_1_kpb_zero_cor_u_w_a}) holds, if and only if, one has
\begin{equation}
\hat{\mathscr{K}} \, {\mathscr{K}}^{-1} \, \widebar{\mathscr{K}} \,  \big(\hat{\mathscr{K}}^{\dag} \, \hat{\mathscr{K}} \, \widebar{\mathscr{O}}_s \,  w\big)  \equiv \zeroindic
\label{constraint_w}
\end{equation}
Next, by first observing that the first identity of system (\ref{id_kp_1_kpb_zero_cor_u_w_dd}) depends on $w$  only through $\hat{\mathscr{K}}^{\dag} \, \hat{\mathscr{K}} \, \widebar{\mathscr{O}}_s \,  w$, second, by injecting into system (\ref{id_kp_1_kpb_zero_cor_u_w_dd}) the constraint (\ref{constraint_w}) after ignoring the term $\big(\idoper - \hat{\mathscr{K}}^{\dag} \, \hat{\mathscr{K}} \big) \, w$, thereby also pruning from system (\ref{id_kp_1_kpb_zero_cor_u_w_a}) its second identity (as the latter necessarily holds provided that constraint (\ref{constraint_w}) is satisfied), one concludes that both terms $\hat{\mathscr{K}}^{\dag} \, \hat{\mathscr{K}} \, u_s$ and $\hat{\mathscr{K}}^{\dag} \, \hat{\mathscr{K}} \, \widebar{\mathscr{O}}_s \,  w$ in
system (\ref{id_kp_1_kpb_zero_cor_u_w_dd}) must be completely identified by the following system:
\begin{equation}
\begin{cases}
\hat{\mathscr{K}}^{\dag} \, \hat{\mathscr{K}} \, u_s +  \hat{\mathscr{K}}^{\dag} \, \hat{\mathscr{K}} \,  \widebar{\mathscr{O}}_s \, {w}   \equiv \hat{\mathscr{K}}^{\dag} \, \oneindic \\
\hat{\mathscr{K}} \, {\mathscr{K}}^{-1} \, \widebar{\mathscr{K}} \,  \big(\hat{\mathscr{K}}^{\dag} \, \hat{\mathscr{K}} \, \widebar{\mathscr{O}}_s \,  w\big)  \equiv \zeroindic \\
u_s \in \mathscr{F}_s, \, w \in {\mathscr{F}}
\end{cases}
\end{equation}
equivalently, both terms $\hat{\mathscr{K}}^{\dag} \, \hat{\mathscr{K}} \, u_s$ and $\hat{\mathscr{K}}^{\dag} \, \hat{\mathscr{K}} \, \widebar{\mathscr{O}}_s \,  w$ in
system (\ref{id_kp_1_kpb_zero_cor_u_w}) must be completely identified by the following system:
\begin{equation}
\begin{cases}
\hat{\mathscr{K}}^{\dag} \, \hat{\mathscr{K}} \, u_s +  \hat{\mathscr{K}}^{\dag} \, \hat{\mathscr{K}} \,  \widebar{w}_s   \equiv \hat{\mathscr{K}}^{\dag} \, \oneindic \\
\hat{\mathscr{K}} \, {\mathscr{K}}^{-1} \, \widebar{\mathscr{K}} \,  \big(\hat{\mathscr{K}}^{\dag} \, \hat{\mathscr{K}} \, \widebar{w}_s \big)  \equiv \zeroindic \\
u_s \in \mathscr{F}_s, \, \widebar{w}_s \in \widebar{\mathscr{F}}
\end{cases}
\label{id_kp_1_kpb_zero_cor_u_w_c}
\end{equation}
Hence, it remains to show that the second identity of system (\ref{id_kp_1_kpb_zero_cor_u_w_c}) is redundant with its first one.

Thus, suppose $(u_s, \widebar{w}_s) \in \mathscr{F}_s, \times \widebar{\mathscr{F}}$ to be a solution to the first identity of system (\ref{id_kp_1_kpb_zero_cor_u_w_c}),
that is,  $(u_s, \widebar{w}_s)$ is a solution of the following identity:
\begin{equation}
\hat{\mathscr{K}}^{\dag} \, \hat{\mathscr{K}} \, u_s +  \hat{\mathscr{K}}^{\dag} \, \hat{\mathscr{K}} \,  \widebar{w}_s   \equiv \hat{\mathscr{K}}^{\dag} \, \oneindic
\label{id_kp_1_kpb_zero_cor_u_w_c_a}
\end{equation}
Then, by using the identity $\hat{\mathscr{K}}^{\dag} \, \oneindic  \equiv \hat{\mathscr{K}}^{\dag} \, \hat{\mathscr{K}} \, p^*$, and by pre-multiplying both hand sides of identity (\ref{id_kp_1_kpb_zero_cor_u_w_c_a}) by $\hat{\mathscr{K}} \, \widebar{\mathscr{K}}^{-1} \, \widebar{\mathscr{K}}$, one finds
\begin{equation}
\begin{array}{lll}
\hat{\mathscr{K}} \, {\mathscr{K}}^{-1} \, \widebar{\mathscr{K}} \, \hat{\mathscr{K}}^{\dag} \, \hat{\mathscr{K}} \,  \widebar{w}_s  \equiv
\hat{\mathscr{K}} \, {\mathscr{K}}^{-1} \, \widebar{\mathscr{K}} \, \hat{\mathscr{K}}^{\dag} \, \hat{\mathscr{K}} \, p^* - \hat{\mathscr{K}} \, {\mathscr{K}}^{-1} \, \widebar{\mathscr{K}} \, \hat{\mathscr{K}}^{\dag} \, \hat{\mathscr{K}} \, u_s
\nonumber
\end{array}
\end{equation}
Furthermore, by Lemma \ref{lemma:i_kk_osq_eq_func_kkosq}, one has
\begin{equation}
\hat{\mathscr{K}} \, {\mathscr{K}}^{-1} \, \widebar{\mathscr{K}} \, \hat{\mathscr{K}}^{\dag} \, \hat{\mathscr{K}} \, u_s  \equiv \zeroindic
\nonumber
\end{equation}
and
\begin{equation}
\begin{array}{lll}
\hat{\mathscr{K}} \, {\mathscr{K}}^{-1} \, \widebar{\mathscr{K}} \, \hat{\mathscr{K}}^{\dag} \, \hat{\mathscr{K}} \, p^* &\equiv&
\hat{\mathscr{K}} \, {\mathscr{K}}^{-1} \, \big({\mathscr{K}} - \hat{\mathscr{K}} \big) \, \hat{\mathscr{K}}^{\dag} \, \hat{\mathscr{K}} \, p^* \\
&\equiv&
\hat{\mathscr{K}} \, \hat{\mathscr{K}}^{\dag} \, \hat{\mathscr{K}} \, p^* -
\hat{\mathscr{K}} \, {\mathscr{K}}^{-1} \, \hat{\mathscr{K}}, \hat{\mathscr{K}} \, \hat{\mathscr{K}}^{\dag} \, \hat{\mathscr{K}} \, p^* \\
&\equiv&
\hat{\mathscr{K}} \, p^* -
\hat{\mathscr{K}} \, {\mathscr{K}}^{-1} , \hat{\mathscr{K}} \, p^*\\
&\equiv&
\oneindic  - \hat{\mathscr{K}} \, \big( {\mathscr{K}}^{-1} , \oneindic \big)\\
&\equiv&
\oneindic  - \hat{\mathscr{K}} \, p^*\\
&\equiv&
\oneindic  - \oneindic \\
&\equiv&
\zeroindic
\end{array}
\nonumber
\end{equation}
implying that
\begin{equation}
\begin{array}{lll}
\hat{\mathscr{K}} \, {\mathscr{K}}^{-1} \, \widebar{\mathscr{K}} \, \hat{\mathscr{K}}^{\dag} \, \hat{\mathscr{K}} \,  \widebar{w}_s   \equiv
\zeroindic
\nonumber
\end{array}
\end{equation}
meaning that the second identity of system (\ref{id_kp_1_kpb_zero_cor_u_w_c}) is redundant with its first identity, in such a way that,
system (\ref{id_kp_1_kpb_zero_cor_u_w_c}) boils down to the following single identity:
\begin{equation}
\begin{cases}
\hat{\mathscr{K}}^{\dag} \, \hat{\mathscr{K}} \, u_s +  \hat{\mathscr{K}}^{\dag} \, \hat{\mathscr{K}} \, \widebar{w}_s  \equiv \hat{\mathscr{K}}^{\dag} \, \oneindic \\
u_s \in \mathscr{F}_s, \, \widebar{w}_s \in \widebar{\mathscr{F}}
\end{cases}
\label{id_kp_1_kpb_zero_cor_u_w_e}
\end{equation}
Thus, by using formula (\ref{p_decompose_us_wb}), and the fact that $p^*$ is a solution of system (\ref{id_kp_1_kpb_zero_cor_p_0}),
one concludes that the quantities $\hat{\mathscr{K}}^{\dag} \, \hat{\mathscr{K}} \, \big({\mathscr{O}}_s \, p^*\big)$ and
$\hat{\mathscr{K}}^{\dag} \, \hat{\mathscr{K}} \, \big(\widebar{\mathscr{O}}_s \, p^*\big)$
are totally identified by
identity (\ref{id_kp_1_kpb_zero_cor_u_w_e}), in that, whatever a solution $(u_s,\widebar{w}_s)$ of identity (\ref{id_kp_1_kpb_zero_cor_u_w_e}), one must have
\begin{equation}
\begin{cases}
\hat{\mathscr{K}}^{\dag} \, \hat{\mathscr{K}} \, u_s \equiv \hat{\mathscr{K}}^{\dag} \, \hat{\mathscr{K}} \, \big({\mathscr{O}}_s \, p^*\big) \\
\hat{\mathscr{K}}^{\dag} \, \hat{\mathscr{K}} \, \widebar{w}_s  \equiv \hat{\mathscr{K}}^{\dag} \, \hat{\mathscr{K}} \, \big(\widebar{\mathscr{O}}_s \, p^*\big)
\end{cases}
\nonumber
\end{equation}
In particular, whatever a solution $(u_s,\widebar{w}_s)$ of identity (\ref{id_kp_1_kpb_zero_cor_u_w_e}), one must have
\begin{equation}
\hat{\mathscr{K}}^{\dag} \, \hat{\mathscr{K}} \, u_s \equiv \hat{\mathscr{K}}^{\dag} \, \hat{\mathscr{K}} \, \big({\mathscr{O}}_s \, p^*\big)
\label{id_kp_1_kpb_zero_cor_u_w_f}
\end{equation}
Furthermore, by pre-multiplying both hand sides of identity (\ref{id_kp_1_kpb_zero_cor_u_w_f}) by $\hat{\mathscr{K}}$, one reexpresses identity (\ref{id_kp_1_kpb_zero_cor_u_w_f}) as follows:
\begin{equation}
\hat{\mathscr{K}} \, u_s \equiv \hat{\mathscr{K}} \, \big({\mathscr{O}}_s \, p^*\big)
\nonumber
\end{equation}
which, by Proposition \ref{prop:ks_eq_kos}, is equivalent to the following identity:
\begin{equation}
{\mathscr{K}} \, u_s \equiv {\mathscr{K}} \, \big({\mathscr{O}}_s \, p^*\big)
\nonumber
\end{equation}
which, in turn, since ${\mathscr{K}}$ is an invertible linear operator, is equivalent to the following identity:
\begin{equation}
u_s \equiv {\mathscr{O}}_s \, p^*
\nonumber
\end{equation}
One concludes that whatever a solution $(u_s,\widebar{w}_s)$ of identity (\ref{id_kp_1_kpb_zero_cor_u_w_e}), one must have
\begin{equation}
u_s \equiv {\mathscr{O}}_s \, p^*
\nonumber
\end{equation}
To complete the proof of Proposition, one pre-multiplies both hand sides of identity (\ref{id_kp_1_kpb_zero_cor_u_w_e}) by $\hat{\mathscr{K}}^{\dag}$,
and uses formula (\ref{p_decompose_us_wb}) to retrieve the first identity of system (\ref{id_kp_1_kpb_zero_cor_p_0}) to finally establish that whatever
a solution $p$ of the following identity:
\begin{equation}
\begin{cases}
\hat{\mathscr{K}} \, p \equiv \oneindic \\
p \in \mathscr{F}
\end{cases}
\nonumber
\end{equation}
one must have
\begin{equation}
{\mathscr{O}}_s \, p \equiv {\mathscr{O}}_s \, p^*
\nonumber
\end{equation}
which amounts to saying that the margin of $p$ with respect to $s$ and the margin of $p^*$ with respect to $s$ must coincide,
thus establishing the proof of Theorem \ref{theorem:founding_kp_eq_1}.
\end{proof}

\section{Proof of Theorem \ref{theorem:complexity_upper_bound}}
\label{app:proof_theorem_complexity_upper_bound}

\begin{proof}
Suppose $p \in \mathscr{F}$. Hence, by Lemma \ref{lemma:rectang_reduction_spatial}, one finds that $p$ is a solution of linear identity (\ref{mrf_dist_solu_lin_sys_polyn}),
if and only if, $p$ satisfies the following (over-determined) ordinary system of linear equations:
\begin{equation}
\begin{cases}
{p}_{s\cup {\ne}_{i}^+}(x_{s\cup {\ne}_{i}^+}) - {\pi}_i(x_{{\ne}_{i}^+}) \, \sum_{x_i \in \mathcal{L}} {p}_{s\cup {\ne}_{i}^+}(x_{s\cup {\ne}_{i}^+}) = 0, \,
\forall x_{s\cup {\ne}_{i}^+} \in \mathcal{L}^{|s\cup {\ne}_{i}^+|}, \, \forall i \in \mathcal{V} \\
{p}_{s\cup {\ne}_{i}^+}(x_{s\cup {\ne}_{i}^+}) = \sum_{j \in \mathcal{V}/{s\cup {\ne}_{i}^+}} \sum_{x_j \in \mathcal{L}} p(x_1,\dots,x_j,\dots,x_n), \forall x_{s\cup {\ne}_{i}^+}
\in \mathcal{L}^{|s\cup {\ne}_{i}^+|}, \, \forall i \in \mathcal{V} \\
\sum_{x \in \mathcal{L}^n} p(x) = 1
\end{cases}
\label{rectang_reduction_spatial_b_a}
\end{equation}
where the following set of equations:
\begin{equation}
{p}_{s\cup {\ne}_{i}^+}(x_{s\cup {\ne}_{i}^+}) = \sum_{j \in \mathcal{V}/{s\cup {\ne}_{i}^+}} \sum_{x_j \in \mathcal{L}} p(x_1,\dots,x_j,\dots,x_n), \forall x_{s\cup {\ne}_{i}^+}
\in \mathcal{L}^{|s\cup {\ne}_{i}^+|},  \forall i \in \mathcal{V}
\label{const:marginals}
\end{equation}
simply expresses that for all $i \in  \mathcal{V}$, ${p}_{s\cup {\ne}_{i}^+}$ is the margin of $p$ with respect to $s\cup {\ne}_{i}^+$. Next,
Lemma \ref{lemma:local_consistency} guarantees that one can enforce constraints (\ref{const:marginals}) by simply enforcing the local
consistency constraints between $n$ local functions ${q}^{(i)}: \mathcal{L}^{|s\cup {\ne}_{i}^+|} \rightarrow \mathbf{R}$, $\forall i \in \mathcal{V}$
as follows:
{\small
\begin{equation}
\begin{cases}
{q}^{(i)}(x_{s\cup {\ne}_{i}^+}) - {\pi}_i(x_{{\ne}_{i}^+}) \, \sum_{x_i \in \mathcal{L}} q^{(i)}(x_{s\cup {\ne}_{i}^+}) = 0, \,
\forall x_{s\cup {\ne}_{i}^+} \in \mathcal{L}^{|s\cup {\ne}_{i}^+|}, \, \forall i \in \mathcal{V} \\
\sum_{ j \in { {\ne}_{i}^+} / (s\cup {\ne}_{k}^+)} \sum_{x_j \in \mathcal{L}} \, {q}^{(i)}(x_{s\cup {\ne}_{i}^+}) -
\sum_{ j \in {{\ne}_{k}^+} /  (s\cup {\ne}_{i}^+)}  \sum_{x_j \in \mathcal{L}} \, {q}^{(k)}(x_{s\cup {\ne}_{i}^+}) = 0, \,  \\
\hspace{4cm} \forall x_{s \cup ( {\ne}_{i}^+ \cap {\ne}_{k}^+)} \in \mathcal{L}^{\big|s \cup ( {\ne}_{i}^+ \cap {\ne}_{k}^+ \big|},  \forall i, k \in \mathcal{V} \text{ s.t.; } i \neq k \\
\sum_{x_{s\cup {\ne}_{i}^+} \in \mathcal{L}^{|s\cup {\ne}_{i}^+|}} {q}^{(i)}(x_{s\cup {\ne}_{i}^+}) = 1, \, \forall i \in \mathcal{V}
\end{cases}
\label{rectang_reduction_spatial_b_b}
\end{equation}}
in that, whatever a set of $n$ local functions $ \Big\{ {q}^{(i)}: \mathcal{L}^{|s\cup {\ne}_{i}^+|} \rightarrow \mathbf{R}, \, \forall i \in \mathcal{V} \Big\}$
satisfying constraints (\ref{rectang_reduction_spatial_b_b}), there exists a function $q: \mathcal{L}^{n} \rightarrow \mathbf{R}$ satisfying:
\begin{equation}
{p}_{s\cup {\ne}_{i}^+}(x_{s\cup {\ne}_{i}^+}) = {q}^{(i)}(x_{s\cup {\ne}_{i}^+}), \, \forall x_{s\cup {\ne}_{i}^+} \in \mathcal{L}^{|s\cup {\ne}_{i}^+|}, \, \forall i \in \mathcal{V}
\nonumber
\end{equation}
where $\forall i \in \mathcal{V}$, ${p}_{s\cup {\ne}_{i}^+}$ stands for the margin of $p$ with respect to hypervertex $s\cup {\ne}_{i}^+$.
Hence, one also has
\begin{equation}
{p}_{s}(x_s) = \sum_{j \in {{\ne}_{i}^+}/s}   {q}^{(i)}(x_{s\cup {\ne}_{i}^+}), \,  \forall x_{s} \in \mathcal{L}^{|s|}, \, \forall i  \in \mathcal{V}
\nonumber
\end{equation}
in such a way that, Theorem \ref{theorem:founding_kp_eq_1} guarantees that
\begin{equation}
\begin{array}{lll}
\mathbb{P}\big(\mathcal{X}_s = x_s \big) = {p}_{s}(x_s), \,  \forall x_{s} \in \mathcal{L}^{|s|}
\end{array}
\nonumber 
\end{equation}
One can reexpress linear system (\ref{rectang_reduction_spatial_b_b}) in matrix form as
\begin{equation}
\mathcal{A} \, y = b
\label{matrix_form_main_identity}
\end{equation}
where,
\begin{itemize}
\item $y$ stands for the following real vector of unknowns:
\begin{equation}
\begin{array}{lll}
y = \Big( {q}^{(i)}(x_{s\cup {\ne}_{i}^+}) \Big)_{ x_{s\cup {\ne}_{i}^+} \in \mathcal{L}^{|s\cup {\ne}_{i}^+|}, \, i \in \mathcal{V}}
\end{array}
\nonumber
\end{equation}
of size
\begin{equation}
\begin{array}{lll}
\eta = \sum_{i\in \mathcal{V}} L^{|s\cup {\ne}_{i}^+|}
\end{array}
\nonumber 
\end{equation}
\item $b$ stands for a known real vector of zeros and ones of size
\begin{equation}
\begin{array}{lll}
\mu = \sum_{i\in \mathcal{V}} L^{|s\cup {\ne}_{i}^+|} + \sum_{i,j \in \mathcal{V}, \text{ s.t., } i \neq j} L^{|s\cup ({\ne}_{i}^+ \cap {\ne}_{j}^+) |} + |\mathcal{V}|
\end{array}
\nonumber 
\end{equation}
\item $\mathcal{A}$ stands for a real $2D$ matrix of size  $(\eta,\nu)$
\end{itemize}
Observe that system (\ref{matrix_form_main_identity}) is over-determined, in that, its number of constraints is larger than its number of
unknowns. Consequently, instead of solving for any solution of $q$ of system (\ref{matrix_form_main_identity}), one may rather solve
the following linear system:
\begin{equation}
\mathcal{A}^{T} \, \mathcal{A} \, y = \mathcal{A}^{T} \,b
\label{matrix_form_main_identity_b}
\end{equation}
where $\mathcal{A}^{T}$ stands for the transpose of $\mathcal{A}$. The complexity of the transformation of linear
system (\ref{matrix_form_main_identity}) into linear system (\ref{matrix_form_main_identity_b}) is $C_1 = O(\eta^2 \, \mu)$, in that, it consists of
\begin{itemize}
\item a matrix-matrix multiplication for computing $\mathcal{A}^{T} \, \mathcal{A}$ of which complexity is $O(\eta^2 \, \mu)$
\item a matrix-vector multiplication for computing $\mathcal{A}^{T} \, b$ of which complexity is $O(\eta \, \mu)$
\end{itemize}
observe that $\mathcal{A}^{T} \, \mathcal{A}$ is a $\eta \times \eta$ matrix. Thus, one one may finally solve for $y$
in system (\ref{matrix_form_main_identity_b}) using a standard linear solver, for instance,
the "Gaussian elimination method" of which time complexity is $O(\eta^{3})$.
One finds that the time complexity of computation of the marginal distribution $\mathbb{P}\big(\mathcal{X}_s = \cdot \big)$
by means of the resolution of linear identity (\ref{mrf_dist_solu_lin_sys_polyn}) is $O(\eta^{3} + \eta^2 \, \mu)$ which we thus
propose to bound by some function of $n$. Indeed, under Assumptions \ref{assump_computationability}, one has
\begin{equation}
\begin{array}{lll}
\eta &\leq& n \, L^{|s|} \, L^{ \sup_{i \in \mathcal{V}} \big\{ \big |{\ne}_{i}^+| \big| \big \} } \\
&\leq& n \, \big(L^{\tau_1 \, \log(n) + \tau_2}\big)^2 \\
&=& L^{2\tau_2} n^{1+2\, \tau_1 \, \log(L)}
\end{array}
\nonumber 
\end{equation}
and
\begin{equation}
\begin{array}{lll}
\mu &\leq& n \, L^{|s|} \, L^{ \sup_{i \in \mathcal{V}} \big\{ \big |{\ne}_{i}^+| \big| \big \} } + n^2 \, L^{|s|} \,  L^{ \sup_{i \in \mathcal{V}} \big\{ \big |{\ne}_{i}^+| \big| \big \} } + n \\
&\leq& 3 n^2 \, L^{|s|} \,  L^{ \sup_{i \in \mathcal{V}} \big\{ \big |{\ne}_{i}^+| \big| \big \} } \\
&\leq& 3 n^2 \, \big(L^{\tau_1 \, \log(n) + \tau_2}\big)^2\\
&=& 3 \, L^{2\tau_2} \, n^{2+2\, \tau_1 \, \log(L)}
\end{array}
\nonumber 
\end{equation}
hence,
\begin{equation}
\begin{array}{lll}
\eta^{3} + \eta^2 \, \mu &\leq& \tau_1  \, n^{3+6\, \tau_1 \, \log(L)} + \tau_2 \,  n^{4+6\, \tau_2 \, \log(L)} \\
&\leq& 2 \tau \, n^{\tau'}
\end{array}
\nonumber 
\end{equation}
where it has been assumed that $\tau = \sup\{\tau_1, \tau_2\}$, and $\tau' = 4+6\, \tau_1 \, \log(L)$ which thus stand for two
real positive constants. One concludes that the time complexity of computation of the marginal distribution $\mathbb{P}\big(\mathcal{X}_s = \cdot \big)$
by means of the resolution of linear identity (\ref{mrf_dist_solu_lin_sys_polyn}) is bounded by a polynomial function of $n$, hence
the proof of Theorem \ref{theorem:complexity_upper_bound}.
\end{proof}


\section{Proofs of Corollaries}

\subsection{Proof of Corollary \ref{corol:p_sol_kp_1_bkp0}}
\label{app:proof_corol_p_sol_kp_1_bkp0}.

\begin{proof}
Firstly, by Propositions \ref{prop:kts_p_sar_eq_one} and \ref{prop:rts_p_sar_eq_zero}, joint MRF distribution $\mathbb{P}\big(\mathcal{X} = \cdot \big)$
is a solution to linear system of identities (\ref{identity_kp_1_kpb_zero}). Secondly, assume $p \in \mathscr{F}$ satisfying (\ref{identity_kp_1_kpb_zero}).
It follows that $p$ also satisfies the following identity:
\begin{equation}
\hat{\mathscr{K}} \,p  + \widebar{\mathscr{K}} \,p  \equiv \oneindic
\nonumber
\end{equation}
equivalently, the following identity:
\begin{equation}
{\mathscr{K}} \,p  \equiv \oneindic
\nonumber
\end{equation}
of which unique solution, by Theorem \ref{theo:markov_distrib_unique_pos_definite}, is joint MRF distribution $\mathbb{P}\big(\mathcal{X} = \cdot \big)$, hence
the proof of Corollary \ref{corol:p_sol_kp_1_bkp0}.
\end{proof}


\section{Proofs of Propositions}

\subsection{Proof of Proposition \ref{prop:o_s_ortho}}
\label{app:proof_prop_o_s_ortho}
\begin{proof}
Assume an arbitrary hypervertex $c$ of graph $\mathcal{G}$, and an arbitrary function $u \in \mathscr{F}$. One writes
\begin{equation}
\begin{array}{lll}
\mathscr{O}_c \, \big( \mathscr{O}_c \, u ) &\equiv& \mathscr{O}_c \, \Big(\frac{u_c}{L^{n-|c|}}\Big) \\
&\equiv& \frac{\mathscr{O}_c \, u_c}{L^{n-|c|}}
\end{array}  \nonumber
\end{equation}
Since $\forall x \in \mathcal{L}^n$, $u_c(x_c)$, solely, depends upon $x_c$, one finds that the margin of $u_c$ with respect to $c$
is given by $L^{n-|c|} \, u_c$, in such a way that,
one can write
\begin{equation}
\begin{array}{lll}
\mathscr{O}_c \, \big( \mathscr{O}_c \, u ) &\equiv& \frac{\frac{L^{n-|c|} \, u_c}{L^{n-|c|}}}{L^{n-|c|}}  \equiv \frac{u_c}{L^{n-|c|}} \\
&\equiv& \mathscr{O}_c \, u
\end{array}  \nonumber
\end{equation}
Since function $u$ is assumed to be arbitrary, one finds
\begin{equation}
\begin{array}{lll}
\mathscr{O}_c \, \big( \mathscr{O}_c \, u ) \equiv \mathscr{O}_c \, u, \, \forall u \in \mathscr{F}
\end{array}  \nonumber
\end{equation}
which means that $\mathscr{O}_c$ is a projection. Next, assume two functions $u,v \in \mathscr{F}$, and denote by $u_c$ and $v_c$ the margins of
$u$ and of $v$ with respect to $c$, respectively. One then writes
\begin{equation}
\begin{array}{lll}
\langle v,\mathscr{O}_c \,u \rangle  &=& \sum_{x \in \mathcal{L}^{n}} v(x) \, \big( \frac{u_c(x_c)}{L^{n-|c|}} \big) \\
                               &=&  \sum_{i \in c} \sum_{x_i \in \mathcal{L}} \Big( u_c(x_c) \, \frac{ \sum_{j \in \mathcal{V}/c} \sum_{x_j \in \mathcal{L}}  v(x)}{L^{n-|c|}}  \Big) \\
                               &=&  \sum_{i \in c} \sum_{x_i \in \mathcal{L}} u_c(x_c) \, \big(\frac{v_c(x_c)}{L^{n-|c|}} \big)  \\
                               &=&  \sum_{x \in \mathcal{L}^{n}} u(x) \, \big( \frac{v_c(x_c)}{L^{n-|c|}} \big) \\
                               &=&  \langle u,\mathscr{O}_c \,v \rangle

\end{array}  \nonumber
\end{equation}
which means that $\mathscr{O}_c$ is a self-adjoint linear operator. To sump  up, we have
shown that $\mathscr{O}_c$ is a self-adjoint projection, hence, an orthogonal projection,
thus establishing the proof of Proposition \ref{prop:o_s_ortho}.
\end{proof}

\subsection{Proof of Proposition \ref{prop:F_direct_sum_marginal}}
\label{app:proof_prop_F_direct_sum_marginal}

\begin{proof}
This is a classical linear algebraic result following immediately from Proposition \ref{prop:o_s_ortho} which establishes
that $\mathscr{O}_s$ is an orthogonal projection, in such a way that, $\hat{\mathscr{F}}_s$ and $\widebar{\mathscr{F}}_s$
are two complementary orthogonal subspaces of $\mathscr{F}$.
\end{proof}

\subsection{Proof of Proposition \ref{prop:ki_is_ortho_project}}
\label{app:proof_prop_ki_is_ortho_project}
\begin{proof}
Suppose $\forall i \in \mathcal{V}$, and introduce linear operators $r_i: \mathscr{F}\rightarrow \mathscr{F}$ defined $\forall u \in \mathscr{F}$ as
\begin{equation}
\begin{array}{lll}
\big(r_i \, u \big)(x) = \widehat{\pi}_i(x_{{\ne}_{i}^+}) \, \sum_{x_i\in \mathcal{L}} \widehat{\pi}_i(x_{{\ne}_{i}^+}) \, u(x), \, \forall x \in \mathcal{L}^n, \,  \forall i \in \mathcal{V}
\end{array}
\nonumber
\end{equation}
Then, suppose $i \in \mathcal{V}$, and an arbitrary function $u \in \mathscr{F}$. One then has $\forall x \in \mathcal{L}^n$
\begin{equation}
\begin{array}{lll}
\Big(r_i \, \big( r_i \, u \big) \Big)(x) &=&   \widehat{\pi}_i(x_{{\ne}_{i}^+}) \, \sum_{x_i\in \mathcal{L}} \widehat{\pi}_i(x_{{\ne}_{i}^+}) \,
\Big( \widehat{\pi}_i(x_{{\ne}_{i}^+}) \, \sum_{x_i\in \mathcal{L}} \widehat{\pi}_i(x_{{\ne}_{i}^+}) \, u(x) \Big)\\
&=&   \widehat{\pi}_i(x_{{\ne}_{i}^+}) \, \big(\sum_{x_i\in \mathcal{L}} \widehat{\pi}_i^2(x_{{\ne}_{i}^+}) \big) \, \big( \sum_{x_i\in \mathcal{L}} \widehat{\pi}_i(x_{{\ne}_{i}^+}) \, u(x) \big)
\end{array}
\nonumber
\end{equation}
Since by construction, one has $\sum_{x_i\in \mathcal{L}} \widehat{\pi}_i^2(x_{{\ne}_{i}^+}) = 1$, $\forall x_{{\ne}_{i}} \in \mathcal{L}^{|{\ne}_{i}|}$, one
derives $\forall x \in \mathcal{L}^n$
\begin{equation}
\begin{array}{lll}
\Big(r_i \, \big( r_i \, u \big) \Big)(x) &=&  \widehat{\pi}_i(x_{{\ne}_{i}^+}) \,  \sum_{x_i\in \mathcal{L}} \widehat{\pi}_i(x_{{\ne}_{i}^+}) \, u(x) \\
&=& \big( r_i \, u \big)(x)
\end{array}
\nonumber
\end{equation}
Furthermore, since $u$ is assumed to be arbitrary, one derives
\begin{equation}
\begin{array}{lll}
r_i \, \big( r_i \, u \big) \equiv r_i \, u, \, \forall u \in \mathscr{F}
\end{array}
\nonumber
\end{equation}
which means that $r_i$ is a projection. Next, to show that $r_i$ is self-adjoint, suppose two arbitrary functions $u,v \in \mathscr{F}$. One has
\begin{equation}
\begin{array}{lll}
\langle v, r_i \, u \rangle  &=&  \sum_{x \in \mathcal{L}^n} v(x) \, \Big( \widehat{\pi}_i(x_{{\ne}_{i}^+}) \, \sum_{x_i\in \mathcal{L}} \widehat{\pi}_i(x_{{\ne}_{i}^+}) \, u(x) \Big) \\
&=&  \sum_{x_{\neq i} \in \mathcal{L}^{n-1}} \Big( \sum_{x_i\in \mathcal{L}} \widehat{\pi}_i(x_{{\ne}_{i}^+}) \, v(x) \Big) \, \Big( \sum_{x_i\in \mathcal{L}} \widehat{\pi}_i(x_{{\ne}_{i}^+}) \, u(x) \Big) \\
&=&  \sum_{x \in \mathcal{L}^{n}} \Big( \sum_{x_i\in \mathcal{L}} \widehat{\pi}_i(x_{{\ne}_{i}^+}) \, v(x) \Big) \, \Big( \widehat{\pi}_i(x_{{\ne}_{i}^+}) \, u(x) \Big) \\
&=&  \sum_{x \in \mathcal{L}^{n}}  u(x) \, \Big( \widehat{\pi}_i(x_{{\ne}_{i}^+}) \,  \Big( \sum_{x_i\in \mathcal{L}} \widehat{\pi}_i(x_{{\ne}_{i}^+}) \, v(x) \Big) \Big) \\
&=& \langle u, r_i \, v \rangle
\end{array}
\nonumber
\end{equation}
Since both $u$ and $v$ are assumed to be arbitrary, one derives
\begin{equation}
\begin{array}{lll}
\langle v, r_i \, u \rangle  = \langle u, r_i \, v \rangle , \, \forall u,v \in \mathscr{F}
\end{array}
\nonumber
\end{equation}
implying that $r_i$ is self-adjoint, hence. One concludes that $r_i$ is an orthogonal projection.
Finally, the proof that $k_i$ is an orthogonal projection follows immediately from the formula $k_i = \idoper - r_i$
which thus ends the proof of Proposition \ref{prop:ki_is_ortho_project}.
\end{proof}

\subsection{Proof of Proposition \ref{prop:K_k_i} }
\label{app:proof_prop_K_k_i}
\begin{proof}
The proof of Proposition \ref{prop:K_k_i} is straightforward.
\end{proof}


\subsection{Proof of Proposition \ref{prop:kilde_eq}}
\label{app:proof_prop_kilde_eq}
\begin{proof}
Suppose $i \in \mathcal{V}$, and arbitrary function $u \in \mathscr{F}$, and put $v \equiv \mathscr{O}_{s\cup {\ne}_{i}^+} \, u \equiv \frac{{u}_{s\cup {\ne}_{i}^+}}{L^{n-|s\cup {\ne}_{i}^+|}}$,
with ${u}_{s\cup {\ne}_{i}^+}$ denoting the margin of $u$ with respect to $s\cup {\ne}_{i}^+$. Then, one can  write
\begin{equation}
\begin{array}{lll}
\Big({k}_i \,  \big( \mathscr{O}_{s\cup {\ne}_{i}^+} \, u \big)\big)(x) &=& \big({k}_i \,  v\big)(x)  \\
                                                          &=&  v(x) - \widehat{\pi}_i(x_{{\ne}_{i}^+}) \, \sum_{x_i\in \mathcal{L}}  \widehat{\pi}_i(x_{{\ne}_{i}^+}) \, v(x) \\
                                                          &=& \Big(\frac{{u}_{s\cup {\ne}_{i}^+}(x_{s\cup {\ne}_{i}^+})}{L^{n-|s\cup {\ne}_{i}^+|}} \Big) -
\widehat{\pi}_i(x_{{\ne}_{i}^+}) \, \sum_{x_i\in \mathcal{L}}  \widehat{\pi}_i(x_{{\ne}_{i}^+}) \, \Big(\frac{{u}_{s\cup {\ne}_{i}^+}(x_{s\cup {\ne}_{i}^+})}{L^{n-|s\cup {\ne}_{i}^+|}} \Big) \\
                                                         &=& \big(\widehat{k}_i^{(s)} \, u \big)(x), \, \forall x \in \mathcal{L}^n
\end{array}
\nonumber
\end{equation}
Moreover, since $u$ is assumed to be arbitrary, one derives
\begin{equation}
\begin{array}{lll}
{k}_i \, \mathscr{O}_{s\cup {\ne}_{i}^+} = \widehat{k}_i^{(s)}
\end{array}
\label{fom_ki_eq_ki_osnip}
\end{equation}
Second, put $w \equiv {k}_i \, u$, and denote by ${w}_{s\cup {\ne}_{i}^+}$ the margin of $w$ with respect to $s\cup {\ne}_{i}^+$. One has
\begin{equation}
\begin{array}{lll}
w(x) \equiv u(x) - \widehat{\pi}_i(x_{{\ne}_{i}^+}) \, \sum_{x_i\in \mathcal{L}} \widehat{\pi}_i(x_{{\ne}_{i}^+}) \, u(x), \, \forall x \in \mathcal{L}^n
\end{array}
\nonumber
\end{equation}
Then, by marginalization of both hand sides of the latter formula with respect to $s\cup {\ne}_{i}^+$, one derives
\begin{equation}
\begin{array}{lll}
{w}_{s\cup {\ne}_{i}^+}(x_{s\cup {\ne}_{i}^+}) &=& u_{s\cup {\ne}_{i}^+}(x_{s\cup {\ne}_{i}^+}) -
\sum_{j \in \mathcal{V}/(s\cup{\ne}_{i}^+)} \sum_{x_j \in \mathcal{L}} \big( \widehat{\pi}_i(x_{{\ne}_{i}^+}) \, \sum_{x_i\in \mathcal{L}} \widehat{\pi}_i(x_{{\ne}_{i}^+}) \, u(x)\big) \\
&=& u_{s\cup {\ne}_{i}^+}(x_{s\cup {\ne}_{i}^+}) -
\widehat{\pi}_i(x_{{\ne}_{i}^+}) \, \sum_{x_i\in \mathcal{L}} \widehat{\pi}_i(x_{{\ne}_{i}^+}) \, \big( \sum_{j \in \mathcal{V}/(s\cup{\ne}_{i}^+)} \sum_{x_j \in \mathcal{L}}u(x)\big) \\
&=& u_{s\cup {\ne}_{i}^+}(x_{s\cup {\ne}_{i}^+}) -
\widehat{\pi}_i(x_{{\ne}_{i}^+}) \, \sum_{x_i\in \mathcal{L}} \widehat{\pi}_i(x_{{\ne}_{i}^+}) \, u_{s\cup {\ne}_{i}^+}(x_{s\cup {\ne}_{i}^+}), \, \forall x_{s\cup {\ne}_{i}^+} \in \mathcal{L}^{|s\cup {\ne}_{i}^+|}
\end{array}
\nonumber
\end{equation}
hence,
\begin{equation}
\begin{array}{lll}
\mathscr{O}_{s\cup {\ne}_{i}^+} \, \big(k_i \, u \big)(x_{s\cup {\ne}_{i}^+})  &=&  \frac{{w}_{s\cup {\ne}_{i}^+}(x_{s\cup {\ne}_{i}^+})}{L^{n-|s\cup {\ne}_{i}^+|}} \\
&=& \Big( \, \Big(\frac{{u}_{s\cup {\ne}_{i}^+}(x_{s\cup {\ne}_{i}^+})}{L^{n-|s\cup {\ne}_{i}^+|}} \Big) -
\widehat{\pi}_i(x_{{\ne}_{i}^+}) \, \sum_{x_i\in \mathcal{L}} \widehat{\pi}_i(x_{{\ne}_{i}^+}) \, \Big( \frac{{u}_{s\cup {\ne}_{i}^+}(x_{s\cup {\ne}_{i}^+})}{L^{n-|s\cup {\ne}_{i}^+|}}\Big) \,\Big) \\
&=&  \big(\widehat{k}_i^{(s)} \, u \big)(x_{s\cup {\ne}_{i}^+}), \, \forall x_{s\cup {\ne}_{i}^+} \in \mathcal{L}^{|s\cup {\ne}_{i}^+|}
\end{array}
\nonumber
\end{equation}
Moreover, since $u$ is assumed to be arbitrary, one derives
\begin{equation}
\begin{array}{lll}
\mathscr{O}_{s\cup {\ne}_{i}^+} \, {k}_i  = \widehat{k}_i^{(s)}
\end{array}
\label{fom_ki_eq_ki_osnip_b}
\end{equation}
Combining both formulas (\ref{fom_ki_eq_ki_osnip}) and (\ref{fom_ki_eq_ki_osnip_b}), one finally establishes
\begin{equation}
\begin{array}{lll}
 \widehat{k}_i^{(s)} =  {k}_i \, \mathscr{O}_{s\cup {\ne}_{i}^+} = \mathscr{O}_{s\cup {\ne}_{i}^+} \, {k}_i  =  \mathscr{O}_{s\cup {\ne}_{i}^+} \, {k}_i  \, \mathscr{O}_{s\cup {\ne}_{i}^+}
\end{array}
\nonumber
\end{equation}
thus proving Proposition \ref{prop:kilde_eq}.
\end{proof}


\subsection{Proof of Proposition \ref{prop:kilde_i_is_ortho_project}}
\label{app:proof_prop_kilde_i_is_ortho_project}
\begin{proof}
Suppose $i \in \mathcal{V}$. By Proposition \ref{prop:kilde_eq}, one can write
\begin{equation}
\begin{array}{lll}
\widehat{k}_i^{(s)} \, \widehat{k}_i^{(s)} &=& \big( \mathscr{O}_{s\cup {\ne}_{i}^+} \, {k}_i \big) \, \big( {k}_i \, \mathscr{O}_{s\cup {\ne}_{i}^+})   \\
                                   &=&  \mathscr{O}_{s\cup {\ne}_{i}^+} \, \big({k}_i \,  {k}_i\big) \, \mathscr{O}_{s\cup {\ne}_{i}^+}
\end{array}
\nonumber
\end{equation}
Moreover, by Proposition \ref{prop:ki_is_ortho_project}, one has ${k}_i \,  {k}_i = k_i$, hence,
\begin{equation}
\begin{array}{lll}
\widehat{k}_i^{(s)} \, \widehat{k}_i^{(s)} =  \mathscr{O}_{s\cup {\ne}_{i}^+} \,  {k}_i \, \mathscr{O}_{s\cup {\ne}_{i}^+}
\end{array}
\nonumber
\end{equation}
Then, since by Proposition \ref{prop:kilde_eq}, one has
$\widehat{k}_i^{(s)} =  \mathscr{O}_{s\cup {\ne}_{i}^+} \,  {k}_i \, \mathscr{O}_{s\cup {\ne}_{i}^+}$, one derives
\begin{equation}
\begin{array}{lll}
\widehat{k}_i^{(s)} \, \widehat{k}_i^{(s)} = \widehat{k}_i^{(s)}
\end{array}
\nonumber
\end{equation}
which means that $\widehat{k}_i^{(s)}$ is a projection. Furthermore, by Proposition \ref{prop:kilde_eq}, one has
\begin{equation}
\begin{array}{lll}
\widehat{k}_i^{(s)} = \mathscr{O}_{s\cup {\ne}_{i}^+} \,  {k}_i \, \mathscr{O}_{s\cup {\ne}_{i}^+}
\end{array}
\nonumber
\end{equation}
in such a way that, by equating the transposes of both hand sides of the latter formula, and since both
$k_i$, and $\mathscr{O}_{s\cup {\ne}_{i}^+}$ are orthogonal projections, hence, self-adjoint linear
operators, one deduces that $\widehat{k}_i^{(s)}$ coincides with its own transpose which means that
$\widehat{k}_i^{(s)}$ is a self-adjoint linear operator. In summary, we have shown that $\widehat{k}_i^{(s)}$
is a self-adjoint projection, meaning that $\widehat{k}_i^{(s)}$ is an orthogonal projection, thus proving
Proposition \ref{prop:kilde_i_is_ortho_project}.
\end{proof}


\subsection{Proof of Proposition \ref{prop:K_k_tilde_i}}
\label{app:proof_prop_K_k_tilde_i}
\begin{proof}
The proof of Proposition \ref{prop:K_k_tilde_i} is straightforward.
\end{proof}


\subsection{Proof of Proposition \ref{prop:markov_lin_oper_s}}
\label{app:proof_prop_markov_lin_oper_s}
\begin{proof}
By Proposition \ref{prop:o_s_ortho}, $\mathscr{O}_{\emptyset}$ is an orthogonal projection, hence,
$L^{n}\, \mathscr{O}_{\emptyset}$ is a self-adjoint positive semi-definite linear operator, moreover, by
Proposition \ref{prop:ki_is_ortho_project}, $\forall i \in \mathcal{V}$, $ \widehat{k}_i^{(s)}$ is an orthogonal
projection, hence, a self-adjoint positive semi-definite linear operator. Then, by Proposition \ref{prop:K_k_tilde_i},
$\hat{\mathscr{K}}^{(s)}$ writes as a sum of self-adjoint positive semi-definite linear operators,  hence,
$\hat{\mathscr{K}}^{(s)}$ is, in turn, a self-adjoint positive semi-definite linear operator,
thus proving Proposition \ref{prop:markov_lin_oper_s}.
\end{proof}



\subsection{Proof of Proposition \ref{prop:ks_eq_kos}}
\label{app:proof_prop_ks_eq_kos}

\begin{proof}
Suppose an arbitrary function $u \in \mathscr{F}$, and denote $v \equiv {\mathscr{O}}_s \,u \equiv \frac{u_s}{L^{n-|s|}}$,
with $u_s$ standing for the margin of $u$ with respect to $s$. Furthermore, $\forall i \in \mathcal{V}$, denote $v_{s\cup {\ne}_{i}^+}$
the margin of $v$ with respect to $s\cup {\ne}_{i}^+$. One then easily checks that
\begin{equation}
\begin{array}{lll}
v_{s\cup {\ne}_{i}^+} \equiv {L^{n-|s\cup {\ne}_{i}^+|}} \, \frac{u_s}{L^{n-|s|}}, \, \forall i \in \mathcal{V}
\nonumber
\end{array}
\end{equation}
in such a way that, one writes $\forall x \in \mathcal{L}^n$
\begin{equation}
\begin{array}{lll}
\Big(\hat{\mathscr{K}}^{(s)} \,  \big({\mathscr{O}}_s \,u \big) \Big)(x)  &=& \big(\hat{\mathscr{K}} \, v \big) \big)(x) \\
&=& \sum_{x \in \mathcal{L}^n} v(x) + \sum_{i \in \mathcal{V}} \frac{ v_{s\cup {\ne}_{i}^+}(x_{s\cup {\ne}_{i}^+}) -
\widehat{\pi}_i(x_{{\ne}_{i}^+}) \, \sum_{x_i\in \mathcal{L}} \widehat{\pi}_i(x_{{\ne}_{i}^+}) \, v_{s\cup {\ne}_{i}^+}(x_{s\cup {\ne}_{i}^+})}{L^{n-|s\cup {\ne}_{i}^+|}} \\
&=& \sum_{x \in \mathcal{L}^n} \big(\frac{u_s(x_s)}{L^{n-|s|}}\big) + \sum_{i \in \mathcal{V}} \frac{ {L^{n-|s\cup {\ne}_{i}^+|}} \, \frac{u_s(x_s)}{L^{n-|s|}} -
\widehat{\pi}_i(x_{{\ne}_{i}^+}) \, \sum_{x_i\in \mathcal{L}} \widehat{\pi}_i(x_{{\ne}_{i}^+}) \, \big({L^{n-|s\cup {\ne}_{i}^+|}} \, \frac{u_s(x_s)}{L^{n-|s|}} \big)}{L^{n-|s\cup {\ne}_{i}^+|}} \\
&=& \sum_{x \in \mathcal{L}^n} \big(\frac{u_s(x_s)}{L^{n-|s|}}\big) + \sum_{i} \big( u_{s}(x_{s}) - \widehat{\pi}_i(x_{{\ne}_{i}^+}) \, \sum_{x_i\in \mathcal{L}} \widehat{\pi}_i(x_{{\ne}_{i}^+}) \, u_{s}(x_{s}) \big) \\
&=& \Big({\mathscr{K}}\, \big({\mathscr{O}}_s \,u \big) \Big)(x)
\nonumber
\end{array}
\end{equation}
Since function $u$ is assumed to be arbitrary, one derives
\begin{equation}
\begin{array}{lll}
 \hat{\mathscr{K}}\,  \big({\mathscr{O}}_s \,u \big) \equiv \mathscr{K} \,  \big({\mathscr{O}}_s \,u \big), \, \forall u \in \mathscr{F}
\nonumber
\end{array}
\end{equation}
hence, $ \hat{\mathscr{K}} \, {\mathscr{O}}_s = \mathscr{K} \,  {\mathscr{O}}_s$, thus proving Proposition \ref{prop:ks_eq_kos}.
\end{proof}


\subsection{Proof of Proposition \ref{prop:kts_p_sar_eq_one}}
\label{app:proof_kts_p_sar_eq_one}

\begin{proof}
First, suppose the following system of linear identities where the unknown is $p \in \mathscr{F}$:
\begin{equation}
\begin{array}{l}
p(x) -  \widehat{\pi}_i(x_{{\ne}_{i}^+}) \, \sum_{x_i\in \mathcal{L}} \widehat{\pi}_i(x_{{\ne}_{i}^+}) \, p(x) = 0, \, \forall x \in \mathcal{L}^n,\, \forall i \in \mathcal{V}
\end{array}
\label{lin_sys_unicity_mrf_tilde}
\end{equation}
Second, suppose $p$ is a solution of system of linear identities (\ref{lin_sys_unicity_mrf_tilde}), and denote $\forall i \in \mathcal{V}$
by $p_{s\cup {\ne}_{i}^+}$ the margin of $p$ with respect to hypervertex $s\cup {\ne}_{i}^+$, that is to say, one has
\begin{equation}
\begin{array}{l}
p_{s\cup {\ne}_{i}^+}\big(x_{s\cup {\ne}_{i}^+}\big) = \sum_{j \in V/s\cup {\ne}_{i}^+} \sum_{x_j \in \mathcal{L}} p(x), \, \forall x_{s\cup {\ne}_{i}^+} \in \mathcal{L}^{|{s\cup {\ne}_{i}^+}|}
\end{array}
\nonumber
\end{equation}
Then, assume an arbitrary node $i \in \mathcal{V}$. Thus, $p$ satisfies the following linear identity:
\begin{equation}
\begin{array}{l}
p(x) -  \widehat{\pi}_i(x_{{\ne}_{i}^+}) \, \sum_{x_i\in \mathcal{L}} \widehat{\pi}_i(x_{{\ne}_{i}^+}) \, p(x) = 0, \, \forall x \in \mathcal{L}^n
\end{array}
\label{lin_sys_unicity_mrf_tilde_b}
\end{equation}
By marginalization of both hands side of the latter identity with respect to all the nodes in $V/s\cup {\ne}_{i}^+$, one finds
\begin{equation}
\begin{array}{l}
\sum_{j \in V/s\cup {\ne}_{i}^+} \sum_{x_j \in \mathcal{L}} p(x) -
\sum_{j \in V/s\cup {\ne}_{i}^+}  \sum_{x_j \in \mathcal{L}} \big(\widehat{\pi}_i(x_{{\ne}_{i}^+}) \, \sum_{x_i\in \mathcal{L}} \widehat{\pi}_i(x_{{\ne}_{i}^+}) \, p(x) \big) = 0, \, \forall x_{s\cup {\ne}_{i}^+}
\in \mathcal{L}^{|{s\cup {\ne}_{i}^+}|}
\end{array}
\nonumber %
\end{equation}
Then, observe that, $\forall x_{s\cup {\ne}_{i}^+} \in \mathcal{L}^{|{s\cup {\ne}_{i}^+}|}$, one has
\begin{equation}
\begin{array}{lll}
\sum_{j \in V/s\cup {\ne}_{i}^+}  \sum_{x_j \in \mathcal{L}} \big(\widehat{\pi}_i(x_{{\ne}_{i}^+}) \, \sum_{x_i\in \mathcal{L}} \widehat{\pi}_i(x_{{\ne}_{i}^+}) \, p(x) \big) &=&
\widehat{\pi}_i(x_{{\ne}_{i}^+}) \, \sum_{x_i\in \mathcal{L}} \widehat{\pi}_i(x_{{\ne}_{i}^+}) \, \big( \sum_{j \in V/s\cup {\ne}_{i}^+}  \sum_{x_j \in \mathcal{L}}  p(x) \big) \\
&=& \widehat{\pi}_i(x_{{\ne}_{i}^+}) \, \sum_{x_i\in \mathcal{L}} \widehat{\pi}_i(x_{{\ne}_{i}^+}) \, p_{s\cup {\ne}_{i}^+}\big(x_{s\cup {\ne}_{i}^+}\big)
\end{array}
\nonumber
\end{equation}
in such a way that, linear identity (\ref{lin_sys_unicity_mrf_tilde_b}) can rewrite as
\begin{equation}
\begin{array}{l}
p_{s\cup {\ne}_{i}^+}\big(x_{s\cup {\ne}_{i}^+}\big) -
\widehat{\pi}_i(x_{{\ne}_{i}^+}) \, \sum_{x_i\in \mathcal{L}} \widehat{\pi}_i(x_{{\ne}_{i}^+}) \, p_{s\cup {\ne}_{i}^+}\big(x_{s\cup {\ne}_{i}^+}\big) = 0 , \, \forall x_{s\cup {\ne}_{i}^+} \in \mathcal{L}^{|{s\cup {\ne}_{i}^+}|}
\end{array}
\nonumber 
\end{equation}
Since both $p$ and $i$ are assumed to be arbitrary, one derives that any solution $p$ of system of linear identities (\ref{lin_sys_unicity_mrf_tilde}) is
also a solution of the following system of linear identities:
\begin{equation}
\begin{array}{l}
p_{s\cup {\ne}_{i}^+}\big(x_{s\cup {\ne}_{i}^+}\big) -
\widehat{\pi}_i(x_{{\ne}_{i}^+}) \, \sum_{x_i\in \mathcal{L}} \widehat{\pi}_i(x_{{\ne}_{i}^+}) \, p_{s\cup {\ne}_{i}^+}\big(x_{s\cup {\ne}_{i}^+}\big) = 0 ,
\, \forall x_{s\cup {\ne}_{i}^+} \in \mathcal{L}^{|{s\cup {\ne}_{i}^+}|},\, \forall i \in \mathcal{V}
\end{array}
\nonumber 
\end{equation}
hence, satisfies
\begin{equation}
\begin{array}{lll}
\big(\hat{\mathscr{K}}^{(s)}  \, p\big)(x) &=& \sum_{x\in \mathcal{L}^n} p(x)+ \sum_{i \in \mathcal{V}} \Big( p_{s\cup {\ne}_{i}^+}\big(x_{s\cup {\ne}_{i}^+}\big) -
\widehat{\pi}_i(x_{{\ne}_{i}^+}) \, \sum_{x_i\in \mathcal{L}} \widehat{\pi}_i(x_{{\ne}_{i}^+}) \, p_{s\cup {\ne}_{i}^+}\big(x_{s\cup {\ne}_{i}^+}\big) \Big) \\
&=& \sum_{x\in \mathcal{L}^n} p(x)+ \sum_{i \in \mathcal{V}} 0 \\
&=& \sum_{x\in \mathcal{L}^n} p(x), \, \forall x \in \mathcal{L}^{n}
\end{array}
\nonumber
\end{equation}
In particular, by Lemma \ref{lemma:unicity_mrf}, MRF joint distribution $\mathbb{P}\big(\mathcal{X} = \cdot \big)$ satisfies system
of linear identities (\ref{lin_sys_unicity_mrf_tilde}), hence, satisfies
\begin{equation}
\begin{array}{lll}
\hat{\mathscr{K}}^{(s)}  \, \mathbb{P}\big(\mathcal{X} = \cdot \big) &\equiv& \big( \sum_{x \in \mathcal{L}^n} \mathbb{P}\big(\mathcal{X} = x\big) \big) \, \oneindic \\
 &\equiv& \oneindic
\end{array}
\nonumber
\end{equation}
thus proving Proposition \ref{prop:kts_p_sar_eq_one}.
\end{proof}


\subsection{Proof of Proposition \ref{prop:kbar_i_is_ortho_project}}
\label{app:proof_prop_kbar_i_is_ortho_project}

\begin{proof}
Suppose $i \in \mathcal{V}$. One then easily checks that one can write $\bar{k}_i^{(s)} = {k}_i - \widehat{k}_i^{(s)}$, in such a way that, one can write
\begin{equation}
\begin{array}{lll}
\bar{k}_i^{(s)} \, \bar{k}_i^{(s)}  = {k}_i \, {k}_i - {k}_i \,  \widehat{k}_i^{(s)} -  \widehat{k}_i^{(s)} \, {k}_i +  \widehat{k}_i^{(s)}\,  \widehat{k}_i^{(s)}
\end{array}  \nonumber
\end{equation}
Since by Propositions \ref{prop:ki_is_ortho_project} and \ref{prop:kilde_i_is_ortho_project}, respectively, both $k_i$ and $\widehat{k}_i^{(s)}$ are
orthogonal projections, one has ${k}_i \, {k}_i = {k}_i$, and $\widehat{k}_i^{(s)}\,  \widehat{k}_i^{(s)} = \widehat{k}_i^{(s)}$, hence,
\begin{equation}
\begin{array}{lll}
\bar{k}_i^{(s)} \, \bar{k}_i^{(s)}  = {k}_i  - {k}_i \,  \widehat{k}_i^{(s)} -  \widehat{k}_i^{(s)} \, {k}_i +  \widehat{k}_i^{(s)}
\end{array}  \label{bar_kis_kis}
\end{equation}
Next, by using the following identity (from Proposition \ref{prop:kilde_eq}):
\begin{equation}
\begin{array}{lll}
\widehat{k}_i^{(s)} =   {k}_i \, \mathscr{O}_{s\cup {\ne}_{i}^+} =  \mathscr{O}_{s\cup {\ne}_{i}^+} \, {k}_i
\end{array}  \nonumber
\end{equation}
one derives
\begin{equation}
\begin{array}{lll}
{k}_i \,  \widehat{k}_i^{(s)} &=& {k}_i \,  \big( {k}_i \,  \mathscr{O}_{s\cup {\ne}_{i}^+} \big) \\
                          &=& \big( {k}_i \,{k}_i \big) \,  \mathscr{O}_{s\cup {\ne}_{i}^+} \\
                          &=& {k}_i \,  \mathscr{O}_{s\cup {\ne}_{i}^+} \\
                          &=&  \widehat{k}_i^{(s)}
\end{array}  \nonumber
\end{equation}
and
\begin{equation}
\begin{array}{lll}
\widehat{k}_i^{(s)} \, {k}_i &=& \big( \mathscr{O}_{s\cup {\ne}_{i}^+} \, {k}_i  \big) \, {k}_i \,  \\
                          &=&  \mathscr{O}_{s\cup {\ne}_{i}^+}  \, \big( {k}_i \,{k}_i \big) \\
                          &=& \mathscr{O}_{s\cup {\ne}_{i}^+} \, {k}_i \\
                          &=&  \widehat{k}_i^{(s)}
\end{array}  \nonumber
\end{equation}
in such a way that, by substituting the expressions both of ${k}_i \,  \widehat{k}_i^{(s)}$ and $\widehat{k}_i^{(s)} \, {k}_i $ with $\widehat{k}_i^{(s)}$
in formula (\ref{bar_kis_kis}), one derives
\begin{equation}
\begin{array}{lll}
\bar{k}_i^{(s)} \, \bar{k}_i^{(s)}  &=& {k}_i  -  \widehat{k}_i^{(s)} -  \widehat{k}_i^{(s)}  + \widehat{k}_i^{(s)} \\
                                     &=& {k}_i  -  \widehat{k}_i^{(s)} \\
                                     &=& \bar{k}_i^{(s)}
\end{array}
\end{equation}
which means that $\bar{k}_i^{(s)}$ is a projection. Furthermore, since $\bar{k}_i^{(s)}$ writes as the difference
of ${k}_i$ and $\widehat{k}_i^{(s)}$ which both are self-adjoint linear operators, one finds that $\bar{k}_i^{(s)}$ is,
in turn, self-adjoint. One concludes that $\bar{k}_i^{(s)}$ is a self-adjoint projection, hence, an orthogonal
projection, thus proving Proposition \ref{prop:kbar_i_is_ortho_project}.
\end{proof}


\subsection{Proof of Proposition \ref{prop:kbar_i_is_sum_ortho_project}}
\label{app:proof_prop_kbar_i_is_sum_ortho_project}
\begin{proof}
The proof of Proposition \ref{prop:kbar_i_is_sum_ortho_project} is straightforward.
\end{proof}


\subsection{Proof of Proposition \ref{prop:markov_lin_oper_s_r}}
\label{app:proof_prop_markov_lin_oper_s_r}

\begin{proof}
The proof of Proposition \ref{prop:markov_lin_oper_s_r} follows immediately from the result of Proposition \ref{prop:kbar_i_is_sum_ortho_project}
which states that $\widebar{\mathscr{K}}^{(s)}$ writes as a sum of orthogonal projections, hence, as a sum of self-adjoint positive semi-definite linear
operators.
\end{proof}


\subsection{Proof of Proposition \ref{prop:rs_eq_rosb}}
\label{app:proof_prop_rs_eq_rosb}

\begin{proof}
Assume an arbitrary function $u \in \mathscr{F}$, and put $v \equiv {\mathscr{O}}_{s} \, u$. Next, assume an arbitrary $i \in \mathcal{V}$.
One then easily checks that the margin of $v$ with respect to $s\cup {\ne}_{i}^+$ is given by $v_{s\cup {\ne}_{i}^+} \equiv L^{n-|s\cup {\ne}_{i}^+|} \, v$,
in such a way that, $v - \frac{v_{s\cup {\ne}_{i}^+}}{L^{n-|s\cup {\ne}_{i}^+|}} \equiv \zeroindic$, and hence,
\begin{equation}
\begin{array}{lll}
\bar{k}_i^{(s)} \, \big( {\mathscr{O}}_{s} \, u \big)  \equiv \bar{k}_i^{(s)} \, v \equiv  \zeroindic
\end{array}  \nonumber
\end{equation}
One derives
\begin{equation}
\begin{array}{lll}
\bar{k}_i^{(s)} \, u  &\equiv& \bar{k}_i^{(s)} \, \big( {\mathscr{O}}_{s} \, u + \widebar{\mathscr{O}}_{s} \, u \big)  \\
                      &\equiv& \bar{k}_i^{(s)} \, \big( {\mathscr{O}}_{s} \, u \big) + \bar{k}_i^{(s)} \, \big(\widebar{\mathscr{O}}_{s} \, u \big)  \\
                      &\equiv& \zeroindic + \bar{k}_i^{(s)} \, \big(\widebar{\mathscr{O}}_{s} \, u \big) \\
                      &\equiv& \bar{k}_i^{(s)} \, \big(\widebar{\mathscr{O}}_{s} \, u \big)
\end{array}  \nonumber
\end{equation}
hence,
\begin{equation}
\begin{array}{lll}
\widebar{\mathscr{K}}^{(s)} \, u  &\equiv&  \sum_{i \in \mathcal{V}} \bar{k}_i^{(s)} \, u \\
                               &\equiv&  \sum_{i \in \mathcal{V}} \bar{k}_i^{(s)} \, \big(\widebar{\mathscr{O}}_{s} \, u \big) \\
                               &\equiv&  \big(\sum_{i \in \mathcal{V}} \bar{k}_i^{(s)} \big) \big(\widebar{\mathscr{O}}_{s} \, u \big) \\
                               &\equiv&   \widebar{\mathscr{K}}^{(s)} \, \big(\widebar{\mathscr{O}}_{s} \, u \big)
\end{array}  \nonumber
\end{equation}
Since $u$ is assumed to be arbitrary, one deduces that
\begin{equation}
\begin{array}{lll}
\widebar{\mathscr{K}}^{(s)} = \widebar{\mathscr{K}}^{(s)} \, \widebar{\mathscr{O}}_{s}
\end{array}   \label{bar_ks_os_b_1}
\end{equation}
Moreover, since $\widebar{\mathscr{K}}^{(s)}$, $\widebar{\mathscr{K}}^{(s)}$, and  $\widebar{\mathscr{O}}_{s}$ are self-adjoint
linear operators, by equating the transposes of both hand sides of identity (\ref{bar_ks_os_b_1}), one also finds
\begin{equation}
\begin{array}{lll}
\widebar{\mathscr{K}}^{(s)} = \widebar{\mathscr{O}}_{s} \, \widebar{\mathscr{K}}^{(s)}
\end{array}  \label{bar_ks_os_b_2}
\end{equation}
Finally, by combining both formulas (\ref{bar_ks_os_b_1}) and (\ref{bar_ks_os_b_2}), one establishes
\begin{equation}
\begin{array}{lll}
\widebar{\mathscr{K}}^{(s)} =  \widebar{\mathscr{K}}^{(s)} \, \widebar{\mathscr{O}}_{s} = \widebar{\mathscr{O}}_{s} \, \widebar{\mathscr{K}}^{(s)} =  \widebar{\mathscr{O}}_{s} \, \widebar{\mathscr{K}}^{(s)} \, \widebar{\mathscr{O}}_{s}
\end{array}  \nonumber
\end{equation}
hence the proof of Proposition \ref{prop:rs_eq_rosb}.
\end{proof}


\subsection{Proof of Proposition \ref{prop:rts_p_sar_eq_zero}}
\label{app:proof_prop_rts_p_sar_eq_zero}.

\begin{proof}
By Theorem \ref{theo:markov_distrib_unique_pos_definite} and Proposition \ref{prop:kts_p_sar_eq_one}, respectively, one has
\begin{equation}
\begin{array}{lll}
\widebar{\mathscr{K}}^{(s)} \, \mathbb{P}\big(\mathcal{X} = \cdot \big) &\equiv&  {\mathscr{K}} \, \mathbb{P}\big(\mathcal{X} = \cdot \big) -  \hat{\mathscr{K}}^{(s)} \, \mathbb{P}\big(\mathcal{X} = \cdot \big) \\
 &\equiv&  \oneindic - \oneindic \\
 &\equiv& \zeroindic
\end{array}  \nonumber
\end{equation}
hence, the proof of Proposition \ref{prop:rts_p_sar_eq_zero}.
\end{proof}


\subsection{Proof of Proposition \ref{prop:kkos_eq_kk}}
\label{app:proof_prop_kkos_eq_kk}

\begin{proof}
Denote $p^* = \mathbb{P}\big(\mathcal{X} = \cdot \big)$. The general solution of identity (\ref{mrf_dist_solu_lin_sys_polyn}) is then given by
\begin{equation}
\begin{array}{lll}
p \equiv p^* + \big(\idoper - \hat{\mathscr{K}}^{\dag} \, \hat{\mathscr{K}}\big)\, z, \, z \in \mathscr{F}
\end{array}
\nonumber
\end{equation}
implying that
\begin{equation}
\begin{array}{lll}
{\mathscr{O}}_s \, p \equiv {\mathscr{O}}_s \, p^* + {\mathscr{O}}_s \, \big(\idoper - \hat{\mathscr{K}}^{\dag} \, \hat{\mathscr{K}}\big)\, z, \, z \in \mathscr{F}
\end{array}
\nonumber
\end{equation}
Since, by Theorem \ref{theorem:founding_kp_eq_1}, any solution $p$ of identity (\ref{mrf_dist_solu_lin_sys_polyn}) must satisfy ${\mathscr{O}}_s \, p \equiv {\mathscr{O}}_s \, p^*$,
one concludes that $\forall  z \in \mathscr{F}$, one must have
\begin{equation}
\begin{array}{lll}
{\mathscr{O}}_s \, \big(\idoper - \hat{\mathscr{K}}^{\dag} \, \hat{\mathscr{K}}\big) \, z \equiv \zeroindic
\end{array}
\nonumber
\end{equation}
implying that
\begin{equation}
\begin{array}{lll}
{\mathscr{O}}_s  = {\mathscr{O}}_s \, \big( \hat{\mathscr{K}}^{\dag} \, \hat{\mathscr{K}}\big)
\end{array}
\nonumber
\end{equation}
and by taking the transpose of both hand sides of the latter operator equality, one also finds
\begin{equation}
\begin{array}{lll}
{\mathscr{O}}_s  =  \big( \hat{\mathscr{K}}^{\dag} \, \hat{\mathscr{K}}\big) \, {\mathscr{O}}_s
\end{array}
\nonumber
\end{equation}
hence the proof of Proposition \ref{prop:kkos_eq_kk}.
\end{proof}


\subsection{Proof of Proposition \ref{prop:marginal_qp}}
\label{app:proof_prop_marginal_qp}

\begin{proof}
First of all, since $p^*$ is a solution of identity (\ref{mrf_dist_solu_lin_sys_polyn}), one finds
\begin{equation}
\hat{\mathscr{K}}^{\dag} \,  \hat{\mathscr{K}} \, p^* \equiv \hat{\mathscr{K}}^{\dag} \, \oneindic
\nonumber
\end{equation}
implying thfat
\begin{equation}
{\mathscr{O}}_s \, \hat{\mathscr{K}}^{\dag} \,  \hat{\mathscr{K}} \, p^* \equiv {\mathscr{O}}_s \, \big( \hat{\mathscr{K}}^{\dag} \, \oneindic\big)
\nonumber
\end{equation}
Furthermore, by using Proposition \ref{prop:kkos_eq_kk}, one finds
\begin{equation}
{\mathscr{O}}_s \, p^* \equiv {\mathscr{O}}_s \, \big( \hat{\mathscr{K}}^{\dag} \, \oneindic\big)
\nonumber
\end{equation}
hence the proof of Proposition \ref{prop:marginal_qp}.
\end{proof}


\section{Proofs of Lemmas}

\subsection{Proof of Lemma \ref{lemma:besag}}
\label{app:proof_lemma_besag}
\begin{proof}
Suppose a second MRF with regard to graph $\mathcal{G}$ denoted by $\mathcal{Y}$, and satisfying
\begin{equation}
\begin{array}{lll}
 \mathbb{P} \big(\mathcal{Y}_i = x_i/\mathcal{Y}_{{\ne}_{i}} = x_{{\ne}_{i}}\big) = \mathbb{P} \big(\mathcal{X}_i = x_i/\mathcal{X}_{{\ne}_{i}} = x_{{\ne}_{i}}\big), \, \forall x_{{\ne}_{i}}^+ \in \mathcal{L}^{|{{\ne}_{i}}^+|}, \, \forall i \in \mathcal{V}
\end{array}
\label{mrf_unik_a}
\end{equation}
Since by assumption, both $\mathcal{X}$ and $\mathcal{Y}$ are assumed to be MRFs with regard to $\mathcal{G}$, they satisfy the Markovian property, that is
\begin{equation}
\forall i \in \mathcal{V}, \, \forall x \in \mathcal{L}^{n}, \,
\begin{cases}
\mathbb{P}\big(\mathcal{X}_i = x_i / \mathcal{X}_{\neq i} = x_{\neq i} \big) = \mathbb{P} \big(\mathcal{X}_i = x_i/\mathcal{X}_{{\ne}_{i}} = x_{{\ne}_{i}}\big)  \\
\mathbb{P}\big(\mathcal{Y}_i = x_i / \mathcal{Y}_{\neq i} = x_{\neq i} \big) = \mathbb{P} \big(\mathcal{Y}_i = x_i/\mathcal{Y}_{{\ne}_{i}} = x_{{\ne}_{i}}\big)
\end{cases}
\label{mrf_unik_b}
\end{equation}
Thus, by combining both formulas (\ref{mrf_unik_a}) and (\ref{mrf_unik_b}), one derives
\begin{equation}
\begin{array}{lll}
\mathbb{P}\big(\mathcal{X}_i = x_i / \mathcal{X}_{\neq i} = x_{\neq i} \big) = \mathbb{P}\big(\mathcal{Y}_i = x_i / \mathcal{Y}_{\neq i} = x_{\neq i} \big), \,
\forall x \in \mathcal{L}^{n}, \, \forall i \in \mathcal{V}
\end{array}
\label{mrf_unik_c}
\end{equation}
Since by assumption, one has $\mathbb{P}\big(\mathcal{X} = x \big) >0, \forall x \in \mathcal{L}^{n}$,  one can reexpress formula (\ref{mrf_unik_c}) as
\begin{equation}
\begin{array}{lll}
\frac{\mathbb{P}\big(\mathcal{Y} = x \big)}{\mathbb{P}\big(\mathcal{X} = x \big)} = \frac{\mathbb{P}\big(\mathcal{Y}_{\neq i} = x_{\neq i} \big)}{\mathbb{P}\big(\mathcal{X}_{\neq i} = x_{\neq i}\big)},
\, \forall x \in \mathcal{L}^{n}, \, \forall i \in \mathcal{V}
\end{array}
\label{mrf_unik_d}
\end{equation}
Denote
\begin{equation}
\begin{cases}
R(x) = \frac{\mathbb{P}\big(\mathcal{Y} = x \big)}{\mathbb{P}\big(\mathcal{X} = x \big)}, \forall x \in \mathcal{L}^{n} \\
R_i(x_{\neq i}) = \frac{\mathbb{P}\big(\mathcal{Y}_{\neq i} = x_{\neq i} \big)}{\mathbb{P}\big(\mathcal{X}_{\neq i} = x_{\neq i}\big)}, \, \forall x_{\neq i} \in \mathcal{L}^{n-1}, \, \forall i \in \mathcal{V}
\end{cases}
\nonumber
\end{equation}
Then, one can rewrite formula (\ref{mrf_unik_d}) as
\begin{equation}
\begin{array}{lll}
R(x) = R_i(x_{\neq i}), \, \, \forall x \in \mathcal{L}^{n}, \, \forall i \in \mathcal{V}
\end{array}
\nonumber 
\end{equation}
meaning that, $\forall i \in \mathcal{V}$, function $R(x)$ does not depend on variable $x_i$, implying in turn, that
function $R(x)$ does not depend on any of the variables $x_i$,  $\forall i \in \mathcal{V}$ which is only possible if
there exists a real constant $c$, such that, $R(x) = c, \forall x \in \mathcal{L}^{n}$, that is
\begin{equation}
\begin{array}{lll}
\frac{\mathbb{P}\big(\mathcal{Y} = x \big)}{\mathbb{P}\big(\mathcal{X} = x \big)} = c, \, \forall x \in \mathcal{L}^{n}
\end{array}
\nonumber
\end{equation}
hence,
\begin{equation}
\begin{array}{lll}
\mathbb{P}\big(\mathcal{Y} = x \big) = c \, \mathbb{P}\big(\mathcal{X} = x \big), \, \forall x \in \mathcal{L}^{n}
\end{array}
\label{mrf_unik_f}
\end{equation}
implying that
\begin{equation}
\begin{array}{lll}
\sum_{ x \in \mathcal{L}^{n}} \mathbb{P}\big(\mathcal{Y} = x \big) = c \, \sum_{ x \in \mathcal{L}^{n}} \mathbb{P}\big(\mathcal{X} = x \big), \, \forall
\end{array}
\nonumber
\end{equation}
and one finds $c = 1$. By substituting $c$ with its value $1$ in formula (\ref{mrf_unik_f}), one derives
\begin{equation}
\begin{array}{lll}
\mathbb{P}\big(\mathcal{Y} = x \big) = \mathbb{P}\big(\mathcal{X} = x \big), \, \forall x \in \mathcal{L}^{n}
\end{array}
\nonumber
\end{equation}
thus establishing the proof of Lemma \ref{lemma:besag}.
\end{proof}



\subsection{Proof of Lemma \ref{lemma:equiv_pi_ki}}
\label{app:proof_lemma_equiv_pi_ki}

\begin{proof}
Suppose $i \in \mathcal{V}$, and let us show the following equivalence:
\begin{equation}
\begin{array}{l}
\forall u \in \mathscr{F}, u(x) -  \pi_{i}(x_{{\ne}_{i}^+})  \sum_{x_i \in \mathcal{L}} u(x) = 0, \forall x \in \mathcal{L}^n
\Leftrightarrow  u(x) - \widehat{\pi}_i(x_{{\ne}_{i}^+}) \sum_{x_i\in \mathcal{L}} \widehat{\pi}_i(x_{{\ne}_{i}^+}) \, u(x)  = 0,  \forall x \in \mathcal{L}^n
\end{array}
\label{equiv_mrf_id_hat}
\end{equation}
Thus, assume an arbitrary function $u \in \mathscr{F}$, and let us first assume that $u$ satisfies
\begin{equation}
\begin{array}{l}
u(x) =  \pi_{i}(x_{{\ne}_{i}^+}) \, \sum_{x_i \in \mathcal{L}} u(x) , \, \forall x \in \mathcal{L}^n
\end{array}
\label{u_identity_pi}
\end{equation}
By using formula (\ref{u_identity_pi}), and the following two identities
\begin{equation}
\begin{array}{lll}
\widehat{\pi}_i(x_{{\ne}_{i}^+})  = \frac{{\pi}_i(x_{{\ne}_{i}^+})}{\sqrt{\sum_{x_i\in \mathcal{L}} {\pi}_i^2(x_{{\ne}_{i}^+})}} \, , \, \forall x_{{\ne}_{i}^+} \in \mathcal{L}^{|{\ne}_{i}^+|}
\end{array} \nonumber
\end{equation}
and
\begin{equation}
\begin{array}{lll}
\sum_{x_i \in \mathcal{L}} \widehat{\pi}_i^2(x_{{\ne}_{i}^+}) = 1, \, \forall x_{{\ne}_{i}} \in \mathcal{L}^{|{\ne}_{i}|}
\end{array} \nonumber
\end{equation}
one derives $\forall x \in \mathcal{L}^n$
\begin{multline}
u(x) - \widehat{\pi}_i(x_{{\ne}_{i}^+}) \, \sum_{x_i\in \mathcal{L}} \widehat{\pi}_i(x_{{\ne}_{i}^+}) \, u(x)  = \pi_{i}(x_{{\ne}_{i}^+}) \, \sum_{x_i \in \mathcal{L}} u(x) -
 \widehat{\pi}_i(x_{{\ne}_{i}^+}) \, \sum_{x_i\in \mathcal{L}} \widehat{\pi}_i(x_{{\ne}_{i}^+}) \, \Big( \pi_{i}(x_{{\ne}_{i}^+}) \, \sum_{x_i \in \mathcal{L}} u(x) \Big) \\
= \pi_{i}(x_{{\ne}_{i}^+}) \, \sum_{x_i \in \mathcal{L}} u(x) -
 \widehat{\pi}_i(x_{{\ne}_{i}^+}) \, \sum_{x_i\in \mathcal{L}} \widehat{\pi}_i(x_{{\ne}_{i}^+}) \, \Big( \sqrt{\sum_{x_i\in \mathcal{L}} {\pi}_i^2(x_{{\ne}_{i}^+})} \, \widehat{\pi}_i(x_{{\ne}_{i}^+}) \, \sum_{x_i \in \mathcal{L}} u(x) \Big) \\
= \pi_{i}(x_{{\ne}_{i}^+}) \, \sum_{x_i \in \mathcal{L}} u(x) -  \sqrt{\sum_{x_i\in \mathcal{L}} {\pi}_i^2(x_{{\ne}_{i}^+})} \, \widehat{\pi}_i(x_{{\ne}_{i}^+}) \, \sum_{x_i\in \mathcal{L}} \widehat{\pi}_i^2(x_{{\ne}_{i}^+}) \, \sum_{x_i \in \mathcal{L}} u(x) \\
= \pi_{i}(x_{{\ne}_{i}^+}) \, \sum_{x_i \in \mathcal{L}} u(x) -  {\pi}_i(x_{{\ne}_{i}^+}) \, \sum_{x_i \in \mathcal{L}} u(x) = 0 \nonumber
\end{multline}
thus showing the right implication ($\Rightarrow$) in formula (\ref{equiv_mrf_id_hat}). Next, to show the left
implication ($\Leftarrow$) in the latter, assume now that $u$ satisfies
\begin{equation}
\begin{array}{l}
u(x) = \widehat{\pi}_i(x_{{\ne}_{i}^+}) \, \sum_{x_i\in \mathcal{L}} \widehat{\pi}_i(x_{{\ne}_{i}^+}) \, u(x), \, \forall x \in \mathcal{L}^n
\end{array}
\label{u_identity_pi_b}
\end{equation}
Then, by summing both hand sides of the latter identity across coordinate $x_i$, and by using the identity
\begin{equation}
\begin{array}{lll}
\sum_{x_i \in \mathcal{L}} {\pi}_i(x_{{\ne}_{i}^+}) = 1, \, \forall x_{{\ne}_{i}} \in \mathcal{L}^{|{\ne}_{i}|}
\end{array} \nonumber
\end{equation}
hence,
\begin{equation}
\begin{array}{lll}
\sum_{x_i \in \mathcal{L}} \widehat{\pi}_i(x_{{\ne}_{i}^+}) = \frac{1}{\sqrt{\sum_{x_i\in \mathcal{L}} {\pi}_i^2(x_{{\ne}_{i}^+})}}, \, \forall x_{{\ne}_{i}} \in \mathcal{L}^{|{\ne}_{i}|}
\end{array} \nonumber
\end{equation}
one derives
\begin{equation}
\begin{array}{lll}
\sum_{x_i \in \mathcal{L}} u(x) &=& \Big(\sum_{x_i \in \mathcal{L}} \widehat{\pi}_i(x_{{\ne}_{i}^+}) \Big) \, \sum_{x_i\in \mathcal{L}} \widehat{\pi}_i(x_{{\ne}_{i}^+}) \, u(x)\\
&=&  \frac{1}{\sqrt{\sum_{x_i\in \mathcal{L}} {\pi}_i^2(x_{{\ne}_{i}^+})}} \sum_{x_i\in \mathcal{L}} \widehat{\pi}_i(x_{{\ne}_{i}^+}) \, u(x), \, \forall x_{\neq i} \in \mathcal{L}^{n-1}
\end{array}
\nonumber
\end{equation}
hence,
\begin{equation}
\begin{array}{lll}
\sum_{x_i\in \mathcal{L}} \widehat{\pi}_i(x_{{\ne}_{i}^+}) \, u(x) &=&  {\sqrt{\sum_{x_i\in \mathcal{L}} {\pi}_i^2(x_{{\ne}_{i}^+})}} \, \sum_{x_i \in \mathcal{L}} u(x)
\end{array}
\label{u_identity_pi_d}
\end{equation}
in such a way that, by substituting in formula (\ref{u_identity_pi_b}) the expression $\sum_{x_i\in \mathcal{L}} \widehat{\pi}_i(x_{{\ne}_{i}^+}) \, u(x)$ with
the right hand side expression of the latter formula (\ref{u_identity_pi_d}), one derives
\begin{equation}
\begin{array}{lll}
u(x) &=& \widehat{\pi}_i(x_{{\ne}_{i}^+}) \, \Big ( {\sqrt{\sum_{x_i\in \mathcal{L}} {\pi}_i^2(x_{{\ne}_{i}^+})}}  \, \sum_{x_i \in \mathcal{L}} u(x)  \Big) \\
     &=& \Big( \widehat{\pi}_i(x_{{\ne}_{i}^+}) \,  {\sqrt{\sum_{x_i\in \mathcal{L}} {\pi}_i^2(x_{{\ne}_{i}^+})}} \Big) \,  \sum_{x_i \in \mathcal{L}} u(x)  \\
     &=& {\pi}_i(x_{{\ne}_{i}^+}) \,  \sum_{x_i \in \mathcal{L}} u(x), \, \forall x \in \mathcal{L}^n
\end{array}
\nonumber
\end{equation}
thus showing also the left implication ($\Leftarrow$) in formula (\ref{equiv_mrf_id_hat}), hence, the proof of Lemma \ref{lemma:equiv_pi_ki}.
\end{proof}

\subsection{Proof of Lemma \ref{lemma:unicity_mrf}}
\label{app:proof_lemma_unicity_mrf}

\begin{proof}
We will first show that MRF distribution $\mathbb{P}(\mathcal{X} = \cdot)$ is the unique solution of the following system of linear equations with a strict
positivity constraint:
\begin{equation}
\begin{array}{l}
\begin{cases}
p(x) -  \pi_{i}(x_{{\ne}_{i}^+}) \, \sum_{x_i \in \mathcal{L}} p(x) = 0, \, \forall x \in \mathcal{L}^n,\, \forall i \in \mathcal{V} \\
\sum_{x \in \mathcal{L}^n} p(x) = 1 \\
p(x)> 0, \, \forall x \in \mathcal{L}^n
\end{cases}
\end{array}
\label{lin_sys_unicity_mrf_strict}
\end{equation}
One firstly checks that MRF distribution $\mathbb{P}(\mathcal{X} = \cdot)$ is a solution to system (\ref{lin_sys_unicity_mrf_strict}).
Next, assume a random field $\mathcal{Y} = \big(\mathcal{Y}_i\big)_{i\in \mathcal{V}}$ of which joint distribution $\mathbb{P}(\mathcal{Y} = \cdot)$
satisfies system (\ref{lin_sys_unicity_mrf_strict}), and put $p(x) = \mathbb{P}(\mathcal{Y} = x), \forall x \in \mathcal{L}^n$. Since, by assumption,
one has $p(x)>0$, $\forall x \in \mathcal{L}^{n}$, one finds $\sum_{x_i \in \mathcal{L}} p(x)>0$,  $\forall x_{\neq i} \in \mathcal{L}^{n-1}$,
$\forall i \in \mathcal{V}$, in such a way, that one can write
\begin{equation}
\begin{array}{l}
\frac{p(x)}{\sum_{x_i \in \mathcal{L}} p(x)} =  \pi_{i}(x_{{\ne}_{i}^+}), \, \forall x \in \mathcal{L}^n,\, \forall i \in \mathcal{V}
\end{array}
\nonumber
\end{equation}
equivalently,
\begin{equation}
\begin{array}{l}
\mathbb{P} \big(\mathcal{Y}_i = x_i / \mathcal{Y}_{\neq i} = x_{\neq i} \big) =  \pi_{i}(x_{{\ne}_{i}^+}), \, \forall x \in \mathcal{L}^n,\, \forall i \in \mathcal{V}
\end{array}
\label{mrf_y_a_a}
\end{equation}
this is on the one hand. On the other hand, by marginalization of both hand sides of formula of the top-most identity of system (\ref{lin_sys_unicity_mrf_strict})
with regard to all the variables whose indices are in $\mathcal{V}$, but, not in ${\ne}_{i}^+$, one finds
\begin{equation}
\begin{array}{l}
\sum_{k \in \mathcal{V}/{\ne}_{i}^+} \sum_{x_k\in \mathcal{L}} \big( p(x) -  \pi_{i}(x_{{\ne}_{i}^+}) \, \sum_{x_i \in \mathcal{L}} p(x)\big)= 0,
\, \forall x_{{\ne}_{i}^+} \in \mathcal{L}^{|{\ne}_{i}^+|},\, \forall i \in \mathcal{V}
\end{array}
\nonumber
\end{equation}
hence,
\begin{equation}
\begin{array}{l}
p_{{\ne}_{i}^+}(x_{{\ne}_{i}^+}) -  \pi_{i}(x_{{\ne}_{i}^+}) \, p_{{\ne}_{i}}(x_{{\ne}_{i}}) = 0, \, \forall x_{{\ne}_{i}^+} \in \mathcal{L}^{|{\ne}_{i}^+|},\, \forall i \in \mathcal{V}
\end{array}
\label{mrf_y_b}
\end{equation}
where it has been assumed that $\forall i \in \mathcal{V}$, $p_{{\ne}_{i}^+}$, and $p_{{\ne}_{i}}$ are the margins of $p$ with respect to ${\ne}_{i}^+$ and ${\ne}_{i}$, respectively.
Since by assumption, one has $p(x)> 0, \, \forall x \in \mathcal{L}^n$,  one also has $p_{{\ne}_{i}}(x_{{\ne}_{i}})>0$,
$\forall x_{{\ne}_{i}} \in \mathcal{L}^{|{\ne}_{i}|},\, \forall i \in \mathcal{V}$, in such a way that, one can reexpress formula (\ref{mrf_y_b}) as
\begin{equation}
\begin{array}{l}
\frac{p_{{\ne}_{i}^+}(x_{{\ne}_{i}^+})}{p_{{\ne}_{i}}(x_{{\ne}_{i}})} = \pi_{i}(x_{{\ne}_{i}^+}), \, \forall x_{{\ne}_{i}^+} \in \mathcal{L}^{|{\ne}_{i}^+|},\, \forall i \in \mathcal{V}
\end{array}
\nonumber
\end{equation}
which in turn reexpresses as
\begin{equation}
\begin{array}{l}
\mathbb{P} \big(\mathcal{Y}_i = x_i/\mathcal{Y}_{{\ne}_{i}} = x_{{\ne}_{i}}\big) = \pi_{i}(x_{{\ne}_{i}^+}), \, \forall x_{{\ne}_{i}^+} \in \mathcal{L}^{|{\ne}_{i}^+|},\, \forall i \in \mathcal{V}
\end{array}
\label{mrf_y_b_b}
\end{equation}
Hence, by combining both formulas (\ref{mrf_y_a_a}) and (\ref{mrf_y_b_b}), one derives
\begin{equation}
\begin{array}{l}
\mathbb{P} \big(\mathcal{Y}_i = x_i / \mathcal{Y}_{\neq i} = x_{\neq i} \big) = \mathbb{P} \big(\mathcal{Y}_i = x_i/\mathcal{Y}_{{\ne}_{i}} = x_{{\ne}_{i}}\big), \, \forall x \in \mathcal{L}^n, \, \forall i \in \mathcal{V}
\end{array}
\nonumber 
\end{equation}
moreover, since $\mathbb{P}(\mathcal{Y} = \cdot)$ is assumed to be a strictly positive joint distribution, one concludes that
$\mathbb{P}(\mathcal{Y} = \cdot)$ is a MRF distribution with regard to graph $\mathcal{G}$.

Finally, by formula (\ref{mrf_y_b_b}), all the individual conditional distributions of MRF distribution $\mathbb{P}(\mathcal{Y} = \cdot)$
coincide with the respective individual conditional distributions of MRF distribution $\mathbb{P}(\mathcal{X} = \cdot)$, in such a way that,
by Lemma \ref{lemma:besag}, one establishes that $\mathbb{P}(\mathcal{Y} = \cdot) \equiv \mathbb{P}(\mathcal{X} = \cdot)$, implying
that $\mathbb{P}(\mathcal{X} = \cdot)$ is the unique solution of system (\ref{lin_sys_unicity_mrf_strict}).

Next, let us show that system (\ref{lin_sys_unicity_mrf_strict}) (with a strict positivity constraint) is equivalent to the
following one (with a nonnegativity constraint):
\begin{equation}
\begin{array}{l}
\begin{cases}
p(x) -  \pi_{i}(x_{{\ne}_{i}^+}) \, \sum_{x_i \in \mathcal{L}} p(x) = 0, \, \forall x \in \mathcal{L}^n,\, \forall i \in \mathcal{V} \\
\sum_{x \in \mathcal{L}^n} p(x) = 1 \\
p(x) \geq 0, \, \forall x \in \mathcal{L}^n
\end{cases}
\end{array}
\label{lin_sys_unicity_mrf_soft}
\end{equation}
For this purpose, we will first assume that the opposite of the latter claim is true, i.e.; there exists  $p \in \mathscr{F}$ satisfying
system (\ref{lin_sys_unicity_mrf_soft}), and, there exists $x^0 \in \mathcal{L}^n$, such that, $p(x^0) = 0$, and we will show that
this leads to a contradiction. Indeed, denote
\begin{equation}
\begin{array}{l}
\pi_0 = \inf_{i \in \mathcal{V}, y_{{\ne}_{i}^+} \in \mathcal{L}^{|{\ne}_{i}^+|}} \Big\{ \big\{ \pi_{i}(x_{{\ne}_{i}^+})  \big \} \, \Big\}
\end{array}
\nonumber
\end{equation}
which, by assumption, is a strictly positive number. From the top row of system (\ref{lin_sys_unicity_mrf_soft}), one finds
\begin{equation}
\begin{array}{lll}
\sum_{x_i \in \mathcal{L}} p(x) &=&  p(x) \, \frac{1}{\pi_{i}(x_{{\ne}_{i}^+})}\\
                                &\leq& p(x) \, \sup_{i \in \mathcal{V}, y_{{\ne}_{i}^+} \in \mathcal{L}^{|{\ne}_{i}^+|}} \Big\{ \big\{ \frac{1}{\pi_{i}(y_{{\ne}_{i}^+})}  \big \} \, \Big\} \\
                                &=& p(x) \, \frac{1}{\inf_{i \in \mathcal{V}, y_{{\ne}_{i}^+} \in \mathcal{L}^{|{\ne}_{i}^+|}} \Big\{ \big\{ {\pi_{i}(y_{{\ne}_{i}^+})}  \big \} \, \Big\}} \\
                                 &=& \frac{1}{\pi_0} \, p(x), \, \forall x \in \mathcal{L}^n, \, \forall i \in \mathcal{V}
\end{array}
\label{px_pi_px_i_bound}
\end{equation}

Thus, by applying formula (\ref{px_pi_px_i_bound}) with $x=x_0$, one finds $\forall i \in \mathcal{V}$
\begin{equation}
\begin{array}{lll}
\sum_{x_i \in \mathcal{L}} p(x_1^0,\dots,x_{i-1}^0,x_i,x_{i+1}^0,\dots,x_n^0) &\leq& \frac{1}{\pi_0} \, p(x_0) \\
&=& 0
\end{array}
\nonumber
\end{equation}
implying
\begin{equation}
\begin{array}{l}
p(x_1^0,\dots,x_{i-1}^0,x_i,x_{i+1}^0,\dots,x_n^0) = 0, \, \forall x_i \in \mathcal{L}, \, \forall i\in  \mathcal{V}
\label{form_y_zero_b}
\end{array}
\end{equation}
as by assumption, $p(x_1^0,\dots,x_{i-1}^0,x_i,x_{i+1}^0,\dots,x_n^0)\geq 0, \forall i \in \mathcal{V}$, $\forall x_i \in \mathcal{L}$.

Now, suppose $k \in \mathcal{V}$, and $x_k \in \mathcal{L}$. Hence, by formula (\ref{form_y_zero_b}), one has
\begin{equation}
\begin{array}{l}
p(x_1^0,\dots,x_{k-1}^0,x_k,x_{k+1}^0,\dots,x_n^0) = 0
\nonumber
\end{array}
\end{equation}
in such a way that, by applying again formula (\ref{px_pi_px_i_bound}) with $x=(x_1^0,\dots,x_{k-1}^0,x_k,x_{k+1}^0,\dots,x_n^0)$ this time,
and since both $k$ and $x_k$ are assumed to be arbitrary, one derives
\begin{equation}
\begin{array}{l}
p(x_1^0, \dots, x_{i-1}^{0}, x_{i+1}^{0}, \dots, x_{k-1}^{0}, x_{k}, x_{k+1}^{0}, \dots,x_n^0 ) = 0, \, \forall x_i, x_k \in \mathcal{L}, \, \forall i,k \in  \mathcal{V} \text{ s.t.; } i<k
\nonumber
\end{array}
\end{equation}
By repeating this process as many times as necessary, one finally derives
\begin{equation}
\begin{array}{l}
p(x_1,\dots, x_n) = 0, \, \forall x_{1},\dots, x_{n} \in \mathcal{L}
\nonumber
\end{array}
\end{equation}
implying that $\sum_{x \in \mathcal{L}^n} p(x) = 0$ which is an absurdity, as by assumption, $p$ satisfies system (\ref{lin_sys_unicity_mrf_strict}),
hence, satisfies, in particular, that $\sum_{x \in \mathcal{L}^n} p(x) = 1$.

One concludes that systems (\ref{lin_sys_unicity_mrf_strict}) and (\ref{lin_sys_unicity_mrf_soft}) are equivalent.
Finally, equivalence between system (\ref{lin_sys_unicity_mrf_soft}) and the following one:
\begin{equation}
\begin{cases}
p(x) - \widehat{\pi}_i(x_{{\ne}_{i}^+}) \, \sum_{x_i\in \mathcal{L}} \widehat{\pi}_i(x_{{\ne}_{i}^+}) \,  p(x) = 0, \, \forall x \in \mathcal{L}^n, \, \forall i \in \mathcal{V} \\
\sum_{x \in \mathcal{L}^{n}} p(x) = 1 \\
p(x) \geq 0, \, \forall x \in \mathcal{L}^n
\label{lin_sys_unicity_mrf_soft_hat}
\end{cases}
\end{equation}
follows immediately from Lemma \ref{lemma:equiv_pi_ki}, hence, the proof of Lemma \ref{lemma:unicity_mrf}.
\end{proof}


\subsection{Proof of Lemma \ref{lemma:max_inf_min_sup_vi}}
\label{subap:proof_lemma_max_inf_min_sup_vi}

\begin{proof}
Assume an arbitrary function $f \in \mathscr{F}$. By Lemma \ref{lemma:unicity_mrf}, one finds that MRF distribution $\mathbb{P}\big(\mathcal{X} = \cdot\big)$ is the unique
feasible solution of following both LPs:
\begin{equation}
\begin{array}{l}
\underset{p}{\min} \Big\{ \, \sum_{x \in \mathcal{L}^n} f(x) \, p(x)  \, \Big\} \\
\begin{cases}
p(x) -  \widehat{\pi}_i(x_{{\ne}_{i}^+}) \, \sum_{x_i\in \mathcal{L}} \widehat{\pi}_i(x_{{\ne}_{i}^+}) \,  p(x) = 0, \, \forall x \in \mathcal{L}^n,\, \forall i \in \mathcal{V} \\
\sum_{x \in \mathcal{L}^n} p(x) = 1 \\
p(x)\geq 0 , \, \forall x \in \mathcal{L}^n
\end{cases}
\end{array}
\label{lin_prog_unicity_mrf_lemma}
\end{equation}
and
\begin{equation}
\begin{array}{l}
\underset{p}{\max} \Big\{ \,  \sum_{x \in \mathcal{L}^n} f(x) \, p(x)  \, \Big\} \\
\begin{cases}
p(x) - \widehat{\pi}_i(x_{{\ne}_{i}^+}) \, \sum_{x_i\in \mathcal{L}} \widehat{\pi}_i(x_{{\ne}_{i}^+}) \,  p(x) = 0, \, \forall x \in \mathcal{L}^n,\, \forall i \in \mathcal{V} \\
\sum_{x \in \mathcal{L}^n} p(x) = 1 \\
p(x)\geq 0 , \, \forall x \in \mathcal{L}^n
\end{cases}
\end{array}
\label{lin_prog_unicity_mrf_lemma_max}
\end{equation}
One concludes that both LP (\ref{lin_prog_unicity_mrf_lemma}) and LP (\ref{lin_prog_unicity_mrf_lemma_max}) achieve
an optimal objective value equal to
\begin{equation}
\begin{array}{l}
\sum_{x \in \mathcal{L}^n} f(x) \, \mathbb{P}\big(\mathcal{X} = x \big)  = \mathbb{E}\big[f(\mathcal{X})\big]
\end{array}
\label{expect_f_optimal}
\end{equation}
This is on the one hand. On the other hand, one first finds that the dual of LP (\ref{lin_prog_unicity_mrf_lemma}) writes as
\begin{equation}
\begin{array}{lll}
\underset{\mu, \, v_i\mathscr{F}, \forall i \in \mathcal{V}}{\max} \big\{\,  \mu \, \big\} \\
\begin{cases}
\mu \leq f(x) + \sum_{i \in  \mathcal{V}} \big(\, v_i(x) - \widehat{\pi}_i(x_{{\ne}_{i}^+}) \, \sum_{x_i\in \mathcal{L}} \widehat{\pi}_i(x_{{\ne}_{i}^+}) \,  v_i(x) \, \big), \, \forall x\in \mathcal{L}^n \\
v_i: \mathcal{L}^n \leftarrow \mathcal{R}, \, \forall i \in \mathcal{V}\\
\mu  \in \mathcal{R}
\end{cases}
\end{array}
\label{lin_prog_unicity_mrf_lemma_dual}
\end{equation}
One finds that the optimal value of $\mu$ in dual LP (\ref{lin_prog_unicity_mrf_lemma_dual}) in terms of functions $f$ and $v_i, \forall i \in \mathcal{V}$ is given by the expression
\begin{equation}
\begin{array}{lll}
\underset{x \in \mathcal{L}^n}{\inf}  \big\{ \, f(x) + \sum_{i \in  \mathcal{V}} \big(\, v_i(x) - \widehat{\pi}_i(x_{{\ne}_{i}^+}) \, \sum_{x_i\in \mathcal{L}} \widehat{\pi}_i(x_{{\ne}_{i}^+}) \,  v_i(x) \, \big)\, \big\}
\end{array}
\nonumber
\end{equation}
One rewrites dual LP (\ref{lin_prog_unicity_mrf_lemma_dual}) accordingly as follows:
\begin{equation}
\begin{array}{lll}
 \underset{ v_i\mathscr{F}, \forall i \in \mathcal{V} }{\max}
\, \underset{x \in \mathcal{L}^n}{\inf} \Big\{ \,
f(x) + \sum_{i \in  \mathcal{V}} \big( \, v_i(x) - \widehat{\pi}_i(x_{{\ne}_{i}^+}) \, \sum_{x_i\in \mathcal{L}} \widehat{\pi}_i(x_{{\ne}_{i}^+}) \,  v_i(x)\,\big)\, \Big\}
\end{array}
\label{mrf_expect_as_min_dual_s_b}
\end{equation}
One proceeds in a similar way, by assuming LP (\ref{lin_prog_unicity_mrf_lemma_max}) instead of LP (\ref{lin_prog_unicity_mrf_lemma}), for establishing that
the dual of LP (\ref{lin_prog_unicity_mrf_lemma_max}) can write as follows:
\begin{equation}
\begin{array}{lll}
\underset{v_i\mathscr{F}, \forall i \in \mathcal{V}}{\min} \, \underset{x \in \mathcal{L}^n}{\sup}
\Big\{  \, f(x) + \sum_{i \in  \mathcal{V}} \big(\, v_i(x) - \widehat{\pi}_i(x_{{\ne}_{i}^+}) \, \sum_{x_i\in \mathcal{L}} \widehat{\pi}_i(x_{{\ne}_{i}^+}) \,  v_i(x)\, \big)\, \Big\}
\end{array}
\label{mrf_expect_as_max_dual_s_c}
\end{equation}
Since both dual LPs (\ref{mrf_expect_as_min_dual_s_b}) and  (\ref{mrf_expect_as_max_dual_s_c}) achieve the same optimal objective value as their respective
primal LPs (\ref{lin_prog_unicity_mrf_lemma}) and LP (\ref{lin_prog_unicity_mrf_lemma_max}) which we have already established in formula (\ref{expect_f_optimal}),
one finally derives
\begin{equation}
\begin{array}{lll}
\mathbb{E}\big[f(\mathcal{X})\big]  &=& \underset{v_i\mathscr{F}, \forall i \in \mathcal{V} }{\max}
\, \underset{x \in \mathcal{L}^n}{\inf} \Big\{ \,
f(x) + \sum_{i \in  \mathcal{V}} \big( \, v_i(x) -\widehat{\pi}_i(x_{{\ne}_{i}^+}) \, \sum_{x_i\in \mathcal{L}} \widehat{\pi}_i(x_{{\ne}_{i}^+}) \,  v_i(x)\,\big)\, \Big\} \\
&=& \underset{v_i\mathscr{F}, \forall i \in \mathcal{V}}{\min} \, \underset{x \in \mathcal{L}^n}{\sup}
\Big\{  \, f(x) + \sum_{i \in  \mathcal{V}} \big(\, v_i(x) - \widehat{\pi}_i(x_{{\ne}_{i}^+}) \, \sum_{x_i\in \mathcal{L}} \widehat{\pi}_i(x_{{\ne}_{i}^+}) \,  v_i(x)\, \big)\, \Big\}
\end{array}
\nonumber 
\end{equation}
hence, the proof of Lemma \ref{lemma:max_inf_min_sup_vi}.
\end{proof}


\subsection{Proof of Lemma \ref{lemma:founding_result_mrf}}
\label{app:proof_lemma_founding_result_mrf}

\begin{proof}
Suppose an arbitrary function $f \in \mathscr{F}$, and assume the following optimization problem:
\begin{equation}
\begin{array}{lll}
\underset{\big\{ v_i\in \mathscr{F}, \forall i \in \mathcal{V} \big\}}{\max} \,
\underset{x \in \mathcal{L}^n}{\inf} \Big\{ f(x) + \sum_{i \in \mathcal{V}} \Big(\, v_i(x) - \widehat{\pi}_i(x_{{\ne}_{i}^+}) \, \sum_{x_i\in \mathcal{L}} \widehat{\pi}_i(x_{{\ne}_{i}^+}) \, v_i(x) \, \Big) \Big\}
\end{array}
\label{mrf_expect_min_sup_a}
\end{equation}
Then, by Lemma \ref{lemma:max_inf_min_sup_vi}, optimization problem (\ref{mrf_expect_min_sup_a}) is bounded, and its optimal solution
denoted by $(v_i^0)_{i \in \mathcal{V}}$ satisfies
\begin{equation}
\begin{array}{lll}
\mathbb{E}\big[f(\mathcal{X})\big] &=&  \underset{x \in \mathcal{L}^n}{\inf} \Big\{ f(x) + \sum_{i \in \mathcal{V}} \big(\, v_i^0(x) -
\widehat{\pi}_i(x_{{\ne}_{i}^+}) \, \sum_{x_i\in \mathcal{L}} \widehat{\pi}_i(x_{{\ne}_{i}^+}) \,  v_i^0(x) \, \big) \Big\}
\end{array}
\nonumber
\end{equation}
hence,
\begin{equation}
\begin{array}{lll}
\underset{x \in \mathcal{L}^n}{\inf} \Big\{ f(x) - \mathbb{E}\big[f(\mathcal{X})\big] + \sum_{i \in \mathcal{V}}  \big( v_i^0(x) - \widehat{\pi}_i(x_{{\ne}_{i}^+}) \, \sum_{x_i\in \mathcal{L}} \widehat{\pi}_i(x_{{\ne}_{i}^+}) \,  v_i^0(x) \, \big) \Big\} = 0
\end{array}
\label{mrf_expect_max_inf_c}
\end{equation}
Next, denote
\begin{equation}
\begin{array}{lll}
r(x) &=&  f(x)-\mathbb{E}\big[f(\mathcal{X})\big] + \sum_{i \in \mathcal{V}} \big( v_i^0(x) - \widehat{\pi}_i(x_{{\ne}_{i}^+}) \, \sum_{x_i\in \mathcal{L}} \widehat{\pi}_i(x_{{\ne}_{i}^+}) \,  v_i^0(x) \, \big), \, \forall x \in \mathcal{L}^n
\end{array}
\nonumber
\end{equation}
Then, on the one hand, by using formula (\ref{mrf_expect_max_inf_c}), one derives
\begin{equation}
\begin{array}{lll}
r(x) \geq 0, \forall x \in \mathcal{L}^n
\end{array}
\label{tilde_h_nonngeative}
\end{equation}
On the other hand, one has
\begin{enumerate}
\item
$\sum_{x \in \mathcal{L}^n} \Big( f(x)-\mathbb{E}\big[f(\mathcal{X})\big] \Big)  \, \mathbb{P}\big(\mathcal{X} = x\big) = \sum_{x \in \mathcal{L}^n} f(x) \, \mathbb{P}\big(\mathcal{X} = x\big) - \mathbb{E}\big[f(\mathcal{X})\big] = 0$
\item
$\sum_{x \in \mathcal{L}^n} \sum_{i \in \mathcal{V}} \Big( v_i^0(x) - \widehat{\pi}_i(x_{{\ne}_{i}^+}) \, \sum_{x_i\in \mathcal{L}} \widehat{\pi}_i(x_{{\ne}_{i}^+}) \,  v_i^0(x) \Big) \, \mathbb{P}\big(\mathcal{X} = x\big) = \sum_{i \in \mathcal{V}} \langle k_i \, v_i^0, \mathbb{P}\big(\mathcal{X} = \cdot \big) \rangle\nonumber
= \sum_{i \in \mathcal{V}} \langle  v_i^0, k_i \, \mathbb{P}\big(\mathcal{X} = \cdot \big) \rangle  = 0
$
\end{enumerate}
implying
\begin{equation}
\begin{array}{lll}
\mathbb{E}\big[r(\mathcal{X})\big] &=& \sum_{x \in \mathcal{L}^n} r(x) \, \mathbb{P}\big(\mathcal{X} = x\big) \\
 &=&  \sum_{x \in \mathcal{L}^n} \Big( f(x)-\mathbb{E}\big[f(\mathcal{X})\big] \Big)  \, \mathbb{P}\big(\mathcal{X} = x\big)   + \sum_{i \in \mathcal{V}}  \sum_{x \in \mathcal{L}^n} \Big( v_i^0(x) - \widehat{\pi}_i(x_{{\ne}_{i}^+}) \, \sum_{x_i\in \mathcal{L}} \widehat{\pi}_i(x_{{\ne}_{i}^+}) \,  v_i^0(x) \Big) \, \mathbb{P}\big(\mathcal{X} = x\big) \\
 &=& 0
 \end{array}
\label{expec_h_tilde}
\end{equation}
in such a way that, by combining both formulas (\ref{tilde_h_nonngeative}) and (\ref{expec_h_tilde}), and by
using the assumption that $\mathbb{P}\big(\mathcal{X} = x\big)>0$, $\forall x \in \mathcal{L}^n$, one establishes
\begin{equation}
\begin{array}{lll}
r(x) = f(x)-\mathbb{E}\big[f(\mathcal{X})\big] + \sum_{i \in \mathcal{V}} \big( v_i^0(x) - \widehat{\pi}_i(x_{{\ne}_{i}^+}) \, \sum_{x_i\in \mathcal{L}} \widehat{\pi}_i(x_{{\ne}_{i}^+}) \,  v_i^0(x) \, \big) = 0, \, \forall x \in \mathcal{L}^n
\end{array}
\nonumber
\end{equation}
equivalently
\begin{equation}
\begin{array}{lll}
\mathbb{E}\big[f(\mathcal{X})\big] = f(x) + \sum_{i \in \mathcal{V}} \big( v_i^0(x) -  \widehat{\pi}_i(x_{{\ne}_{i}^+}) \, \sum_{x_i\in \mathcal{L}} \widehat{\pi}_i(x_{{\ne}_{i}^+}) \, v_i^0(x) \, \big), \, \forall x \in \mathcal{L}^n
\end{array}
\nonumber
\end{equation}
as otherwise, i.e.; if there existed $x \in \mathcal{L}^n$, such that, $r(x)>0$ then one would find $\mathbb{E}\big[r(\mathcal{X})\big]>0$ which is an absurdity.
Then, by putting $v_i \equiv -v_i^0$, $\forall i \in \mathcal{V}$, one finds
\begin{equation}
\begin{array}{lll}
f(x)  = \mathbb{E}\big[f(\mathcal{X})\big] +  \sum_{i \in \mathcal{V}} \big( v_i(x) -  \widehat{\pi}_i(x_{{\ne}_{i}^+}) \, \sum_{x_i\in \mathcal{L}} \widehat{\pi}_i(x_{{\ne}_{i}^+}) \, v_i(x) \, \big), \, \forall x \in \mathcal{L}^n
\end{array}
\nonumber
\end{equation}
Finally, since $f$ is assumed to be arbitrary, one concludes that $\forall f \in \mathscr{F}$, there exists $n$ functions $v_i \in \mathscr{F}$, $\forall i \in \mathcal{V}$, such that, one can write
\begin{equation}
\begin{array}{lll}
f(x)  = \mathbb{E}\big[f(\mathcal{X})\big] +  \sum_{i \in \mathcal{V}} \big( v_i(x) -  \widehat{\pi}_i(x_{{\ne}_{i}^+}) \, \sum_{x_i\in \mathcal{L}} \widehat{\pi}_i(x_{{\ne}_{i}^+}) \, v_i(x) \, \big),
\, \forall x \in \mathcal{L}^n
\end{array}
\nonumber
\end{equation}
hence, the proof of Lemma \ref{lemma:founding_result_mrf}.
\end{proof}


\subsection{Proof of Lemma \ref{lemma:mrf_syst_nondegenerate}}
\label{app:proof_lemma_mrf_syst_nondegenerate}

\begin{proof}
Assume the opposite, that is, the following system:
\begin{equation}
\begin{array}{lll}
\begin{cases}
\lambda(x) - \widehat{\pi}_i(x_{{\ne}_{i}^+}) \, \sum_{x_i\in \mathcal{L}} \widehat{\pi}_i(x_{{\ne}_{i}^+}) \,  \lambda(x) = 0, \, \forall x \in \mathcal{L}^n, \, \forall i \in \mathcal{V} \\
\sum_{x \in \mathcal{L}^n} \lambda(x) = 0
\end{cases}
\end{array}
\label{lin_comb_quad_cond_prf}
\end{equation}
has a non-zero solution denoted by $\lambda_1$, implying that $\sup_{x \in \mathcal{L}^n} \big\{ \, \big | \lambda_1(x) \big|  \, \big\}>0$.
Furthermore, put
\begin{equation}
\begin{cases}
a_0 = \inf_{x \in \mathcal{L}^n} \big\{ \mathbb{P}\big(\mathcal{X} = x \big) \big\}\\
b_0 = \sup_{x \in \mathcal{L}^n} \big\{ \, \big | \lambda_1(x) \big|  \, \big\} \\
q \equiv   \mathbb{P}\big(\mathcal{X} = \cdot \big) +  \frac{\alpha_0}{b_0} \, \lambda_1
\end{cases}
\nonumber
\end{equation}
Since by assumption, $a_0 >0$, and $\lambda_1 \not \equiv \zeroindic$, one derives $\frac{\alpha_0}{b_0} \, \lambda_1 \not \equiv \zeroindic$,
implying that $q \not \equiv  \mathbb{P}\big(\mathcal{X} = x \big)$, this is on the one hand. On the other hand, one easily checks that $q$
satisfies
\begin{equation}
\begin{cases}
q(x) - \widehat{\pi}_i(x_{{\ne}_{i}^+}) \, \sum_{x_i\in \mathcal{L}} \widehat{\pi}_i(x_{{\ne}_{i}^+}) \,  q(x) = 0, \, \forall x \in \mathcal{L}^n, \, \forall i \in \mathcal{V} \\
\sum_{x \in \mathcal{L}^{n}} q(x) = 1 \\
q(x) \geq 0, \, \forall x \in \mathcal{L}^n
\nonumber
\end{cases}
\end{equation}
which is an absurdity, as by Lemma \ref{lemma:unicity_mrf}, $ \mathbb{P}\big(\mathcal{X} = \cdot \big)$ is the unique solution of the following system:
\begin{equation}
\begin{cases}
p(x) - \widehat{\pi}_i(x_{{\ne}_{i}^+}) \, \sum_{x_i\in \mathcal{L}} \widehat{\pi}_i(x_{{\ne}_{i}^+}) \,  p(x) = 0, \, \forall x \in \mathcal{L}^n, \, \forall i \in \mathcal{V} \\
\sum_{x \in \mathcal{L}^{n}} p(x) = 1 \\
p(x) \geq 0, \, \forall x \in \mathcal{L}^n
\nonumber
\end{cases}
\end{equation}
One concludes that $\lambda \equiv \zeroindic$ is the unique solution of system (\ref{lin_comb_quad_cond_prf}), hence, the proof of Lemma \ref{lemma:mrf_syst_nondegenerate}.
\end{proof}

\subsection{Proof of Lemma \ref{lemma:markov_lin_oper}}
\label{app:proof_lemma_markov_lin_oper}

\begin{proof}
By Proposition \ref{prop:o_s_ortho}, $\mathscr{O}_{\emptyset}$ is an orthogonal projection, hence, $L^{n}\, \mathscr{O}_{\emptyset}$ is a
self-adjoint positive semi-definite linear operator. Furthermore, by Proposition \ref{prop:ki_is_ortho_project},
$\forall i \in \mathcal{V}$, $k_i$ is an orthogonal projection, hence a positive semi-definite linear operator.
Then, by Propositions \ref{prop:K_k_i}, $\mathscr{K}$ writes as a sum of self-adjoint positive semi-definite
linear operators, implying that $\mathscr{K}$ is, in turn, a self-adjoint positive semi-definite linear operator.

Thus, it remains to show that $\mathscr{K}$ is invertible. For this purpose, suppose an arbitrary function
$f \in \mathscr{F}$, and assume the following quadratic program (QP):
\begin{equation}
\begin{array}{lll}
\underset{a \in \mathbb{R}, v_i\in \mathscr{F}, \forall i \in \mathcal{V}}{\min} \Big\{ \frac{a ^2}{2} + \sum_{i \in \mathcal{V}}  \frac{\| v_i \|^2}{2} \, \Big\} \\
\text{s.t., } a + \sum_{i \in \mathcal{V}} \big( \, v_i(x) - \widehat{\pi}_i(x_{{\ne}_{i}^+}) \, \sum_{x_i\in \mathcal{L}} \widehat{\pi}_i(x_{{\ne}_{i}^+})\, v_i(x) \, \big) = f(x), \, \forall x \in \mathcal{L}^n
\end{array}
\label{quad_prog_uniq}
\end{equation}
Then, by Lemma \ref{lemma:founding_result_mrf}, there exists $n$ functions $v_i^0 \in \mathscr{F}, \forall i \in \mathcal{V}$, such that, one can write
\begin{equation}
\begin{array}{lll}
f(x) = \mathbb{E}[f(\mathcal{X})] + \sum_{i \in \mathcal{V}} \big( \, v_i^0(x) - \widehat{\pi}_i(x_{{\ne}_{i}^+}) \, \sum_{x_i\in \mathcal{L}} \widehat{\pi}_i(x_{{\ne}_{i}^+})\, v_i^0(x) \, \big), \, \forall x \in \mathcal{L}^n
\end{array}
\nonumber
\end{equation}
meaning that, at least, the solution which is given by $a = \mathbb{E}[f(\mathcal{X})]$, and $v_i \equiv v_i^0, \forall i \in \mathcal{V}$
is a feasible solution of QP (\ref{quad_prog_uniq}), implying that QP (\ref{quad_prog_uniq}) is feasible. Second, QP (\ref{quad_prog_uniq}) is
convex, since its objective function is a positive definite quadratic form. Third, the Lagrangian functional of QP (\ref{quad_prog_uniq})
with $\lambda \in \mathscr{F}$ as a Lagrange multiplier function writes $\forall a \in \mathbb{R}$, $\forall  v_i \in \mathscr{F}$,
$\forall i \in \mathcal{V}$, and $\forall \lambda \in \mathscr{F}$ as
\begin{equation}
\begin{array}{lll}
\mathscr{L}_{\lambda}\big(a,(v_i)_{i \in \mathcal{V}}\big) &=& \frac{a ^2}{2} + \sum_{i \in \mathcal{V}}  \frac{\| v_i \|^2}{2} -
\sum_{x \in \mathcal{L}^n} \Big( a + \sum_{i \in \mathcal{V}} \big( \, v_i(x) - \widehat{\pi}_i(x_{{\ne}_{i}^+}) \, \sum_{x_i\in \mathcal{L}} \widehat{\pi}_i(x_{{\ne}_{i}^+}) \, v_i(x) \, \big) \Big) \, \lambda(x) \\
&=& \frac{a ^2}{2} + \sum_{i \in \mathcal{V}}  \frac{\| v_i \|^2}{2}
- \big(\sum_{x \in \mathcal{L}^n}  \lambda(x) \big) \, a - \sum_{i \in \mathcal{V}} \sum_{x \in \mathcal{L}^n} \Big( \, \lambda(x) - \widehat{\pi}_i(x_{{\ne}_{i}^+}) \, \sum_{x_i\in \mathcal{L}}
\widehat{\pi}_i(x_{{\ne}_{i}^+}) \,  \lambda(x) \, \Big) \,  v_i(x)
\end{array}
\nonumber
\end{equation}
Moreover, one has
\begin{equation}
\begin{cases}
\frac{\partial \mathscr{L}_{\lambda}\big(a,(v_i)_{i \in \mathcal{V}}\big)}{\partial a} = a -  \big(\sum_{x \in \mathcal{L}^n}  \lambda(x) \big)\\
\frac{\partial \mathscr{L}_{\lambda}\big(a,(v_i)_{i \in \mathcal{V}}\big)}{\partial v_i} = v_i(x) - \big(\lambda(x) -  \widehat{\pi}_i(x_{{\ne}_{i}^+}) \,
\sum_{x_i\in \mathcal{L}} \widehat{\pi}_i(x_{{\ne}_{i}^+}) \, \lambda(x)\big), \,  \forall x \in \mathcal{L}^n, \, \forall i \in \mathcal{V}
\end{cases}
\nonumber
\end{equation}
Suppose $\big(a^*,(v_i^*)_{\forall i \in \mathcal{V}}\big)$ is an optimal solution of QP (\ref{quad_prog_uniq}). Thus, by Karush-Kuhn-Tucker Theorem,
one finds that there exists a Lagrange multiplier function $\lambda \in \mathscr{F}$, such that, one has
\begin{equation}
\begin{cases}
\frac{\partial \mathscr{L}_{\lambda}\big(a^*,(v_i^*)_{i \in \mathcal{V}}\big)}{\partial a} = 0\\
\frac{\partial \mathscr{L}_{\lambda}\big(a^*,(v_i^*)_{i \in \mathcal{V}}\big)}{\partial v_i} = 0, \,  \forall x \in \mathcal{L}^n, \, \forall i \in \mathcal{V}
\end{cases}
\nonumber
\end{equation}
equivalently, there exists $\lambda \in \mathscr{F}$, such that, one has
\begin{equation}
\begin{cases}
a^* -  \big(\sum_{x \in \mathcal{L}^n}  \lambda(x) \big) = 0\\
v_i^*(x) - \big(\lambda(x) -  \widehat{\pi}_i(x_{{\ne}_{i}^+}) \, \sum_{x_i\in \mathcal{L}} \widehat{\pi}_i(x_{{\ne}_{i}^+}) \, \lambda(x)\big) = 0, \,  \forall x \in \mathcal{L}^n, \, \forall i \in \mathcal{V}
\end{cases}
\nonumber
\end{equation}
One finds that there exists $\lambda \in \mathscr{F}$, such that, one can write
\begin{equation}
\begin{cases}
a^* = \sum_{x \in \mathcal{L}^n} \lambda(x) \\
v_i^*(x) = \lambda(x) -  \widehat{\pi}_i(x_{{\ne}_{i}^+}) \, \sum_{x_i\in \mathcal{L}} \widehat{\pi}_i(x_{{\ne}_{i}^+}) \, \lambda(x), \,  \forall x \in \mathcal{L}^n, \, \forall i \in \mathcal{V}
\end{cases}
\nonumber
\end{equation}
equivalently
\begin{equation}
\begin{array}{lll}
\begin{cases}
v_i^* \equiv k_i \, \lambda, \, \forall i \in \mathcal{V} \\
a^* =\sum_{x \in \mathcal{L}^n} \lambda(x)
\end{cases}
\end{array}
\nonumber
\end{equation}
in a way which enables to write $f$ as
\begin{equation}
\begin{array}{lll}
f &\equiv& \big(\sum_{x \in \mathcal{L}^n}\lambda(x)\big) \, \oneindic + \sum_{i \in \mathcal{V}} k_i^2 \, \lambda \\
   &\equiv& \big(\sum_{x \in \mathcal{L}^n}\lambda(x)\big) \, \oneindic + \sum_{i \in \mathcal{V}} k_i \, \lambda  \\
   &\equiv& \mathscr{K} \, \lambda \nonumber
\end{array}
\end{equation}
Moreover, by Proposition \ref{prop:ki_is_ortho_project}, $k_i$ is an orthogonal projection, hence, $k_i^2 = k_i$, $\forall i \in \mathcal{V}$,
thus one can write
\begin{equation}
\begin{array}{lll}
f  &\equiv& \big(\sum_{x \in \mathcal{L}^n}\lambda(x)\big) \, \oneindic + \sum_{i \in \mathcal{V}} k_i \, \lambda  \\
   &\equiv& \mathscr{K} \, \lambda \nonumber
\end{array}
\end{equation}
Finally, since $f$ is assumed to be arbitrary, one concludes that $\forall f \in \mathscr{F}$, there exists a function $\lambda \in \mathscr{F}$, such that, one can write
$f \equiv \mathscr{K} \, \lambda$. Next, assume the following linear identity:
\begin{equation}
\begin{array}{lll}
\mathscr{K} \, \lambda \equiv \zeroindic
\end{array}
\label{id_Klambda_zero}
\end{equation}
and let us show that its unique solution is given by $\lambda \equiv \zeroindic$. First, one can rewrite linear identity (\ref{id_Klambda_zero}) as
\begin{equation}
\begin{array}{lll}
\big(\sum_{x \in \mathcal{L}^n}\lambda(x)\big) \, \oneindic + \sum_{i \in \mathcal{V}} k_i \, \lambda  \equiv \zeroindic
\end{array}
\nonumber
\end{equation}
Thus, by convolving both hand sides of the latter identity with function $p^* \equiv \mathbb{P}\big(\mathcal{X} = \cdot\big)$, one finds
\begin{equation}
\begin{array}{lll}
0 &=&\big(\sum_{x \in \mathcal{L}^n}\lambda(x)\big) \, \big(\sum_{x \in \mathcal{L}^n} p^*(x)\big) + \sum_{i \in \mathcal{V}} \langle k_i \,\lambda, p^*\rangle  \\
&=&\big(\sum_{x \in \mathcal{L}^n}\lambda(x)\big) \, \big(\sum_{x \in \mathcal{L}^n} p^*(x)\big) + \sum_{i \in \mathcal{V}} \langle k_i \,p^*, \lambda \rangle
\end{array}
\end{equation}
Since one has $k_i \,p^* \equiv \zeroindic$, $\forall i \in \mathcal{V}$, one finally derives $ \big(\sum_{x \in \mathcal{L}^n}\lambda(x)\big) =0 $,
in such a way that, linear identity (\ref{id_Klambda_zero}) is equivalent to the following system:
\begin{equation}
\begin{cases}
\sum_{i \in \mathcal{V}} k_i \, \lambda  \equiv \zeroindic \\
\sum_{x \in \mathcal{L}^n}\lambda(x) = 0
\end{cases}
\nonumber
\end{equation}
Then, since $\forall i \in \mathcal{V}$, $k_i$ is an orthogonal projection, hence, a positive semi-definite linear operator,
one can rewrite the latter system as
\begin{equation}
\begin{cases}
k_i \, \lambda  \equiv \zeroindic, \, \forall i \in \mathcal{V} \\
\sum_{x \in \mathcal{L}^n}\lambda(x) = 0
\end{cases}
\nonumber
\end{equation}
which, in turn, can be reexpressed as
\begin{equation}
\begin{cases}
\lambda(x) - \widehat{\pi}_i(x_{{\ne}_{i}^+}) \, \sum_{x_i\in \mathcal{L}} \widehat{\pi}_i(x_{{\ne}_{i}^+}) \,  \lambda(x) = 0, \, \forall x \in \mathcal{L}^n, \, \forall i \in \mathcal{V} \\
\sum_{x \in \mathcal{L}^n} \lambda(x) = 0
\end{cases}
\nonumber
\end{equation}
and of which, by Lemma \ref{lemma:mrf_syst_nondegenerate}, $\lambda \equiv \zeroindic$ is the unique solution. One deduces that the unique
solution of linear identity (\ref{id_Klambda_zero}) is, in turn, given by $\lambda \equiv \zeroindic$, thus implying that $\mathscr{K}$ is
an invertible linear operator which ends the proof of Lemma \ref{lemma:markov_lin_oper}.
\end{proof}




\subsection{Proof of Lemma \ref{lemma:kq_eq_kbar_q}}
\label{app:proof_lemma_kq_eq_kbar_q}
\begin{proof}
Suppose $q \in \mathscr{F}$. One writes
\begin{equation}
\begin{array}{lll}
\mathscr{K} \, \big(\idoper - \hat{\mathscr{K}}^{\dag} \, \hat{\mathscr{K}} \big) \,  {q} &\equiv&
 \hat{\mathscr{K}} \, \big(\idoper - \hat{\mathscr{K}}^{\dag} \, \hat{\mathscr{K}} \big) \,  {q} +
\widebar{\mathscr{K}} \, \big(\idoper - \hat{\mathscr{K}}^{\dag} \, \hat{\mathscr{K}} \big) \,  {q}
\end{array}
\nonumber
\end{equation}
Observe that
\begin{equation}
\begin{array}{lll}
\hat{\mathscr{K}} \, \big(\idoper - \hat{\mathscr{K}}^{\dag} \, \hat{\mathscr{K}} \big) \,  {q} &\equiv&
\Big(\hat{\mathscr{K}} \,  - \hat{\mathscr{K}} \, \big(\hat{\mathscr{K}}^{\dag}\, \hat{\mathscr{K}} \big) \Big) \,  {q} \\
&\equiv& \big(\hat{\mathscr{K}} \,  - \hat{\mathscr{K}} \big) \,  {q} \\
&\equiv& \zeroindic \\
\end{array}
\nonumber
\end{equation}
implying that
\begin{equation}
\begin{array}{lll}
\mathscr{K} \, \big(\idoper - \hat{\mathscr{K}}^{\dag} \, \hat{\mathscr{K}} \big) \,  {q} &\equiv& \widebar{\mathscr{K}} \, \big(\idoper - \hat{\mathscr{K}}^{\dag} \, \hat{\mathscr{K}} \big) \,  {q}
\end{array}
\nonumber
\end{equation}
then follows the proof of Lemma \ref{lemma:kq_eq_kbar_q} as $q$ is assumed to be arbitrary.
\end{proof}


\subsection{Proof of Lemma \ref{lemma:i_kk_osq_eq_func_kkosq}}
\label{app:proof_lemma_i_kk_osq_eq_func_kkosq}

\begin{proof}
Suppose $q \in \mathscr{F}$. To show the first statement of Lemma\ref{lemma:i_kk_osq_eq_func_kkosq}, one writes
\begin{equation}
\begin{array}{lll}
\mathscr{K}^{-1} \, \widebar{\mathscr{K}} \, \hat{\mathscr{K}}^{\dag} \, \hat{\mathscr{K}} \,  \big({\mathscr{O}}_s \, q \big) &\equiv&
\mathscr{K}^{-1} \, \big({\mathscr{K}} - \hat{\mathscr{K}} \big) \, \hat{\mathscr{K}}^{\dag} \, \hat{\mathscr{K}} \,  \big({\mathscr{O}}_s \, q \big) \\
 &\equiv&
\mathscr{K}^{-1} \, {\mathscr{K}} \, \hat{\mathscr{K}}^{\dag} \, \hat{\mathscr{K}} \,  \big({\mathscr{O}}_s \, q \big)
-\mathscr{K}^{-1} \, \hat{\mathscr{K}} \, \hat{\mathscr{K}}^{\dag} \, \hat{\mathscr{K}} \,  \big({\mathscr{O}}_s \, q \big)
\\
&\equiv&
\hat{\mathscr{K}}^{\dag} \, \hat{\mathscr{K}} \,  \big({\mathscr{O}}_s \, q \big)
-
\mathscr{K}^{-1} \, \hat{\mathscr{K}}  \big({\mathscr{O}}_s \, q \big)
\label{system_i_kk_osq_eq_func_kkosq}
\end{array}
\end{equation}
Furthermore, by Proposition \ref{prop:ks_eq_kos}, one has
\begin{equation}
\hat{\mathscr{K}} \,  \big({\mathscr{O}}_s \, q \big) \equiv {\mathscr{K}} \big({\mathscr{O}}_s \, q \big)
\nonumber
\end{equation}
Thus, by substituting $\hat{\mathscr{K}} \,  \big({\mathscr{O}}_s \, q \big)$ with ${\mathscr{K}} \big({\mathscr{O}}_s \, q \big)$ in
the bottommost formula above, one finds
\begin{equation}
\begin{array}{lll}
\mathscr{K}^{-1} \, \widebar{\mathscr{K}} \, \hat{\mathscr{K}}^{\dag} \, \hat{\mathscr{K}} \,  \big({\mathscr{O}}_s \, q \big) &\equiv&
\hat{\mathscr{K}}^{\dag} \, \hat{\mathscr{K}} \,  \big({\mathscr{O}}_s \, q \big)
-
\mathscr{K}^{-1} \, {\mathscr{K}}  \big({\mathscr{O}}_s \, q \big) \\
&\equiv&
\hat{\mathscr{K}}^{\dag} \, \hat{\mathscr{K}} \,  \big({\mathscr{O}}_s \, q \big)
-
\big({\mathscr{O}}_s \, q \big)\\
&\equiv&
-\big(\idoper - \hat{\mathscr{K}}^{\dag} \, \hat{\mathscr{K}} \big) \,   \big({\mathscr{O}}_s \, q \big)
\nonumber
\end{array}
\end{equation}
The second statement of Lemma \ref{lemma:i_kk_osq_eq_func_kkosq} follows immediately by pre-multiplying both hand sides of the first identity by
$\hat{\mathscr{K}}$ and using the operator equality $\hat{\mathscr{K}} = \hat{\mathscr{K}} \, \big( \hat{\mathscr{K}}^{\dag} \, \hat{\mathscr{K}} \big)$.
This concludes the proof of Lemma \ref{lemma:i_kk_osq_eq_func_kkosq} as $q$ is assumed to be arbitrary.
\end{proof}

\subsection{Proof of Lemma \ref{lemma:local_consistency}}
\label{app:proof_lemma_local_consistency}

\begin{proof}
Suppose $c_1, \dots, c_m$ are $m$ distinct hypervertices of graph $\mathcal{G}$ for some natural number greater than, or equal to $2$,
and assume $m$ local functions of the form $q_i: \mathcal{L}^{|c_i|} \rightarrow \mathbb{R}$, $\forall i\in \big\{1,\dots, m\big\}$
satisfying
\begin{equation}
\begin{array}{lll}
\sum_{k \in \mathcal{V}/c_i} \sum_{x_k \in \mathcal{L}} q_i(x_{c_i}) = \sum_{k \in \mathcal{V}/c_j} \sum_{x_k \in \mathcal{L}} q_j(x_{c_j}), \, \forall x_{c_i \cap c_j} \in \mathcal{L}^{|c_i \cap c_j|},
\forall i,j \in \big\{1,\dots, m\big\}, \text{s.t.,} i \neq j
\end{array}
\label{local_consist_lemma_app}
\end{equation}
and let us show that there exists a global function ${q}: \mathcal{L}^{n} \rightarrow \mathbb{R}$ satisfying
\begin{equation}
\begin{array}{lll}
q_{c_i}(x_{c_i}) = \sum_{k \in \mathcal{V}/c_i} \sum_{x_k \in \mathcal{L}} f(x) = q_i(x_{c_i}), \, \forall x_{c_i} \in \mathcal{L}^{|c_i|}, \forall i \in \big\{1,\dots, m\big\}
\end{array}
\nonumber
\end{equation}
Thus, consider the function series ${f}^{(i)}: \mathcal{L}^{n} \rightarrow \mathbb{R}$, $\forall i\in \big\{1,\dots, m\big\}$ defined as
\begin{equation}
\begin{cases}
{f}^{(1)}(x) = \frac{q_1(x_{c_1})}{L^{n-|c_1|}}, \, \forall x \in \mathcal{L}^n \\
{f}^{(i+1)}(x) = {f}^{(i)}(x) + \frac{q_{i+1}(x_{c_{i+1}})}{L^{n-|c_{i+1}|}} - \frac{{f}_{c_{i+1}}^{(i)}(x_{c_{i+1}})}{L^{n-|c_{i+1}|}}, \, \forall x \in \mathcal{L}^n, \, \forall i\in \big\{1,\dots, m-1\big\}
\end{cases}
\label{function_series}
\end{equation}
and let us show that the following identities:
\begin{equation}
\begin{array}{lll}
{f}_{c_j}^{(i)}(x_{c_j}) = q_j(x_{c_j}), \, \forall x_{c_j} \in \mathcal{L}^{|c_j|}, \, \forall j\in \big\{1,\dots, i\big\}, \, \forall i\in \big\{1,\dots, m\big\}
\end{array}
\label{induction_formula}
\end{equation}
hold, where ${f}_{c_j}^{(i)}$ denotes the margin of ${f}^{(i)}$ with respect to $c_j$, $\forall j\in \big\{1,\dots, i\big\}$, $\forall i\in \big\{1,\dots, m\big\}$.
We will show formula (\ref{induction_formula}) by induction. Thus, first, suppose $i=1$. By the topmost identity in formula (\ref{function_series}), one finds
\begin{equation}
\begin{array}{lll}
{f}_{c_1}^{(1)}(x_{c_1}) &=& \sum_{k \in \mathcal{V}/c_1} \sum_{x_k \in \mathcal{L}} \frac{q_1(x_{c_1})}{L^{n-|c_1|}} \\
                          &=& q_1(x_{c_1}), \, \forall x_{c_1} \in \mathcal{L}^{|c_1|}
\end{array}
\nonumber
\end{equation}
thus showing formula (\ref{induction_formula}) for order $i=1$.
Next, suppose an arbitrary $i \in \big\{1,\dots, m-1\big\}$, and assume that the following identities (up to order $i$):
\begin{equation}
\begin{array}{lll}
{f}_{c_j}^{(i)}(x_{c_j}) = q_j(x_{c_j}), \, \forall x_{c_j} \in \mathcal{L}^{|c_j|}, \, \forall j\in \big\{1,\dots, i\big\}
\end{array}
\label{fomul_recursi_a}
\end{equation}
hold, and let us show that the following identities (up to order $i+1$):
\begin{equation}
\begin{array}{lll}
{f}_{c_j}^{(i+1)}(x_{c_j}) = q_j(x_{c_j}), \, \forall x_{c_j} \in \mathcal{L}^{|c_j|}, \, \forall j\in \big\{1,\dots, i+1\big\}
\end{array}
\label{fomul_recursi_b}
\end{equation}
hold too. Thus, first, suppose $j\in \big\{1,\dots, i\big\}$. By marginalization of both hand sides of the bottommost identity in
formula (\ref{function_series}) with respect to $c_j$, one finds
\begin{equation}
\begin{array}{lll}
{f}_{c_j}^{(i+1)}(x_{c_j}) = {f}_{c_j}^{(i)}(x_{c_j}) + \frac{ \sum_{k \in \mathcal{V}/c_j} \sum_{x_k \in \mathcal{L}}  q_{i+1}(x_{c_{i+1}})}{L^{n-|c_{i+1}|}} -
 \frac{ \sum_{k \in \mathcal{V}/c_j} \sum_{x_k \in \mathcal{L}} {f}_{c_{i+1}}^{(i)}(x_{c_{i+1}})}{L^{n-|c_{i+1}|}}
\end{array}
\label{id_ficj_1}
\end{equation}
One writes
\begin{equation}
\begin{array}{lll}
\sum_{k \in \mathcal{V}/c_j} \sum_{x_k \in \mathcal{L}} {f}_{c_{i+1}}^{(i)}(x_{c_{i+1}}) &=&
\sum_{k \in \mathcal{V}/c_j} \sum_{x_k \in \mathcal{L}} \big( \sum_{t \in \mathcal{V}/c_{i+1}} \sum_{x_t \in \mathcal{L}} f^{(i)}(x) \big) \\
&=& \sum_{t \in \mathcal{V}/c_{i+1}} \sum_{x_t \in \mathcal{L}} \big( \sum_{k \in \mathcal{V}/c_j} \sum_{x_k \in \mathcal{L}}  f^{(i)}(x) \big) \\
&=& \sum_{t \in \mathcal{V}/c_{i+1}} \sum_{x_t \in \mathcal{L}} f_{c_j}^{(i)}(x_{c_j}) \\
&=& \sum_{k \in \mathcal{V}/c_{i+1}} \sum_{x_k \in \mathcal{L}} q_{j}(x_{c_j}), \, \forall x_{c_j} \in \mathcal{L}^{|c_j|}
\end{array}
\nonumber
\end{equation}
as by assumption (formula (\ref{fomul_recursi_a})), one has $f_{c_j}^{(i)}(x_{c_j}) \equiv q_{j}(x_{c_j}), \, \forall x_{c_j} \in \mathcal{L}^{|c_j|}$.
Hence, by substituting  in formula (\ref{id_ficj_1}) above the left hand side expression of the latter formula with its right hand side one, one derives
\begin{equation}
\begin{array}{lll}
{f}_{c_j}^{(i+1)}(x_{c_j})= {f}_{c_j}^{(i)}(x_{c_j}) + \frac{1}{L^{n-|c_{i+1}|}} \, \Big( \sum_{k \in \mathcal{V}/c_j} \sum_{x_k \in \mathcal{L}}  q_{i+1}(x_{c_{i+1}})-
\sum_{k \in \mathcal{V}/c_{i+1}} \sum_{x_k \in \mathcal{L}} q_{j}(x_{c_j})\Big), \, \forall x_{c_j} \in \mathcal{L}^{|c_j|}
\end{array}
\nonumber
\end{equation}
Furthermore, by assumption (formula (\ref{local_consist_lemma_app})), one has
\begin{equation}
\begin{array}{lll}
\sum_{k \in \mathcal{V}/c_j} \sum_{x_k \in \mathcal{L}}  q_{i+1}(x_{c_{i+1}}) - \sum_{k \in \mathcal{V}/c_{i+1}} \sum_{x_k \in \mathcal{L}} q_{j}(x_{c_j}) = 0, \, \forall
 x_{c_{i+1} \cap c_j} \in \mathcal{L}^{|c_{i+1} \cap c_j|}
\end{array}
\nonumber
\end{equation}
One derives
\begin{equation}
\begin{array}{lll}
\end{array}
{f}_{c_{j}}^{(i+1)}(x_{c_{j}}) = f_j(x_{c_j}),  \forall x_j \in \mathcal{L}^{|c_{j}|}, \, \forall j\in \big\{1,\dots, i\big\}
\nonumber
\end{equation}
Thus, it remains to show formula (\ref{fomul_recursi_b}) for $j=i+1$. By marginalization of both hand sides of the the bottomost identity in formula (\ref{function_series})
with respect to $c_{i+1}$, one has
\begin{equation}
\begin{array}{lll}
f_{c_{i+1}}^{(i+1)}(x_{c_{i+1}}) &=& \sum_{k \in \mathcal{V}/c_{i+1}} \sum_{x_k \in \mathcal{L}} {f}^{(i+1)}(x) \\
&=& {f}_{c_{i+1}}^{(i)}(x_{c_{i+1}}) + \frac{ \sum_{k \in \mathcal{V}/c_{i+1}} \sum_{x_k \in \mathcal{L}} q_{i+1}(x_{c_{i+1}})}{L^{n-|c_{i+1}|}}
                         - \frac{\sum_{k \in \mathcal{V}/c_{i+1}} \sum_{x_k \in \mathcal{L}} {f}_{c_{i+1}}^{(i)}(x_{c_{i+1}})}{L^{n-|c_{i+1}|}}\\
                         &=& {f}_{c_{i+1}}^{(i)}(x_{c_{i+1}}) + q_{i+1}(x_{c_{i+1}}) - {f}_{c_{i+1}}^{(i)}(x_{c_{i+1}})\\
                         &=&  q_{i+1}(x_{c_{i+1}}), \,  \forall x_{c_{i+1}} \in \mathcal{L}^{|c_{i+1}|}
\end{array}
\nonumber
\end{equation}
which means that formula (\ref{fomul_recursi_b}) also holds for $j=i+1$, hence the proof of formula (\ref{induction_formula}).
Consequently, by choosing $q$ as $q \equiv f^{(m)}$, one finds by formula (\ref{induction_formula}) assuming $i=m$
\begin{equation}
\begin{array}{lll}
q_{c_i}(x_{c_i}) = q_i(x_{c_i}), \, \forall x_{c_i} \in \mathcal{L}^{|c_i|}, \forall i \in \big\{1,\dots, m\big\}
\end{array}
\nonumber
\end{equation}
where ${q}_{c_i}$ denotes the margin of $q$ with respect to $c_i$, $\forall i\in \big\{1,\dots, m\big\}$.
Therefore, we have constructed a function $q: \mathcal{L}^{n} \rightarrow \mathbb{R}$, such that, $\forall i\in \big\{1,\dots, m\big\}$,
the margin of $q$ with respect to $c_i$ coincides with $q_i$ exactly, hence the proof of Lemma \ref{lemma:local_consistency}.
\end{proof}


\subsection{Proof of Lemma \ref{lemma:equiv_pi_ki_tilde}}
\label{app:proof_lemma_equiv_pi_ki_tilde}

\begin{proof}
The result of Lemma \ref{lemma:equiv_pi_ki_tilde} follows immediately from the result of Lemma \ref{lemma:equiv_pi_ki_tilde} by assuming local
functions of the form $u(x_{s\cup {\ne}_{i}^+})$, $ \forall i \in \mathcal{V}$.
\end{proof}


\subsection{Proof of Lemma \ref{lemma:rectang_reduction_spatial}}
\label{app:proof_lemma_rectang_reduction_spatial}

\begin{proof}
Let us first show that linear identity (\ref{mrf_dist_solu_lin_sys_polyn}) is equivalent to either of following two systems of linear equations:
\begin{equation}
\begin{cases}
{q}_{s\cup {\ne}_{i}^+}(x_{s\cup {\ne}_{i}^+}) - \widehat{\pi}_i(x_{{\ne}_{i}^+}) \, \sum_{x_i\in \mathcal{L}} \widehat{\pi}_i(x_{{\ne}_{i}^+}) \,  {q}_{s\cup {\ne}_{i}^+}(x_{s\cup {\ne}_{i}^+}) = 0, \,
\forall x_{s\cup {\ne}_{i}^+} \in \mathcal{L}^{|s\cup {\ne}_{i}^+|}, \, \forall i \in \mathcal{V} \\
{q}_{s\cup {\ne}_{i}^+}(x_{s\cup {\ne}_{i}^+}) = \sum_{j \in \mathcal{V}/{s\cup {\ne}_{i}^+}} \sum_{x_j \in \mathcal{L}} q(x_1,\dots,x_j,\dots,x_n), \forall x_{s\cup {\ne}_{i}^+}
\in \mathcal{L}^{|s\cup {\ne}_{i}^+|},  \forall i \in \mathcal{V} \\
\sum_{x \in \mathcal{L}^n} q(x) = 1
\end{cases}
\label{rectang_reduction_spatial_b_hat}
\end{equation}
or
\begin{equation}
\begin{cases}
{q}_{s\cup {\ne}_{i}^+}(x_{s\cup {\ne}_{i}^+}) - {\pi}_i(x_{{\ne}_{i}^+}) \, \sum_{x_i \in \mathcal{L}} {q}_{s\cup {\ne}_{i}^+}(x_{s\cup {\ne}_{i}^+}) = 0, \,
\forall x_{s\cup {\ne}_{i}^+} \in \mathcal{L}^{|s\cup {\ne}_{i}^+|}, \, \forall i \in \mathcal{V} \\
{q}_{s\cup {\ne}_{i}^+}(x_{s\cup {\ne}_{i}^+}) = \sum_{j \in \mathcal{V}/{s\cup {\ne}_{i}^+}} \sum_{x_j \in \mathcal{L}} q(x_1,\dots,x_j,\dots,x_n), \forall x_{s\cup {\ne}_{i}^+}
\in \mathcal{L}^{|s\cup {\ne}_{i}^+|},  \forall i \in \mathcal{V} \\
\sum_{x \in \mathcal{L}^n} q(x) = 1
\end{cases}
\label{rectang_reduction_spatial_b}
\end{equation}

First of all, rewrite linear identity (\ref{mrf_dist_solu_lin_sys_polyn}) as
\begin{equation}
\begin{array}{lll}
\big( \sum_{x \in \mathcal{L}^n } q(x) \big) \, \oneindic +  \sum_{i \in \mathcal{V}} \widehat{k}_i \, q \equiv \oneindic
\end{array}
\label{corol_systm_a}
\end{equation}
Next, let us show that linear identity (\ref{mrf_dist_solu_lin_sys_polyn}) is equivalent to the following
linear system:
\begin{equation}
\begin{cases}
\sum_{i \in \mathcal{V}} \widehat{k}_i \, q \equiv \zeroindic \\
\sum_{x \in \mathcal{L}^n } q(x) = 1
\end{cases}
\label{corol_systm_b}
\end{equation}
For this purpose, suppose $q$ is a solution to linear identity (\ref{mrf_dist_solu_lin_sys_polyn}).
Hence, by convolving both hand sides of linear identity (\ref{mrf_dist_solu_lin_sys_polyn}) with
$\mathbb{P}\big(\mathcal{X} = \cdot \big)$, and since one has
\begin{equation}
\widehat{k}_i \, \mathbb{P}\big(\mathcal{X} = \cdot \big) \equiv \zeroindic, \, \forall i \in \mathcal{V}
\nonumber
\end{equation}
one finds
\begin{equation}
\begin{array}{lll}
\sum_{x \in \mathcal{L}^n} q(x) &=& \sum_{x \in \mathcal{L}^n} \mathbb{P}\big(\mathcal{X} = x \big) \\
 &=& 1
\end{array}
\nonumber
\end{equation}
implying that $q$ is a solution of linear system (\ref{corol_systm_b}).

Conversely, suppose $q$ is a solution of linear system (\ref{corol_systm_b}), hence a solution to the following linear system:
\begin{equation}
\begin{cases}
\sum_{i \in \mathcal{V}} \widehat{k}_i \, q \equiv \zeroindic \\
\sum_{x \in \mathcal{L}^n } q(x) \, \oneindic  =  \oneindic
\end{cases}
\label{corol_systm_c}
\end{equation}
too. By adding both identities of linear system (\ref{corol_systm_c}) sidewise, one finds that $q$ satisfies linear identity (\ref{mrf_dist_solu_lin_sys_polyn}).
Since $q$ is assumed to be arbitrary, one concludes that linear identity (\ref{mrf_dist_solu_lin_sys_polyn}) and linear system (\ref{corol_systm_c}) are equivalent.

Now, suppose suppose $q$ is a solution of linear (\ref{corol_systm_b}). It follows that, by convolving both hand sides of the topmost identity of linear system
(\ref{corol_systm_b}) with $q$, one finds
\begin{equation}
\sum_{i \in \mathcal{V}} \langle q, \widehat{k}_i \, q \rangle \equiv \zeroindic
\end{equation}
But, since by Proposition \ref{prop:kilde_i_is_ortho_project}, $\forall i \in \mathcal{V}$, $\widehat{k}_i^{(s)}$ is an orthogonal projection,
hence, a positive semi-definite linear operator, one finds that the topmost identity of system (\ref{corol_systm_b}) is equivalent to the
following system of linear identities:
\begin{equation}
\begin{array}{lll}
\widehat{k}_i \, q \equiv \zeroindic, \, \forall i \in \mathcal{V}
\end{array}
\nonumber
\end{equation}
in such a way that, system (\ref{corol_systm_b}), equivalently, linear identity (\ref{corol_systm_a}) is equivalent to the following system:
\begin{equation}
\begin{cases}
\widehat{k}_i \, q \equiv \zeroindic, \, \forall i \in \mathcal{V}  \\
\sum_{x \in \mathcal{L}^n } q(x) = 1
\end{cases}
\nonumber
\end{equation}
and the latter is nothing else than linear system (\ref{rectang_reduction_spatial_b_hat}).

One concludes that system (\ref{corol_systm_c}) can rewrite as
\begin{equation}
\begin{cases}
{q}_{s\cup {\ne}_{i}^+}(x_{s\cup {\ne}_{i}^+}) - {\pi}_i(x_{{\ne}_{i}^+}) \, \sum_{x_i \in \mathcal{L}} {q}_{s\cup {\ne}_{i}^+}(x_{s\cup {\ne}_{i}^+}) = 0, \,
\forall x_{s\cup {\ne}_{i}^+} \in \mathcal{L}^{|s\cup {\ne}_{i}^+|}, \, \forall i \in \mathcal{V} \\
{q}_{s\cup {\ne}_{i}^+}(x_{s\cup {\ne}_{i}^+}) = \sum_{j \in \mathcal{V}/{s\cup {\ne}_{i}^+}} \sum_{x_j \in \mathcal{L}} q(x_1,\dots,x_j,\dots,x_n), \forall x_{s\cup {\ne}_{i}^+}
\in \mathcal{L}^{|s\cup {\ne}_{i}^+|},  \forall i \in \mathcal{V} \\
\sum_{x \in \mathcal{L}^n} q(x) = 1
\end{cases}
\nonumber 
\end{equation}
thus establishing equivalence between both systems of linear equations (\ref{mrf_dist_solu_lin_sys_polyn}) and
(\ref{rectang_reduction_spatial_b_hat}). Last, by Lemma \ref{lemma:equiv_pi_ki_tilde}, one finds that
both systems of linear equations (\ref{rectang_reduction_spatial_b_hat}) and (\ref{rectang_reduction_spatial_b})
are, in turn, equivalent, thus concluding the proof of Lemma \ref{lemma:rectang_reduction_spatial}.
\end{proof}

\end{appendices}

\bibliographystyle{ACM-Reference-Format}
\bibliography{ibechar_mybibfile_jai_two}
\end{document}